\documentclass[leqno]{article} 
\usepackage{amsmath,amsthm}
\usepackage{amssymb,amscd,latexsym}
\usepackage{boxedminipage}
\usepackage{eepic}
\usepackage{epic}
\usepackage{wrapfig}
\usepackage{graphicx}

\theoremstyle{plain}

\newtheorem*{def-theo}{Definition-Theorem}

\theoremstyle{definition}

\newtheorem*{definition}{Definition}

\theoremstyle{remark}

\newtheorem{thm}{Theorem}[section]

\newtheorem{lem}[thm]{Lemma}

\theoremstyle{definition}

\newtheorem{defn}[thm]{Definition}

\newtheorem{rem}[thm]{Remark}

\numberwithin{equation}{section}
\newcommand{\bp}{\begin{pmatrix}}
\newcommand{\ep}{\end{pmatrix}}
\newcommand{\bps}{\begin{smallmatrix}}
\newcommand{\eps}{\end{smallmatrix}}
\def\C{{\mathbb C}}

\def\R{{\mathbb R}}
\def\Z{{\mathbb Z}}

\def\O{{\mathcal O}}

\def \0{{\bf 0}}
\def \1{{\bf 1}}

\def \mf#1#2#3#4{
\xymatrix{{#1}\  \ar@<0.4ex>[r]^{{#2}} & \ {#4}
\ar@<0.4ex>[l]^{{#3}}
}
}

\def \mfs#1#2#3#4{\!
\xymatrix@C=1.5em{{#1} \! \ar@<0.2ex>[r]^{{#2}} & \! {#4}
\ar@<0.2ex>[l]^{{#3}}
}
\!}

\def \mfl#1#2#3#4{
\xymatrix@C=2.6em{{#1}\  \ar@<0.4ex>[r]^{{#2}} &\  {#4}
\ar@<0.2ex>[l]^{{#3}}
}
}

\def \mfss#1#2#3#4{\!
\xymatrix@C=1.5em{{#1} \ar@<0.3ex>[r]^{{#2}} & {#4}
\ar@<0.3ex>[l]^{{#3}}
}
\!}

\begin{document}

\title{{\Large Monoids in the fundamental groups of the complement of 
logarithmic free divisors in $\C^3$
\footnote{
The present paper is a complete version with proofs of [S-I1].}
}
}

\author{Kyoji Saito and Tadashi Ishibe}
\date{}
\maketitle
\vspace{ -0.6cm}

 {\renewcommand{\baselinestretch}{0.1}

\vspace{ -0.4cm}
\tableofcontents

\begin{abstract}
\noindent
We study monoids generated by Zariski-van Kampen generators in the 17 fundamental groups of the complement of logarithmic free divisors in $\C^3$ listed by Sekiguchi (Theorem~1).
Five of them are 
Artin monoids 
and eight of them are free abelian monoids. 
{\it The remaining four monoids are not Gau\ss ian and, hence, are neither Garside nor  Artin} (Theorem~2). 
However, {\it we introduce, similarly to Artin monoids, fundamental elements and show their existence {\rm (Theorem 3)}}.
One of the four non-Gaussian monoids satisfies the cancellation condition
{\rm (Theorem 4)}. 
\end{abstract}

 \section{Introduction}

A hypersurface $D$ in $\C^l$ ($l\!\in\!\Z_{\ge0}$) 
is called a {\it logarithmic free divisor} (\cite{[S1]}), if the associated module $Der_{\C^l}(-{\rm log}(D))$ 
of logarithmic vector fields is a free $\O_{\C^l}$-module.
Classical example of logarithmic free divisors is the discriminant loci of 
a finite reflection group ([S1,2,3,4]). 
The fundamental group of the complement of the discriminant loci 
is presented (Brieskorn [B]) by certain positive homogeneous relations, called Artin braid relations.\
The group (resp.\ monoid) defined by that presentation is called an 
{\it Artin group} (resp.\ {\it Artin monoid}) of finite type \cite{[B-S]}, for which the word problem and other problems are solved using a particular element $\Delta$, 
the {\it fundamental elements}, in the monoids
 (\cite{[B-S]},\cite{[D]},\cite{[G]}).

In [Se1], Sekiguchi listed up 17 weighted homogeneous  polynomials,  
defining logarithmic free divisors in $\C^3$, whose weights coincide with those of the discriminant of 
types $\mathrm{A_3}$, $\mathrm{B_3}$ or $\mathrm{H_3}$.
Then, the fundamental groups of the complements of the divisors are presented by Zariski-van Kampen method by \cite{[I]} (we recall the result in \S3). It turns out that the defining relations can be reformulated by a system of positive homogeneous relations in the sense explained in \S4 of the present paper, so that we can introduce monoids defined by them. We show that, among 17 monoids, 5 are Artin monoids, and 8 are free abelian monoids. However, four remaining monoids are not Gaussian, and hence are neither Garside nor Artin (\S5). Nevertheless, we show that they carry certain particular elements similar to the fundamental elements in Artin monoids (\S6). 

Let us explain more details of the contents.
The 17 Sekiguchi-polynomials $\Delta_X(x,y,z)$  are labeled by the type $X\in \{\mathrm{A_i,A_{ii},B_i,B_{ii},B_{iii},B_{iv},B_v,B_{vi},}$ $\mathrm{B_{vii},H_i,}$ $\mathrm{H_{ii},H_{iii},}$ $\mathrm{H_{iv},H_v},$ $\mathrm{H_{vi},H_{vii},H_{viii}}\}$ (\S2). They are monic polynomials of degree 3 in the variable $z$. We calculate the fundamental group of the complement of the divisor $D_X:=\{\Delta_X(x,y,z)=0\}$ in $\C^3$  by choosing Zariski-pencils $l$ 
in $z$-coordinate direction, which intersect with the divisor $D_X$ by 3 points. Zariski-van Kampen method gives a presentation of the fundamental group $\pi_1(\C^3\setminus D_X,*)$ with respect to three generators $a$, $b$ and $c$ presented by a choice of paths in the pencil turning once around each of three intersection points.
%
%

We rewrite the Zariski-van Kampen relations into a system of positive homogeneous relations (not unique, \S4 Theorem 1), and study the group $G_X$ and the monoid $M_X$ defined by the relations as well as the localization homomorphism $M_X\!\to\! G_X$, where $G_X$ is naturally isomorphic to $\pi_1(\C^3\setminus D_X,*)$. 
We denote by $G_X^+$ the image of $M_X$ in $G_X$, that is, the monoid generated by the Zariski-van Kampen generators $\{a,b,c\}$ in $\pi_1(\C^3\setminus D_X,*)$. The $G_X^+$ depends  on the choice of generators but not on homogeneous relations, whereas the monoid $M_X$ does. 
It turns out that $M_X$ are Artin monoids for the types ${\mathrm{A_i}},{\mathrm{B_i}},{\mathrm{H_i}},{\mathrm{A_{ii}}},\mathrm{B_{iv}}$, and are free abelian monoid  for  
the types ${\mathrm{B_v}}, {\mathrm{B_{vii}}}, {\mathrm{H_{iv}}}, {\mathrm{H_{v}}},$ ${\mathrm{H_{vi}}}, {\mathrm{H_{vii}}}, {\mathrm{H_{viii}}},{\mathrm{B_{iii}}}$ so that one has natural isomorphisms: $M_X\simeq G_X^+$.
%
However, for any of the remaining four  types $\mathrm{B_{ii}}, \mathrm{B_{vi}}, \mathrm{H_{ii}}, \mathrm{H_{iii}}$, 
the  monoids $G_X^+$ does not admit the divisibility theory 
(see \cite[\S5]{[B-S]}, or \S5 Theorem 2 of present paper). 
That is, they are not Gaussian groups \cite[\S2]{[D-P]}, and, hence, they are neither Artin nor Garside groups (actually, we have an isomorphism $M_{\mathrm{B_{vi}}}\simeq M_{\mathrm{H_{iii}}}$ and hence $G^+_{\mathrm{B_{vi}}}\simeq G^+_{\mathrm{H_{iii}}}$).

On the other hand, as one  main result of the present paper, we show that the monoid $M_X$ carries some distinguished elements, 
which we call {\it fundamental} (\S6 Theorem 3).
Namely, we call an element $\Delta\!\in\! M_X$ fundamental if
there exists a permutation $\sigma_\Delta$ of the set 
$\{a,b,c\}/\sim$ (see \S6 )
such that for any $d\in \{a,b,c\}/\!\sim$, there exists 
$\Delta_d\!\in\! M_X$ such that the following relation holds:
\vspace{-0.2cm}
\[
\Delta\ =\ d\cdot\Delta_d\ =\ \Delta_d \cdot \sigma_\Delta(d).
\vspace{-0.2cm}\] 
The set $\mathcal{F}(M_X)$ of fundamental elements in $M_X$ form a submonoid of $M_X$ such that $\mathcal{QZ}(M_X)\mathcal{F}(M_X)\!=\!\mathcal{F}(M_X)\mathcal{QZ}(M_X)\!=\!\mathcal{F}(M_X)$ (see \S6 {\bf Fact 3.}) where $\mathcal{QZ}(M_X)$ is the quasi-center of $M_X$.\footnote
{An element $\Delta\!\in\!\! M_X$ is called quasi-central (\cite[\!7.1]{[B-S]})
if $d\!\cdot\!\Delta\!=\! \Delta\!\cdot\! \sigma_\Delta(d)$ 
for $d\in\{a,b,c\}$.\!\!}
For an Artin monoid of finite type, $\mathcal{F}(M_X)$ is generated by a single element $\Delta$ and $\mathcal{F}(M_X)\!=\!\Delta^{\Z_{\ge 1}}$ ([B-S]). Since the 

\newpage
\noindent
localization morphism induces a map
$\mathcal{F}(M_X)\!\to\!\mathcal{F}(G_X^+)$, the fact $\mathcal{F}(M_X)\!\not=\!\emptyset$ for all 17 monoids (\S6 Theorem3) implies $\mathcal{F}(G_X^+)\!\not=\!\emptyset$. We ask, more generally, {\it whether the monoid generated by Zariski-van Kampen generators in the local fundamental group of the complement of a free divisor 
 has always a fundamental element} (see \S6 Remark 6.4). In the 4 types $\mathrm{B_{ii}},{\mathrm{B_{vi}}},{\mathrm{H_{ii}}},{\mathrm{H_{iii}}}$, we observe that $\mathcal{F}(G_X^+)$ is not singly generated. Therefore, we ask, also, {\it whether the set of fundamental elements $\mathcal{F}(G_X^+)$ is finitely generated over $\mathcal{QZ}(G_X^+)$ or not}.

In \S7, we discuss about the cancellation condition on the monoid $M_X$. In fact, this condition together with the existence of fundamental elements (shown in \S6), imply that the localization morphism  $M_X\! \to\! G_X^+$ is an isomorphism. 
An Artin monoid or a free abelian monoid satisfies already this condition ([B-S]).  We show that {\it the monoid  $M_{\mathrm{B_{ii}}}$ satisfies the cancellation condition} (Theorem4).  
For the remaining three types 
${\mathrm{B_{vi}}},{\mathrm{H_{ii}}},{\mathrm{H_{iii}}}$, 
we do not know whether the localization map is $M_X\! \to\! G_X^+$ is injective or not. That is, we don't know whether we have sufficiently many defining relations to assert the cancellation condition or not.

Finally in \S8, we construct non-abelian representations of the groups 
$G_{\mathrm{B_{ii}}}, G_{\mathrm{B_{vi}}},$ 
$G_{\mathrm{H_{ii}}}$ and $ G_{\mathrm{H_{iii}}}$ 
into $\mathrm{GL_2}(\C)$ (Theorem 5). Actually, this result is independent of \S5, 6 and 7, and is used in the proof of Theorem 2 in \S5.

\section{Sekiguchi's Polynomial}
J. Sekiguchi [Se1,2] listed the following 17 weighted homogeneous polynomials 
$\Delta$ in three variables $(x,y,z)$ 
satisfying freeness criterion by K.Saito \cite{[S1]}. 
\begin{eqnarray}\nonumber
\Delta_{\mathrm{A_i}} (x,y,z)&:=&
-4x^3y^2-27y^4+16x^4z+144xy^2z-128x^2z^2+256z^3 \\ \nonumber
\Delta_{\mathrm{A_{ii}}} (x,y,z)&:=&2x^6-3x^4z+18x^3y^2-18xy^2z+27y^4+z^3 \\ \nonumber
\Delta_{\mathrm{B_i}}(x,y,z)&:=&z(x^2y^2-4y^3-4x^3z+18xyz-27z^2) \\ \nonumber
\Delta_{\mathrm{B_{ii}}}(x,y,z)&:=&z(-2y^3+4x^3z+18xyz+27z^2) \\ \nonumber
\Delta_{\mathrm{B_{iii}}}(x,y,z)&:=&z(-2y^3+9xyz+45z^2) \\ \nonumber
\Delta_{\mathrm{B_{iv}}}(x,y,z)&:=&z(9x^2y^2-4y^3+18xyz+9z^2) \\ \nonumber
\Delta_{\mathrm{B_{v}}}(x,y,z)&:=&xy^4+y^3z+z^3 \\ \nonumber
\Delta_{\mathrm{B_{vi}}}(x,y,z)&:=&9xy^4+6x^2y^2z-4y^3z+x^3z^2-12xyz^2+4z^3 \\ \nonumber
\Delta_{\mathrm{B_{vii}}}(x,y,z)&:=&(1/2)xy^4-2x^2y^2z-y^3z+2x^3z^2+2xyz^2+z^3 \\ \nonumber
\Delta_{\mathrm{H_i}}(x,y,z)&:=&-50z^3+(4x^5-50x^2y)z^2+(4x^7+60x^4y^2+225xy^3)z \\ \nonumber 
& & -(135/2)y^5-115x^3y^4-10x^6y^3-4x^9y^2 \\ \nonumber
\Delta_{\mathrm{H_{ii}}}(x,y,z)&:=&100x^3y^4+y^5+40x^4y^2z-10xy^3z+4x^5z^2-15x^2yz^2+z^3 \\ \nonumber
\Delta_{\mathrm{H_{iii}}}(x,y,z)&:=&8x^3y^4+108y^5-36xy^3z-x^2yz^2+4z^3 \\ \nonumber
\Delta_{\mathrm{H_{iv}}}(x,y,z)&:=&y^5-2xy^3z+x^2yz^2+z^3 \\ \nonumber
\Delta_{\mathrm{H_{v}}}(x,y,z)&:=&x^3y^4-y^5+3xy^3z+z^3 \\ \nonumber
\Delta_{\mathrm{H_{vi}}}(x,y,z)&:=&x^3y^4+y^5-2x^4y^2z-4xy^3z+x^5z^2+3x^2yz^2+z^3\\ \nonumber
\Delta_{\mathrm{H_{vii}}}(x,y,z)&:=&xy^3z+y^5+z^3 \\ \nonumber
\Delta_{\mathrm{H_{viii}}}(x,y,z)&:=&x^3y^4+y^5-8x^4y^2z-7xy^3z+16x^5z^2+12x^2yz^2+z^3. \nonumber
\end{eqnarray}

Here, the polynomials are classified into three types A, B and H 
according as the  numerical data 
$({\rm deg}(x),{\rm deg}(y),{\rm deg}(z);{\rm deg}(\Delta))$ is equal to 
$(2,3,4;12)$, $(2,4,6;18)$ or $(2,6,10;30)$, respectively. In 
each type, the polynomials are numbered by small Roman numerals 
i, ii,\ldots etc. We remark
that, in all cases, the polynomial is a monic polynomial of degree 3 in the variable $z$.

\section{Zariski-van Kampen method}

Let $X$ be one of the 17 types  
$\mathrm{A_i}$,$\mathrm{A_{ii}}$,$\mathrm{B_i}$,\ldots, $\mathrm{B_{vii}}$, 
$\mathrm{H_i}$,\ldots,$\mathrm{H_{viii}}$. 
In the present section, we recall from \cite{[I]} the calculation of  
the fundamental group $\pi_1(S_X \setminus D_X, *_X)$ 
 of the complement of the free divisor $D_X$ in the space $S_X$ by Zarisik-van Kampen method, where we put $S_X:=\C^3$ and 
\begin{equation}
D_X:=\{(x,y,z)\in \C^3\mid \Delta_X(x,y,z)=0\}. 
\end{equation}

The first step is the following  reduction from the space $S_X$ to a plane $H_X$.

\medskip
\begin{lem}[Lefschetz Theorem [H-L{]}]
{\it  Let $H_X\!\subset\! S_X$ be a hyperplane defined by $x\!=\!\varepsilon$ 
for a general $\varepsilon\!\in\! \C^\times$.
Then, the natural inclusion induces an isomorphism:
\begin{equation}
\pi_1(H_X \setminus {(H_X \cap D_X)}, *_X) 
\stackrel{\sim}{\rightarrow} 
\pi_1(S_X \setminus D_X, *_X) 
\end{equation}
 for any choice of a base point $*_X\in H_X \setminus {(H_X \cap D_X)}$.}
\end{lem}
\medskip

The second step is to apply Zariski-van Kampen method, using pencils.

To define pencils,  
we  consider the projection map $\pi$ 
from $S_X$ to the space $T_X:=\C^2$ of coordinates $x,y$ by 
forgetting the coordinate $z$. The fibers of the projection $\pi$ shall be 
called the Zariski-pencils.
The $\pi|D_X$ is a triple covering map, whose branching loci 
(or, bifurcation set) $B_X$ is defined by  
\begin{equation}
B_X:=\{(x,y)\in T_X=\C^2 \mid \omega_X(x,y)=0\}, 
\end{equation}
where $\omega_X(x,y):= \delta \Bigl( \Delta_{X},\frac{\partial {\Delta_{X}}}{\partial z}\Bigl)$ is the resultant of $\Delta_{X}$ 
and $\frac{\partial {\Delta_{X}}}{\partial z}$ with respect to the 
variable $z$. 
In fact, $\omega_X$ is a weighted homogeneous polynomial which is monic 
in the variable $y$. As we can see explicitly from Table 
of the equations below, the restriction of $\omega_X$ to
the line $L_X\!:=\!\{(x,y)\!\in\! T_X \mid x\!=\!\varepsilon\}$, where $\varepsilon\!=\!-1$ 
for the type A and $\varepsilon\!=\!1$ for the types B and H, 
is totally real, i.e.~all roots of the equation 
$\omega_X(\varepsilon,y)=0$ in $y$ are real numbers,
except for $\mathrm{B_{vii}}$ and $\mathrm{H_{vi}}$.

\[
\begin{array}{rllrll}
\omega_{\mathrm{A_i}}(-1,y) &=& -c y^2(27y^2-8)^3, &
\omega_{\mathrm{A_{ii}}}(-1,y) &=& c y^6(27y^2-4),\\
\omega_{\mathrm{B_i}}(1,y) &=& c y^4(1-4y)^2(1-3y)^3, &
\omega_{\mathrm{B_{ii}}}(1,y) &=& c y^6(2+3y)^2(1+3y),\\
\omega_{\mathrm{B_{iii}}}(1,y) &=& c y^8(9+40y), &
\omega_{\mathrm{B_{iv}}}(1,y) &=& c y^7(9-4y)^2,\\
\omega_{\mathrm{B_{v}}}(1,y) &=& c y^8(27+4y), &
\omega_{\mathrm{B_{vi}}}(1,y) &=& c y^5(3-64y)^2(2-y)^3,\\
\omega_{\mathrm{B_{vii}}}(1,y) &=& c y^7(16y^2+13y+8), &
\omega_{\mathrm{H_{i}}}(1,y) &=& c y^2(2-5y)^5(2+27y)^3,\\
\omega_{\mathrm{H_{ii}}}(1,y) &=& c y^5(4-27y)^5(12-y)^4, &
\omega_{\mathrm{H_{iii}}}(1,y) &=& c y^7(1-54y)^3,\\
\omega_{\mathrm{H_{iv}}}(1,y) &=& -c y^9(4+27y), &
\omega_{\mathrm{H_{v}}}(1,y) &=& -c y^8(1+y)^2,\\
\omega_{\mathrm{H_{vi}}}(1,y) &=& -c y^8(27y^2+14y+3),&
\omega_{\mathrm{H_{vii}}}(1,y) &=& c y^9(4+27y),\\
\omega_{\mathrm{H_{viii}}}(1,y) &=& -c y^7(3+y)^2(32+27y).
\end{array}
\]

 \newpage

Remember that $H_X=\pi^{-1}(L_X)$.
We apply, now, Zariski-van Kampen method (see [Ch],[T-S] for instance) 
to calculate the fundamental group 
$\pi_1(H_X\setminus {C_X},*_X)$ of the complement of the plane curve 
$C_X:=H_X \cap D_X$ in the $yz$-plane. 
Let us explain the step wisely the process more in details.


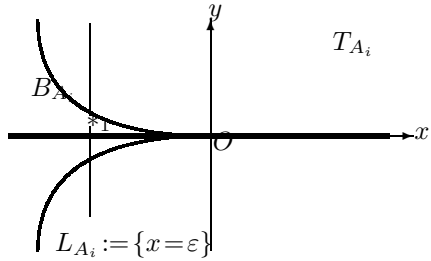
\begin{figure}
\setlength{\unitlength}{0.0060333in}
\begingroup\makeatletter\ifx\SetFigFont\undefined%
\gdef\SetFigFont#1#2#3#4#5{%
  \reset@font\fontsize{#1}{#2pt}%
  \fontfamily{#3}\fontseries{#4}\fontshape{#5}%
  \selectfont}%
\fi\endgroup%
{\renewcommand{\dashlinestretch}{30}
\begin{picture}(100,200)(-200,-0)
\put(-100,100){\line(1,0){350}}
\put(250,100){\vector(1,0){1}}
\put(252,97){$x$}
\put(75,0){\line(0,1){200}}
\put(73,204){$y$}
\put(75,200){\vector(0,1){1}}
\put(77,88){$O$}

\put(181,175){$T_{A_i}$}
\put(-81,135){$B_{A_i}$}
\put(-30,30){\line(0,1){78}}
\put(-30,118){\line(0,1){80}}
\put(-60,-0){$L_{A_i}\!:=\!\{x\!=\!\varepsilon\}$}
\put(-33,108){$*_{1}$}
\thicklines
\qbezier(75,100)(-75,100)(-75,200)
\qbezier(75,100)(-75,100)(-75,0)
\put(-100,100){\line(1,0){330}}
\put(-100,101){\line(1,0){330}}
\put(-100,99){\line(1,0){330}}
\end{picture}
\caption{bifurcation set $B_{A_i}$ in $T_{A_i}$}
\label{fig1}
}
\vspace{-0.3cm}
\end{figure}

\bigskip
1. Choose a base point $*_X$ in $L_X\setminus (L_X\cap B_X)$ and
call the associated pencil $l_{*_1}:=\pi^{-1}(*_1)$ the {\it basic pencil}. 

\smallskip
2. 
Choose and 
fix  (i) the base point $*_X\in l_{*_1}\! \setminus\!(l_{*_1}\!\cap\! D_X)$ and (ii)
three mutually disjoint (except at $*_X$) path connecting $*_X$ with the 
three points $l_{*_1}\!\cap\! D_X$ in the basic pencil. Accordingly, fix three 
the generators, say $a,\ b$ and $c$,  of the free group $F_3\!:=\!\pi_1(l_{*_1}\! \setminus\!(l_{*_1}\!\cap\! D_X),*_X)$ (they are presented by the movements
from $*_X$ 
to close to the points on $D_X$ along paths, then turn once around the end points of the paths
counterclockwise, and then return to $*_X$ along the paths).

\smallskip
3. Move the pencils $l_{t}:=\pi^{-1}(t)$ by moving $t$ along a closed 
path $\gamma$ in $L_X\setminus (L_X\cap B_X)$ turning around a bifurcation point 
in $L_X\cap B_X$. This induces a (braid) action $\gamma_*: F_3\to F_3$, 
and we define the relations: $\gamma_*(a)=a,\gamma_*(b)=b,\gamma_*(c)=c$.
Running $\gamma$ over all generators of $\pi_1(L_x\setminus(L_X\setminus (L_X\cap B_X),*_1)$, we obtain all list of defining relations of the group 
$\pi_1(H_X\setminus {C_X},*_X)$.

\medskip
Actually, in most of the cases except for the cases $X\in\{\mathrm{B_{vii}, H_{ii},H_{vi}},\mathrm{H_{viii}}\}$,
we can find totally real region in $L_X$ in the sense that, if $*_1\in\{$totally real 
region$\}$, three roots $l_{*_1}\cap D_X$ of the equation 
$\Delta_X(\varepsilon,*_1,z)=0$ 
with respect to the coordinate $z$ of the pencil are real numbers.
In such case, we choose the base point $*_X$ and then the  paths $a,b,c$ in 
the basic pencil $l_{*_1}$ as in Figure 2. along three intervals connecting $*_X$ with the points in  $l_{*_1}\cap D_X$. 

\bigskip
\centerline{
\setlength{\unitlength}{0.00039333in}
\begingroup\makeatletter\ifx\SetFigFont\undefined%
\gdef\SetFigFont#1#2#3#4#5{%
  \reset@font\fontsize{#1}{#2pt}%
  \fontfamily{#3}\fontseries{#4}\fontshape{#5}%
  \selectfont}%
\fi\endgroup%
{\renewcommand{\dashlinestretch}{30}
\begin{picture}(3849,2975)(10,100)
\put(1329,663){\ellipse{34}{34}}
\put(2919,667){\ellipse{34}{34}}
\put(483,659){\ellipse{34}{34}}
\put(1646,2881){\ellipse{144}{144}}
\path(1333,1274)(1358,1184)(1403,1259)
\path(1592,2825)(464,758)
\path(1597,2819)(542,753)
\path(673,1177)(658,1092)(733,1132)
\path(689,1066)(749,1116)(759,1041)
\path(1646,2802)(1370,747)
\path(1386,1185)(1436,1245)(1476,1175)
\path(1276,534)(1331,561)(1287,601)
\path(1641,2813)(1639,2812)(1641,2808)(1300,755)
\path(1695,2821)(2945,755)
\path(1687,2821)(2885,749)
\path(2881,538)(2936,565)(2892,605)
\path(2704,1097)(2699,1172)(2774,1132)
\path(2576,1172)(2651,1117)(2661,1192)
\path(1598,2824)(1687,2935)
\path(1603,2933)(1689,2827)
\path(427,536)(482,563)(438,603)
\path(12,659)(3837,659)
\path(3717.000,629.000)(3837.000,659.000)(3717.000,689.000)
\path(1301,758)(1300,758)(1297,756)
	(1290,751)(1279,745)(1266,736)
	(1253,726)(1242,714)(1232,701)
	(1225,685)(1222,666)(1224,650)
	(1228,636)(1233,624)(1237,614)
	(1242,606)(1246,600)(1250,593)
	(1256,587)(1264,580)(1274,573)
	(1287,566)(1303,562)(1318,562)
	(1331,564)(1343,568)(1353,572)
	(1361,576)(1368,581)(1375,585)
	(1381,589)(1387,595)(1394,601)
	(1402,609)(1410,618)(1418,630)
	(1424,643)(1427,659)(1427,674)
	(1423,687)(1417,699)(1410,711)
	(1401,721)(1392,731)(1384,739)
	(1378,745)(1374,749)(1372,751)
\path(2888,743)(2886,742)(2882,740)
	(2875,737)(2866,732)(2857,727)
	(2848,720)(2839,713)(2831,704)
	(2825,694)(2821,681)(2820,667)
	(2823,652)(2827,638)(2832,627)
	(2837,617)(2841,609)(2845,602)
	(2849,595)(2855,588)(2862,581)
	(2872,574)(2885,567)(2901,563)
	(2916,563)(2929,565)(2941,569)
	(2952,573)(2961,577)(2968,581)
	(2975,585)(2982,590)(2989,595)
	(2996,601)(3003,609)(3011,619)
	(3018,630)(3022,644)(3022,664)
	(3016,682)(3006,699)(2993,714)
	(2979,728)(2964,741)(2952,751)
	(2944,758)(2940,761)(2939,761)
\path(464,762)(463,762)(460,760)
	(453,755)(443,749)(431,740)
	(418,730)(407,718)(398,705)
	(391,689)(389,670)(391,654)
	(396,640)(400,628)(405,618)
	(409,610)(413,604)(418,597)
	(423,591)(431,584)(441,577)
	(454,570)(470,566)(484,566)
	(498,568)(510,572)(520,576)
	(528,580)(535,585)(541,589)
	(547,594)(553,599)(560,605)
	(568,613)(576,623)(584,634)
	(591,647)(595,662)(595,677)
	(592,690)(587,702)(580,713)
	(572,724)(564,733)(557,741)
	(551,747)(547,750)(545,752)
\put(1723,3099){\makebox(0,0)[lb]{\smash{{{\SetFigFont{12}{14.4}{\familydefault}{\mddefault}{\updefault}$*_X$}}}}}

\put(3923,549){\makebox(0,0)[lb]{\smash{{{\SetFigFont{12}{14.4}{\familydefault}{\mddefault}{\updefault}$l_{*_1,\R}$}}}}}

\put(523,399){\makebox(0,0)[lb]{\smash{{{\SetFigFont{12}{14.4}{\familydefault}{\mddefault}{\updefault}$t_3$}}}}}
\put(1418,394){\makebox(0,0)[lb]{\smash{{{\SetFigFont{12}{14.4}{\familydefault}{\mddefault}{\updefault}$t_2$}}}}}
\put(3033,389){\makebox(0,0)[lb]{\smash{{{\SetFigFont{12}{14.4}{\familydefault}{\mddefault}{\updefault}$t_1$}}}}}
\put(578,1529){\makebox(0,0)[lb]{\smash{{{\SetFigFont{12}{14.4}{\familydefault}{\mddefault}{\updefault}$c$}}}}}
\put(1108,1389){\makebox(0,0)[lb]{\smash{{{\SetFigFont{12}{14.4}{\familydefault}{\mddefault}{\updefault}$b$}}}}}
\put(2103,1499){\makebox(0,0)[lb]{\smash{{{\SetFigFont{12}{14.4}{\familydefault}{\mddefault}{\updefault}$a$}}}}}

\put(2998,2774){\makebox(0,0)[lb]{\smash{{{\SetFigFont{12}{14.4}{\familydefault}{\mddefault}{\updefault}$z$-plane }}}}}
\put(2600,2354){\makebox(0,0)[lb]{\smash{{{\SetFigFont{12}{14.4}{\familydefault}{\mddefault}{\updefault} = the complex pencil $l_{*_1,\C}$}}}}}

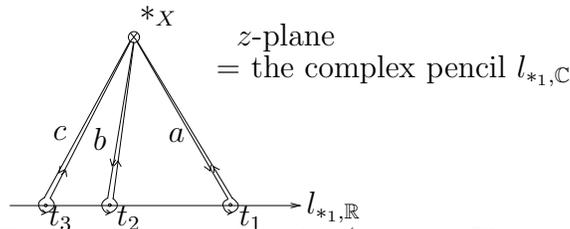
\put(-2100,54){\makebox(0,0)[lb]{\smash{{{\SetFigFont{12}{14.4}{\familydefault}{\mddefault}{\updefault}Figure 2. The generators $a$, $b$ and $c$ (see also Figures 3.1-3.11). }}}}}
\end{picture}
}\qquad \qquad\qquad \qquad \qquad \qquad 
}

 \newpage
The following Figure 3. briefly describes the  
real plane curve $C_{X,\R}=H_{X,\R}\cap D_X$ and the real basic pencil $l_{*_1,\R}$ 
inside the real plane $H_{X,\R}$ for all $X$ except for the cases 
$X\in\{ \mathrm{B_{vii}, H_{ii},H_{vi}},\mathrm{H_{viii}}\}$.

{\tiny
\ \\ 
\bigskip
\bigskip
\bigskip
\setlength{\unitlength}{0.00053333in}
\begingroup\makeatletter\ifx\SetFigFont\undefined%
\gdef\SetFigFont#1#2#3#4#5{%
  \reset@font\fontsize{#1}{#2pt}%
  \fontfamily{#3}\fontseries{#4}\fontshape{#5}%
  \selectfont}%
\fi\endgroup%
{\renewcommand{\dashlinestretch}{30}
\begin{picture}(2076,2095)(0,-10)
\put(1102.197,-1209.478){\arc{3587.153}{4.2141}{5.2060}}
\path(156,1768)(2056,284)
\path(150,292)(2064,1785)
\path(1844,1925)(1844,270)(1844,168)
\put(1862,1717){\makebox(0,0)[lb]{\smash{{{\SetFigFont{8}{9.6}{\familydefault}{\mddefault}{\updefault}$a$}}}}}
\put(1879,461){\makebox(0,0)[lb]{\smash{{{\SetFigFont{8}{9.6}{\familydefault}{\mddefault}{\updefault}$b$}}}}}
\put(1904,200){\makebox(0,0)[lb]{\smash{{{\SetFigFont{8}{9.6}{\familydefault}{\mddefault}{\updefault}$c$}}}}}
\put(1755,-71){\makebox(0,0)[lb]{\smash{{{\SetFigFont{8}{9.6}{\familydefault}{\mddefault}{\updefault}$l_{*_1,\R}$}}}}}
\put(0,1924){\makebox(0,0)[lb]{\smash{{{\SetFigFont{12}{14.4}{\familydefault}{\mddefault}{\updefault}$\mathrm{A_i}$.}}}}}
\end{picture}
}\qquad
\setlength{\unitlength}{0.00053333in}
\begingroup\makeatletter\ifx\SetFigFont\undefined%
\gdef\SetFigFont#1#2#3#4#5{%
  \reset@font\fontsize{#1}{#2pt}%
  \fontfamily{#3}\fontseries{#4}\fontshape{#5}%
  \selectfont}%
\fi\endgroup%
{\renewcommand{\dashlinestretch}{30}
\begin{picture}(2297,1660)(0,-10)
\put(1170,766){\ellipse{2222}{748}}
\path(1685,1482)(1685,173)(1686,166)
\path(35,1039)(37,1040)(41,1041)
	(48,1043)(58,1046)(73,1050)
	(91,1055)(112,1061)(136,1067)
	(162,1074)(188,1080)(216,1086)
	(244,1091)(273,1096)(302,1101)
	(331,1105)(362,1109)(394,1112)
	(428,1115)(464,1118)(502,1120)
	(543,1122)(572,1124)(601,1125)
	(632,1126)(664,1127)(697,1128)
	(732,1129)(768,1130)(806,1131)
	(844,1132)(884,1132)(925,1133)
	(967,1133)(1009,1134)(1052,1134)
	(1096,1134)(1140,1134)(1184,1134)
	(1228,1134)(1272,1134)(1315,1134)
	(1358,1134)(1401,1133)(1442,1133)
	(1483,1132)(1522,1132)(1560,1131)
	(1597,1130)(1633,1129)(1667,1128)
	(1700,1127)(1732,1126)(1762,1125)
	(1791,1124)(1818,1122)(1862,1120)
	(1903,1117)(1940,1114)(1975,1111)
	(2008,1107)(2038,1103)(2068,1098)
	(2097,1092)(2125,1086)(2152,1080)
	(2178,1073)(2202,1066)(2225,1059)
	(2244,1053)(2259,1048)(2271,1044)
	(2279,1041)(2283,1040)(2285,1039)
\put(1702,1285){\makebox(0,0)[lb]{\smash{{{\SetFigFont{9}{10.8}{\familydefault}{\mddefault}{\updefault}$a$}}}}}
\put(1702,969){\makebox(0,0)[lb]{\smash{{{\SetFigFont{9}{10.8}{\familydefault}{\mddefault}{\updefault}$b$}}}}}
\put(1699,479){\makebox(0,0)[lb]{\smash{{{\SetFigFont{9}{10.8}{\familydefault}{\mddefault}{\updefault}$c$}}}}}
\put(1605,-70){\makebox(0,0)[lb]{\smash{{{\SetFigFont{9}{10.8}{\familydefault}{\mddefault}{\updefault}$l_{*_1,\R}$}}}}}
\put(0,1489){\makebox(0,0)[lb]{\smash{{{\SetFigFont{12}{14.4}{\familydefault}{\mddefault}{\updefault}$\mathrm{A_{ii}}$.}}}}}
\end{picture}
}\qquad
\setlength{\unitlength}{0.00053333in}
\begingroup\makeatletter\ifx\SetFigFont\undefined%
\gdef\SetFigFont#1#2#3#4#5{%
  \reset@font\fontsize{#1}{#2pt}%
  \fontfamily{#3}\fontseries{#4}\fontshape{#5}%
  \selectfont}%
\fi\endgroup%
{\renewcommand{\dashlinestretch}{30}
\begin{picture}(2390,1758)(0,-10)
\put(1231.834,1878.383){\arc{2449.168}{0.7676}{2.4862}}
\path(2203,1139)(1307,12)
\path(75,654)(2378,654)
\path(1682,1671)(1682,280)
\put(1707,547){\makebox(0,0)[lb]{\smash{{{\SetFigFont{9}{10.8}{\familydefault}{\mddefault}{\updefault}$b$}}}}}
\put(1729,358){\makebox(0,0)[lb]{\smash{{{\SetFigFont{9}{10.8}{\familydefault}{\mddefault}{\updefault}$c$}}}}}
\put(1717,858){\makebox(0,0)[lb]{\smash{{{\SetFigFont{9}{10.8}{\familydefault}{\mddefault}{\updefault}$a$}}}}}
\put(1603,23){\makebox(0,0)[lb]{\smash{{{\SetFigFont{9}{10.8}{\familydefault}{\mddefault}{\updefault}$l_{*_1,\R}$}}}}}
\put(0,1587){\makebox(0,0)[lb]{\smash{{{\SetFigFont{12}{14.4}{\familydefault}{\mddefault}{\updefault}$\mathrm{B_i}$.}}}}}
\end{picture}
}\\\\
\bigskip
\setlength{\unitlength}{0.00053333in}
\begingroup\makeatletter\ifx\SetFigFont\undefined%
\gdef\SetFigFont#1#2#3#4#5{%
  \reset@font\fontsize{#1}{#2pt}%
  \fontfamily{#3}\fontseries{#4}\fontshape{#5}%
  \selectfont}%
\fi\endgroup%
{\renewcommand{\dashlinestretch}{30}
\begin{picture}(2263,1585)(0,-10)
\path(1542,1506)(1542,92)
\path(75,647)(2251,647)
\path(2150,848)(2149,848)(2146,847)
	(2140,845)(2132,843)(2121,840)
	(2108,837)(2092,833)(2073,828)
	(2053,823)(2032,817)(2009,811)
	(1984,805)(1958,798)(1929,791)
	(1898,784)(1865,775)(1828,767)
	(1787,757)(1745,747)(1702,737)
	(1663,728)(1630,720)(1603,715)
	(1583,710)(1569,707)(1561,706)
	(1557,705)(1556,705)(1555,706)
	(1552,706)(1546,705)(1535,702)
	(1517,698)(1491,693)(1455,685)
	(1409,674)(1353,662)(1290,647)
	(1242,636)(1194,625)(1147,614)
	(1101,604)(1058,595)(1018,587)
	(981,580)(946,574)(914,568)
	(884,564)(856,559)(829,556)
	(804,552)(779,549)(755,546)
	(731,543)(706,539)(681,535)
	(655,530)(629,524)(601,517)
	(574,510)(545,501)(517,491)
	(491,480)(466,468)(445,456)
	(429,443)(419,430)(416,417)
	(417,406)(422,395)(429,387)
	(438,379)(448,373)(457,368)
	(467,364)(477,361)(488,359)
	(498,357)(509,355)(521,353)
	(534,351)(549,349)(567,346)
	(587,343)(612,339)(641,335)
	(675,329)(715,323)(760,316)
	(813,308)(871,300)(935,292)
	(987,286)(1040,281)(1093,276)
	(1147,272)(1199,268)(1252,264)
	(1303,261)(1354,258)(1404,256)
	(1453,253)(1501,252)(1549,250)
	(1597,248)(1644,247)(1690,246)
	(1735,245)(1780,244)(1823,244)
	(1865,243)(1905,243)(1943,242)
	(1979,242)(2011,242)(2041,241)
	(2067,241)(2089,241)(2107,241)
	(2122,241)(2133,241)(2141,241)
	(2146,241)(2149,241)(2150,241)
\put(1594,292){\makebox(0,0)[lb]{\smash{{{\SetFigFont{8}{9.6}{\familydefault}{\mddefault}{\updefault}$c$}}}}}
\put(1594,546){\makebox(0,0)[lb]{\smash{{{\SetFigFont{8}{9.6}{\familydefault}{\mddefault}{\updefault}$b$}}}}}
\put(1594,798){\makebox(0,0)[lb]{\smash{{{\SetFigFont{8}{9.6}{\familydefault}{\mddefault}{\updefault}$a$}}}}}
\put(1452,-100){\makebox(0,0)[lb]{\smash{{{\SetFigFont{8}{9.6}{\familydefault}{\mddefault}{\updefault}$l_{*_1,\R}$}}}}}
\put(0,1414){\makebox(0,0)[lb]{\smash{{{\SetFigFont{12}{14.4}{\familydefault}{\mddefault}{\updefault}$\mathrm{B_{ii}}$.}}}}}
\end{picture}
}\qquad
\setlength{\unitlength}{0.00053333in}
\begingroup\makeatletter\ifx\SetFigFont\undefined%
\gdef\SetFigFont#1#2#3#4#5{%
  \reset@font\fontsize{#1}{#2pt}%
  \fontfamily{#3}\fontseries{#4}\fontshape{#5}%
  \selectfont}%
\fi\endgroup%
{\renewcommand{\dashlinestretch}{30}
\begin{picture}(2283,1615)(0,-10)
\path(75,651)(2271,651)
\path(1301,1519)(1301,192)
\path(2220,805)(2219,804)(2218,803)
	(2216,801)(2212,797)(2206,792)
	(2199,786)(2189,778)(2179,769)
	(2166,760)(2151,750)(2135,739)
	(2118,728)(2099,717)(2078,707)
	(2055,697)(2030,687)(2003,678)
	(1973,670)(1940,663)(1904,657)
	(1863,652)(1818,649)(1769,648)
	(1715,648)(1658,651)(1606,655)
	(1555,661)(1504,667)(1457,674)
	(1412,681)(1371,687)(1334,694)
	(1300,700)(1269,705)(1242,710)
	(1217,715)(1195,720)(1174,724)
	(1154,728)(1135,732)(1116,737)
	(1096,741)(1076,747)(1055,753)
	(1032,759)(1006,767)(979,776)
	(949,787)(916,799)(881,813)
	(844,828)(804,846)(765,865)
	(725,885)(688,907)(652,931)
	(621,955)(594,978)(572,1000)
	(554,1021)(540,1040)(530,1057)
	(522,1073)(517,1087)(514,1100)
	(512,1112)(512,1124)(512,1135)
	(514,1145)(516,1155)(518,1166)
	(521,1176)(523,1187)(527,1199)
	(530,1211)(533,1223)(537,1237)
	(542,1251)(547,1265)(554,1279)
	(562,1292)(573,1305)(586,1315)
	(605,1324)(626,1329)(647,1333)
	(668,1334)(688,1334)(705,1334)
	(721,1333)(735,1331)(747,1330)
	(758,1328)(768,1326)(778,1324)
	(788,1321)(800,1317)(813,1312)
	(828,1306)(847,1299)(869,1289)
	(895,1276)(927,1260)(963,1241)
	(1004,1219)(1049,1192)(1097,1162)
	(1136,1136)(1175,1108)(1214,1080)
	(1251,1051)(1286,1022)(1320,994)
	(1353,965)(1385,937)(1415,909)
	(1445,881)(1473,853)(1501,826)
	(1528,799)(1555,772)(1580,745)
	(1604,720)(1628,695)(1650,671)
	(1670,649)(1689,629)(1705,610)
	(1720,594)(1732,581)(1741,570)
	(1749,562)(1754,556)(1757,552)
	(1759,550)(1760,549)
\put(1351,702){\makebox(0,0)[lb]{\smash{{{\SetFigFont{9}{10.8}{\familydefault}{\mddefault}{\updefault}$b$}}}}}
\put(1209,-89){\makebox(0,0)[lb]{\smash{{{\SetFigFont{9}{10.8}{\familydefault}{\mddefault}{\updefault}$l_{*_1,\R}$}}}}}
\put(1351,549){\makebox(0,0)[lb]{\smash{{{\SetFigFont{9}{10.8}{\familydefault}{\mddefault}{\updefault}$c$}}}}}
\put(1351,1008){\makebox(0,0)[lb]{\smash{{{\SetFigFont{9}{10.8}{\familydefault}{\mddefault}{\updefault}$a$}}}}}
\put(0,1444){\makebox(0,0)[lb]{\smash{{{\SetFigFont{12}{14.4}{\familydefault}{\mddefault}{\updefault}$\mathrm{B_{iii}}$.}}}}}
\end{picture}
}\qquad
\setlength{\unitlength}{0.00053333in}
\begingroup\makeatletter\ifx\SetFigFont\undefined%
\gdef\SetFigFont#1#2#3#4#5{%
  \reset@font\fontsize{#1}{#2pt}%
  \fontfamily{#3}\fontseries{#4}\fontshape{#5}%
  \selectfont}%
\fi\endgroup%
{\renewcommand{\dashlinestretch}{30}
\begin{picture}(2346,1699)(0,-10)
\put(1158.128,2259.730){\arc{2671.477}{0.6456}{2.0808}}
\path(1580,1561)(1580,147)
\path(75,1142)(2334,1142)
\path(2040,197)(429,1139)
\put(1487,-82){\makebox(0,0)[lb]{\smash{{{\SetFigFont{9}{10.8}{\familydefault}{\mddefault}{\updefault}$l_{*_1,\R}$}}}}}
\put(1635,1195){\makebox(0,0)[lb]{\smash{{{\SetFigFont{9}{10.8}{\familydefault}{\mddefault}{\updefault}$a$}}}}}
\put(1635,880){\makebox(0,0)[lb]{\smash{{{\SetFigFont{9}{10.8}{\familydefault}{\mddefault}{\updefault}$b$}}}}}
\put(1635,513){\makebox(0,0)[lb]{\smash{{{\SetFigFont{9}{10.8}{\familydefault}{\mddefault}{\updefault}$c$}}}}}
\put(0,1528){\makebox(0,0)[lb]{\smash{{{\SetFigFont{12}{14.4}{\familydefault}{\mddefault}{\updefault}$\mathrm{B_{iv}}$.}}}}}
\end{picture}
}\\\\

\bigskip
\setlength{\unitlength}{0.00053333in}
\begingroup\makeatletter\ifx\SetFigFont\undefined%
\gdef\SetFigFont#1#2#3#4#5{%
  \reset@font\fontsize{#1}{#2pt}%
  \fontfamily{#3}\fontseries{#4}\fontshape{#5}%
  \selectfont}%
\fi\endgroup%
{\renewcommand{\dashlinestretch}{30}
\begin{picture}(2391,1936)(0,-10)
\path(565,1716)(565,179)
\path(189,474)(2379,947)
\path(178,1310)(180,1310)(185,1310)
	(194,1310)(207,1311)(223,1311)
	(243,1311)(266,1311)(291,1311)
	(317,1310)(345,1308)(373,1306)
	(402,1303)(433,1299)(465,1293)
	(498,1286)(533,1277)(568,1266)
	(601,1254)(630,1241)(655,1230)
	(676,1220)(692,1213)(705,1206)
	(716,1201)(725,1197)(732,1193)
	(739,1189)(746,1185)(754,1179)
	(761,1171)(770,1162)(778,1149)
	(786,1134)(792,1116)(794,1097)
	(792,1078)(786,1059)(778,1044)
	(769,1031)(761,1020)(753,1012)
	(746,1005)(739,1000)(732,995)
	(724,990)(715,985)(704,979)
	(691,972)(675,963)(654,952)
	(630,940)(601,926)(568,912)
	(534,900)(499,889)(466,881)
	(435,874)(405,869)(376,865)
	(348,862)(321,860)(295,859)
	(271,858)(249,857)(229,857)
	(213,857)(200,857)(192,857)
	(187,857)(185,857)
\put(624,1298){\makebox(0,0)[lb]{\smash{{{\SetFigFont{10}{12.0}{\familydefault}{\mddefault}{\updefault}$a$}}}}}
\put(632,805){\makebox(0,0)[lb]{\smash{{{\SetFigFont{10}{12.0}{\familydefault}{\mddefault}{\updefault}$b$}}}}}
\put(624,647){\makebox(0,0)[lb]{\smash{{{\SetFigFont{10}{12.0}{\familydefault}{\mddefault}{\updefault}$c$}}}}}
\put(451,-85){\makebox(0,0)[lb]{\smash{{{\SetFigFont{10}{12.0}{\familydefault}{\mddefault}{\updefault}$l_{*_1,\R}$}}}}}
\put(0,1765){\makebox(0,0)[lb]{\smash{{{\SetFigFont{12}{14.4}{\familydefault}{\mddefault}{\updefault}$\mathrm{B_v}$.}}}}}
\end{picture}
}\qquad
\setlength{\unitlength}{0.00053333in}
\begingroup\makeatletter\ifx\SetFigFont\undefined%
\gdef\SetFigFont#1#2#3#4#5{%
  \reset@font\fontsize{#1}{#2pt}%
  \fontfamily{#3}\fontseries{#4}\fontshape{#5}%
  \selectfont}%
\fi\endgroup%
{\renewcommand{\dashlinestretch}{30}
\begin{picture}(2263,1861)(0,-10)
\path(1898,1710)(1898,179)
\path(2248,1545)(2246,1545)(2241,1545)
	(2232,1545)(2219,1545)(2203,1545)
	(2183,1545)(2160,1545)(2136,1545)
	(2110,1544)(2083,1542)(2055,1540)
	(2027,1537)(1998,1534)(1968,1529)
	(1937,1523)(1905,1515)(1874,1506)
	(1845,1496)(1820,1486)(1799,1476)
	(1782,1468)(1769,1462)(1759,1457)
	(1752,1453)(1746,1450)(1741,1446)
	(1736,1443)(1732,1440)(1727,1435)
	(1722,1429)(1717,1421)(1711,1411)
	(1706,1398)(1702,1384)(1701,1368)
	(1704,1350)(1711,1334)(1719,1320)
	(1726,1308)(1733,1300)(1740,1293)
	(1746,1288)(1752,1283)(1758,1278)
	(1766,1273)(1776,1267)(1790,1260)
	(1807,1251)(1830,1240)(1858,1227)
	(1890,1215)(1921,1205)(1952,1197)
	(1982,1191)(2011,1186)(2038,1183)
	(2065,1181)(2091,1179)(2116,1178)
	(2140,1178)(2163,1178)(2184,1178)
	(2202,1178)(2218,1179)(2229,1179)
	(2237,1180)(2242,1180)(2244,1180)
\path(2251,402)(2250,402)(2248,403)
	(2243,404)(2236,405)(2227,408)
	(2214,411)(2198,414)(2179,419)
	(2157,423)(2133,429)(2106,435)
	(2077,441)(2046,448)(2013,455)
	(1979,462)(1944,469)(1907,476)
	(1869,483)(1830,490)(1790,497)
	(1748,504)(1704,511)(1659,518)
	(1611,525)(1561,531)(1508,538)
	(1453,544)(1397,550)(1340,556)
	(1279,561)(1221,566)(1167,569)
	(1117,572)(1072,574)(1032,576)
	(997,576)(966,577)(939,577)
	(916,576)(895,575)(877,574)
	(860,573)(844,572)(829,571)
	(813,570)(797,568)(780,567)
	(761,566)(741,565)(719,564)
	(695,563)(670,562)(643,562)
	(615,561)(587,561)(562,560)
	(539,560)(509,559)(497,558)
	(500,558)(513,558)(533,559)
	(557,559)(584,561)(612,561)
	(641,562)(669,561)(697,560)
	(725,556)(751,551)(774,544)
	(792,536)(806,529)(817,522)
	(825,516)(831,511)(835,507)
	(838,504)(840,501)(841,498)
	(842,494)(844,490)(845,484)
	(846,477)(847,468)(847,457)
	(845,446)(840,433)(833,422)
	(827,414)(822,408)(819,405)
	(817,403)(815,401)(813,400)
	(808,398)(801,394)(790,389)
	(772,381)(749,371)(719,359)
	(693,350)(668,341)(645,334)
	(626,329)(610,324)(597,321)
	(587,318)(579,317)(571,315)
	(563,314)(554,312)(543,310)
	(529,307)(511,303)(488,298)
	(460,292)(428,285)(393,277)
	(354,269)(318,261)(285,255)
	(256,249)(229,245)(205,240)
	(182,237)(162,233)(143,230)
	(126,228)(112,226)(101,224)
	(93,223)(89,222)(87,222)
\put(1933,1584){\makebox(0,0)[lb]{\smash{{{\SetFigFont{9}{10.8}{\familydefault}{\mddefault}{\updefault}$a$}}}}}
\put(1933,1219){\makebox(0,0)[lb]{\smash{{{\SetFigFont{9}{10.8}{\familydefault}{\mddefault}{\updefault}$b$}}}}}
\put(1933,508){\makebox(0,0)[lb]{\smash{{{\SetFigFont{9}{10.8}{\familydefault}{\mddefault}{\updefault}$c$}}}}}
\put(1803,-85){\makebox(0,0)[lb]{\smash{{{\SetFigFont{9}{10.8}{\familydefault}{\mddefault}{\updefault}$l_{*_1,\R}$}}}}}
\put(0,1690){\makebox(0,0)[lb]{\smash{{{\SetFigFont{12}{14.4}{\familydefault}{\mddefault}{\updefault}$\mathrm{B_{vi}}$.}}}}}
\end{picture}
}\qquad
\setlength{\unitlength}{0.00053333in}
\begingroup\makeatletter\ifx\SetFigFont\undefined%
\gdef\SetFigFont#1#2#3#4#5{%
  \reset@font\fontsize{#1}{#2pt}%
  \fontfamily{#3}\fontseries{#4}\fontshape{#5}%
  \selectfont}%
\fi\endgroup%
{\renewcommand{\dashlinestretch}{30}
\begin{picture}(1976,2128)(0,-10)
\put(799.263,2595.231){\arc{2418.076}{0.7535}{1.9983}}
\path(79,358)(1856,1917)
\path(189,1586)(1964,140)
\path(1190,2025)(1190,281)(1190,173)
\put(1202,1530){\makebox(0,0)[lb]{\smash{{{\SetFigFont{9}{10.8}{\familydefault}{\mddefault}{\updefault}$a$}}}}}
\put(1224,1237){\makebox(0,0)[lb]{\smash{{{\SetFigFont{9}{10.8}{\familydefault}{\mddefault}{\updefault}$b$}}}}}
\put(1231,750){\makebox(0,0)[lb]{\smash{{{\SetFigFont{9}{10.8}{\familydefault}{\mddefault}{\updefault}$c$}}}}}
\put(1100,-83){\makebox(0,0)[lb]{\smash{{{\SetFigFont{9}{10.8}{\familydefault}{\mddefault}{\updefault}$l_{*_1,\R}$}}}}}
\put(0,1957){\makebox(0,0)[lb]{\smash{{{\SetFigFont{12}{14.4}{\familydefault}{\mddefault}{\updefault}$\mathrm{H_i}$.}}}}}
\end{picture}
}\\\\

\bigskip
\setlength{\unitlength}{0.00053333in}
\begingroup\makeatletter\ifx\SetFigFont\undefined%
\gdef\SetFigFont#1#2#3#4#5{%
  \reset@font\fontsize{#1}{#2pt}%
  \fontfamily{#3}\fontseries{#4}\fontshape{#5}%
  \selectfont}%
\fi\endgroup%
{\renewcommand{\dashlinestretch}{30}
\begin{picture}(2243,1763)(0,-10)
\path(104,617)(1973,1536)
\path(2231,700)(1825,751)(504,812)
\path(1189,1648)(1189,183)
\put(731.519,2192.661){\arc{2784.431}{0.5108}{1.6517}}
\put(1222,1275){\makebox(0,0)[lb]{\smash{{{\SetFigFont{9}{10.8}{\familydefault}{\mddefault}{\updefault}$a$}}}}}
\put(1216,1024){\makebox(0,0)[lb]{\smash{{{\SetFigFont{9}{10.8}{\familydefault}{\mddefault}{\updefault}$b$}}}}}
\put(1230,664){\makebox(0,0)[lb]{\smash{{{\SetFigFont{9}{10.8}{\familydefault}{\mddefault}{\updefault}$c$}}}}}
\put(1003,-82){\makebox(0,0)[lb]{\smash{{{\SetFigFont{9}{10.8}{\familydefault}{\mddefault}{\updefault}$l_{*_1,\R}$}}}}}
\put(0,1592){\makebox(0,0)[lb]{\smash{{{\SetFigFont{12}{14.4}{\familydefault}{\mddefault}{\updefault}$\mathrm{H_{iii}}$.}}}}}
\end{picture}
}\qquad
\setlength{\unitlength}{0.00053333in}
\begingroup\makeatletter\ifx\SetFigFont\undefined%
\gdef\SetFigFont#1#2#3#4#5{%
  \reset@font\fontsize{#1}{#2pt}%
  \fontfamily{#3}\fontseries{#4}\fontshape{#5}%
  \selectfont}%
\fi\endgroup%
{\renewcommand{\dashlinestretch}{30}
\begin{picture}(1506,2080)(0,-10)
\put(1623.391,2191.375){\arc{3133.924}{1.7164}{2.7943}}
\path(1174,1914)(1174,186)
\path(1397,640)(1396,641)(1393,642)
	(1388,645)(1380,650)(1370,656)
	(1357,664)(1341,674)(1322,685)
	(1302,698)(1280,713)(1256,728)
	(1232,744)(1207,762)(1181,780)
	(1155,800)(1128,821)(1100,843)
	(1072,867)(1042,893)(1012,921)
	(981,952)(949,984)(918,1018)
	(886,1055)(856,1091)(831,1124)
	(808,1153)(790,1178)(774,1199)
	(762,1216)(751,1231)(743,1243)
	(735,1253)(729,1263)(723,1272)
	(717,1282)(711,1294)(704,1307)
	(697,1322)(690,1340)(682,1361)
	(675,1385)(669,1411)(664,1439)
	(662,1466)(663,1487)(667,1506)
	(672,1523)(679,1538)(686,1551)
	(693,1562)(700,1572)(708,1580)
	(715,1588)(722,1594)(728,1600)
	(735,1605)(742,1610)(749,1614)
	(757,1618)(765,1622)(774,1625)
	(784,1627)(795,1629)(808,1630)
	(822,1629)(838,1627)(856,1623)
	(875,1616)(896,1607)(918,1594)
	(942,1576)(965,1556)(987,1535)
	(1005,1516)(1022,1497)(1035,1481)
	(1047,1467)(1056,1454)(1063,1444)
	(1069,1435)(1074,1426)(1079,1418)
	(1083,1409)(1088,1399)(1094,1387)
	(1101,1373)(1110,1354)(1122,1332)
	(1135,1304)(1152,1271)(1171,1232)
	(1192,1187)(1215,1137)(1239,1082)
	(1259,1033)(1279,985)(1297,936)
	(1315,889)(1331,844)(1346,800)
	(1360,757)(1373,716)(1385,675)
	(1396,636)(1407,598)(1417,560)
	(1427,523)(1436,488)(1445,454)
	(1453,421)(1461,391)(1468,362)
	(1474,337)(1479,315)(1484,296)
	(1487,280)(1490,268)(1492,260)
	(1493,254)(1494,251)(1494,250)
\put(1026,-80){\makebox(0,0)[lb]{\smash{{{\SetFigFont{10}{12.0}{\familydefault}{\mddefault}{\updefault}$l_{*_1,\R}$}}}}}
\put(1195,1299){\makebox(0,0)[lb]{\smash{{{\SetFigFont{10}{12.0}{\familydefault}{\mddefault}{\updefault}$a$}}}}}
\put(940,976){\makebox(0,0)[lb]{\smash{{{\SetFigFont{10}{12.0}{\familydefault}{\mddefault}{\updefault}$b$}}}}}
\put(906,531){\makebox(0,0)[lb]{\smash{{{\SetFigFont{10}{12.0}{\familydefault}{\mddefault}{\updefault}$c$}}}}}
\put(0,1909){\makebox(0,0)[lb]{\smash{{{\SetFigFont{12}{14.4}{\familydefault}{\mddefault}{\updefault}$\mathrm{H_{iv}}$.}}}}}
\end{picture}
}\qquad\qquad\qquad
\setlength{\unitlength}{0.00053333in}
\begingroup\makeatletter\ifx\SetFigFont\undefined%
\gdef\SetFigFont#1#2#3#4#5{%
  \reset@font\fontsize{#1}{#2pt}%
  \fontfamily{#3}\fontseries{#4}\fontshape{#5}%
  \selectfont}%
\fi\endgroup%
{\renewcommand{\dashlinestretch}{30}
\begin{picture}(1821,1965)(0,-10)
\put(2718.396,4662.907){\arc{7898.193}{1.8021}{2.2663}}
\path(863,1869)(863,181)
\path(1080,1074)(1079,1072)(1078,1069)
	(1076,1063)(1073,1054)(1068,1042)
	(1062,1027)(1054,1008)(1045,988)
	(1036,965)(1025,941)(1013,916)
	(1000,890)(985,863)(970,836)
	(952,807)(933,778)(912,748)
	(888,717)(861,685)(832,652)
	(800,619)(767,588)(735,560)
	(705,535)(678,513)(655,494)
	(635,479)(617,466)(603,455)
	(590,446)(579,437)(569,430)
	(558,423)(547,417)(535,410)
	(521,403)(506,396)(487,389)
	(466,382)(442,376)(416,371)
	(389,368)(363,369)(342,373)
	(324,379)(307,387)(293,397)
	(281,406)(270,416)(261,426)
	(253,436)(247,445)(241,454)
	(236,463)(231,472)(227,481)
	(224,490)(220,500)(217,510)
	(215,521)(213,533)(211,547)
	(210,561)(210,578)(212,596)
	(215,615)(220,636)(227,658)
	(238,681)(254,707)(273,732)
	(293,754)(313,773)(331,789)
	(348,802)(363,812)(377,820)
	(390,826)(402,832)(414,836)
	(427,841)(440,846)(455,851)
	(471,858)(491,867)(514,878)
	(540,890)(571,905)(605,922)
	(641,940)(680,958)(722,976)
	(762,993)(799,1007)(833,1018)
	(865,1028)(893,1036)(920,1043)
	(945,1048)(969,1053)(991,1057)
	(1012,1060)(1030,1063)(1047,1065)
	(1061,1066)(1073,1067)(1081,1068)
	(1086,1069)(1090,1069)(1091,1069)
\put(767,-80){\makebox(0,0)[lb]{\smash{{{\SetFigFont{10}{12.0}{\familydefault}{\mddefault}{\updefault}$l_{*_1,\R}$}}}}}
\put(888,1190){\makebox(0,0)[lb]{\smash{{{\SetFigFont{10}{12.0}{\familydefault}{\mddefault}{\updefault}$a$}}}}}
\put(905,536){\makebox(0,0)[lb]{\smash{{{\SetFigFont{10}{12.0}{\familydefault}{\mddefault}{\updefault}$c$}}}}}
\put(889,915){\makebox(0,0)[lb]{\smash{{{\SetFigFont{10}{12.0}{\familydefault}{\mddefault}{\updefault}$b$}}}}}
\put(0,1794){\makebox(0,0)[lb]{\smash{{{\SetFigFont{12}{14.4}{\familydefault}{\mddefault}{\updefault}$\mathrm{H_{vii}}$.}}}}}
\end{picture}
}\\\\
}
\begin{center}
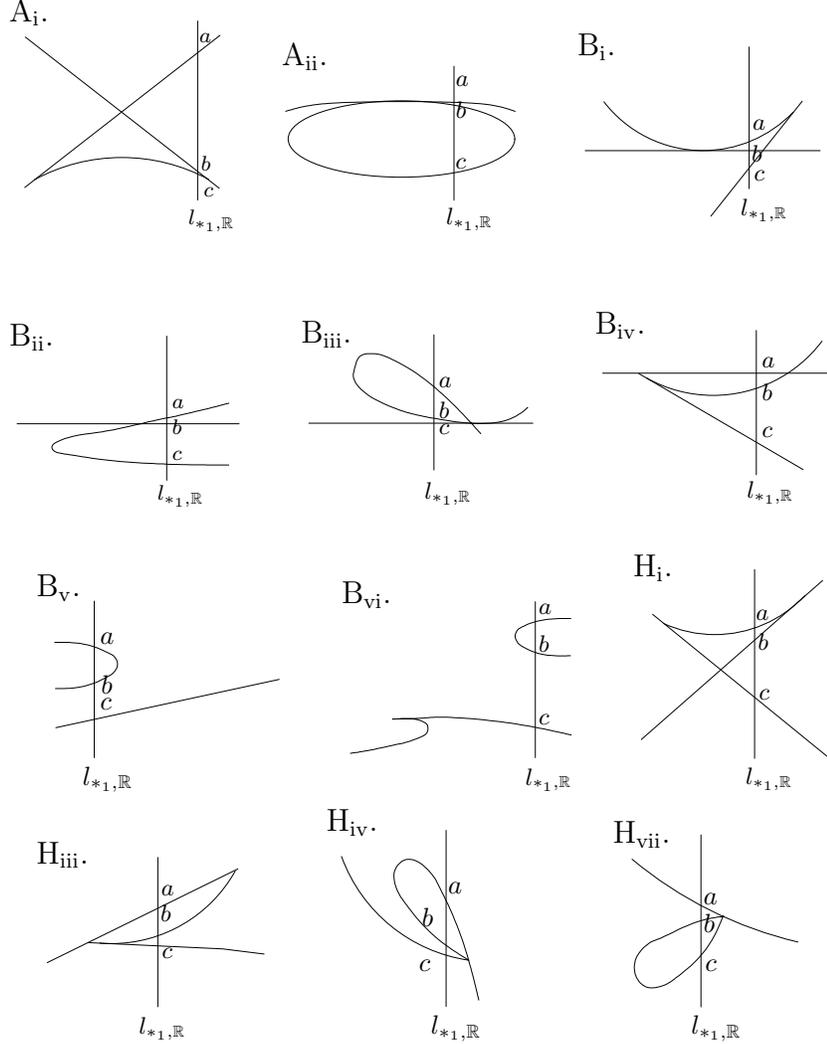
Figure 3. Real plane curve $C_{X,\R}$ and the pencil $l_{*_1,\R}$ in the real plane $H_{X,\R}$.
\end{center}

For the remaining cases $X\in\{\mathrm{B_{vii}, H_{ii},H_{vi},H_{viii}}\}$, 
some more careful considerations are necessary. 
We briefly indicate the choices of $*_1$ in the (complex) line $L_X$, 
the base point $*_X$ and then the paths $a,b,c$ in the basic (complex) pencil 
$l_{*_1}$ along the intervals connecting $*_X$ and the three points $l_{*1}\cap D_X$ as in Figure 4.1-4.5. We indicate also the bifurcation points 
$L_X\cap B_X$ and the paths $\gamma_i$ which shall be used in the 
step 3.~of Zariski-van Kampen method.

{\tiny
\bigskip
\bigskip
\setlength{\unitlength}{0.00063333in}
\begingroup\makeatletter\ifx\SetFigFont\undefined%
\gdef\SetFigFont#1#2#3#4#5{%
  \reset@font\fontsize{#1}{#2pt}%
  \fontfamily{#3}\fontseries{#4}\fontshape{#5}%
  \selectfont}%
\fi\endgroup%
{\renewcommand{\dashlinestretch}{30}
\begin{picture}(3249,2474)(0,-10)
\put(210,1544){\ellipse{52}{52}}
\put(1059,1292){\ellipse{14}{14}}
\put(908,1534){\ellipse{14}{14}}
\put(1070,1825){\ellipse{14}{14}}
\path(231,1524)(1029,1278)
\path(231,1526)(1029,1306)
\path(240,1544)(879,1544)
\path(235,1542)(235,1541)(236,1542)(879,1513)
\path(231,1562)(1036,1831)
\path(231,1558)(1038,1804)
\path(231,1527)(190,1558)
\path(191,1528)(231,1558)
\path(1098,1271)(1088,1290)(1074,1276)
\path(932,1324)(913,1345)(941,1349)
\path(898,1300)(928,1295)(913,1323)
\path(943,1510)(934,1531)(919,1513)
\path(783,1502)(801,1513)(786,1524)
\path(806,1528)(783,1535)(801,1551)
\path(1112,1810)(1102,1831)(1087,1814)
\path(897,1733)(918,1758)(889,1764)
\path(915,1787)(888,1785)(902,1812)
\path(154,2011)(1378,2011)(1378,1126)
	(154,1126)(154,2011)
\path(1021,1279)(1022,1278)(1024,1274)
	(1030,1267)(1036,1260)(1045,1255)
	(1056,1253)(1063,1254)(1069,1256)
	(1073,1259)(1077,1260)(1080,1262)
	(1083,1264)(1085,1267)(1088,1271)
	(1091,1276)(1093,1283)(1092,1292)
	(1090,1299)(1087,1304)(1084,1308)
	(1081,1312)(1077,1316)(1071,1321)
	(1063,1325)(1055,1326)(1048,1325)
	(1041,1322)(1036,1318)(1031,1313)
	(1028,1310)(1026,1308)
\path(870,1520)(872,1517)(874,1513)
	(878,1507)(883,1501)(889,1496)
	(896,1493)(904,1491)(912,1492)
	(919,1495)(924,1497)(928,1500)
	(931,1502)(934,1506)(938,1512)
	(941,1520)(941,1529)(939,1536)
	(937,1541)(935,1545)(932,1548)
	(928,1553)(923,1558)(915,1562)
	(904,1563)(894,1561)(885,1555)
	(878,1550)(874,1547)(873,1546)
\path(1041,1814)(1042,1812)(1044,1808)
	(1047,1803)(1050,1797)(1055,1793)
	(1061,1790)(1068,1789)(1076,1791)
	(1082,1794)(1087,1796)(1091,1798)
	(1095,1801)(1099,1805)(1102,1811)
	(1105,1819)(1105,1828)(1102,1835)
	(1100,1841)(1098,1845)(1095,1850)
	(1091,1854)(1085,1859)(1077,1862)
	(1065,1861)(1056,1855)(1047,1847)
	(1040,1839)(1036,1834)(1035,1833)
\put(1126,1254){\makebox(0,0)[lb]{\smash{{{\SetFigFont{9}{10.8}{\familydefault}{\mddefault}{\updefault}$c$}}}}}
\put(952,1492){\makebox(0,0)[lb]{\smash{{{\SetFigFont{9}{10.8}{\familydefault}{\mddefault}{\updefault}$b$}}}}}
\put(1126,1810){\makebox(0,0)[lb]{\smash{{{\SetFigFont{9}{10.8}{\familydefault}{\mddefault}{\updefault}$a$}}}}}
\put(1662,936){\ellipse{84}{84}}
\put(1662,343){\ellipse{84}{84}}
\put(1905,670){\ellipse{84}{84}}
\put(1662,339){\ellipse{24}{22}}
\put(1900,670){\ellipse{24}{22}}
\put(1658,939){\ellipse{24}{22}}
\path(150,670)(3234,670)
\path(3140.460,646.615)(3234.000,670.000)(3140.460,693.385)
\path(1928.385,2022.460)(1905.000,2116.000)(1881.615,2022.460)
\path(1905,2116)(1905,12)
\path(661,696)(1632,916)
\path(1865,651)(656,661)
\path(654,654)(1624,343)
\path(496,1123)(575,893)
\path(537.690,929.637)(575.000,893.000)(581.923,944.830)
\put(1818,2131){\makebox(0,0)[lb]{\smash{{{\SetFigFont{10}{12.0}{\familydefault}{\mddefault}{\updefault}$I_m$}}}}}
\put(1645,1001){\makebox(0,0)[lb]{\smash{{{\SetFigFont{10}{12.0}{\familydefault}{\mddefault}{\updefault}$\gamma_1$}}}}}
\put(1935,787){\makebox(0,0)[lb]{\smash{{{\SetFigFont{10}{12.0}{\familydefault}{\mddefault}{\updefault}$\gamma_2$}}}}}
\put(1737,289){\makebox(0,0)[lb]{\smash{{{\SetFigFont{10}{12.0}{\familydefault}{\mddefault}{\updefault}$\gamma_3$}}}}}
\put(1967,553){\makebox(0,0)[lb]{\smash{{{\SetFigFont{10}{12.0}{\familydefault}{\mddefault}{\updefault}0}}}}}
\put(496,502){\makebox(0,0)[lb]{\smash{{{\SetFigFont{10}{12.0}{\familydefault}{\mddefault}{\updefault}$\mathrm{B_{vii}}$}}}}}
\put(557,708){\makebox(0,0)[lb]{\smash{{{\SetFigFont{10}{12.0}{\familydefault}{\mddefault}{\updefault}$*_1$}}}}}
\put(3249,619){\makebox(0,0)[lb]{\smash{{{\SetFigFont{10}{12.0}{\familydefault}{\mddefault}{\updefault}$Re$}}}}}
\put(0,2303){\makebox(0,0)[lb]{\smash{{{\SetFigFont{12}{14.4}{\familydefault}{\mddefault}{\updefault}$\mathrm{B_{vii}}$.}}}}}

\put(2650,1719){\makebox(0,0)[lb]{\smash{{{\SetFigFont{10}{12.0}{\familydefault}{\mddefault}{\updefault}$L_{\mathrm{B_{vii}},\C}$}}}}}
\put(220,1629){\makebox(0,0)[lb]{\smash{{{\SetFigFont{10}{12.0}{\familydefault}{\mddefault}{\updefault}$*_{\mathrm{B_{vii}}}$}}}}}
\put(700,2139){\makebox(0,0)[lb]{\smash{{{\SetFigFont{10}{12.0}{\familydefault}{\mddefault}{\updefault}$l_{*_1,\C}$}}}}}
\end{picture}
}
\setlength{\unitlength}{0.00060333in}
\begingroup\makeatletter\ifx\SetFigFont\undefined%
\gdef\SetFigFont#1#2#3#4#5{%
  \reset@font\fontsize{#1}{#2pt}%
  \fontfamily{#3}\fontseries{#4}\fontshape{#5}%
  \selectfont}%
\fi\endgroup%
{\renewcommand{\dashlinestretch}{30}
\begin{picture}(3618,2433)(0,-10)
\put(2430,1207){\ellipse{14}{14}}
\put(2896,2096){\ellipse{56}{56}}
\put(3032,1352){\ellipse{14}{14}}
\put(3317,1352){\ellipse{14}{14}}
\path(2874,2073)(2421,1246)
\path(2877,2071)(2454,1244)
\path(2506,1413)(2501,1379)(2530,1396)
\path(2512,1369)(2536,1391)(2541,1359)
\path(2896,2062)(3052,1387)
\path(2894,2068)(2893,2067)(2894,2066)(3028,1388)
\path(2915,2071)(3322,1391)
\path(2911,2071)(3307,1382)
\path(2878,2073)(2911,2117)
\path(2878,2116)(2911,2073)
\path(2409,1158)(2430,1167)(2412,1184)
\path(3019,1299)(3041,1308)(3023,1325)
\path(3019,1489)(3027,1514)(3045,1494)
\path(2998,1473)(3011,1452)(3023,1468)
\path(3246,1449)(3277,1425)(3281,1456)
\path(3292,1419)(3290,1448)(3321,1432)
\path(3310,1304)(3332,1314)(3313,1331)
\path(2421,1248)(2418,1246)(2414,1243)
	(2408,1239)(2402,1233)(2397,1227)
	(2393,1219)(2392,1210)(2394,1201)
	(2397,1193)(2400,1188)(2402,1184)
	(2405,1180)(2409,1176)(2416,1172)
	(2425,1169)(2435,1170)(2443,1172)
	(2449,1176)(2453,1179)(2457,1182)
	(2462,1187)(2468,1193)(2473,1202)
	(2475,1211)(2473,1219)(2470,1226)
	(2466,1233)(2461,1238)(2457,1242)(2455,1244)
\path(3026,1391)(3025,1390)(3020,1387)
	(3012,1381)(3003,1374)(2997,1364)
	(2994,1352)(2996,1342)(2999,1335)
	(3001,1330)(3004,1326)(3007,1322)
	(3011,1318)(3018,1314)(3027,1311)
	(3037,1311)(3044,1313)(3050,1316)
	(3054,1318)(3058,1321)(3063,1325)
	(3068,1331)(3073,1340)(3074,1350)
	(3073,1358)(3069,1366)(3065,1374)
	(3060,1380)(3056,1384)(3054,1387)
\path(3305,1386)(3302,1385)(3298,1383)
	(3293,1380)(3287,1376)(3282,1371)
	(3279,1364)(3278,1356)(3280,1347)
	(3283,1339)(3286,1334)(3288,1329)
	(3291,1325)(3295,1320)(3302,1316)
	(3311,1313)(3321,1314)(3329,1316)
	(3336,1319)(3341,1323)(3345,1326)
	(3350,1331)(3355,1338)(3359,1347)
	(3357,1360)(3351,1371)(3341,1380)
	(3332,1388)(3326,1392)(3325,1393)
\put(379,585){\ellipse{82}{82}}
\put(911,586){\ellipse{82}{82}}
\put(2622,586){\ellipse{82}{82}}
\path(401.915,1983.340)(379.000,2075.000)(356.085,1983.340)
\path(379,2075)(379,12)
\path(150,585)(3359,585)
\path(3267.340,562.085)(3359.000,585.000)(3267.340,607.915)
\path(2216,2154)(3606,2154)(3606,1036)
	(2216,1036)(2216,2154)
\path(2731,1037)(2673,910)
\path(2671.195,961.208)(2673.000,910.000)(2712.883,942.169)
\path(387,627)(387,628)(387,629)
	(387,630)(387,632)(387,635)
	(387,638)(388,641)(390,645)
	(392,650)(395,654)(400,659)
	(406,664)(413,669)(423,674)
	(435,679)(450,684)(467,689)
	(489,694)(514,699)(544,704)
	(580,709)(621,714)(667,719)
	(719,723)(764,726)(811,729)
	(857,731)(902,733)(944,735)
	(984,736)(1021,737)(1054,738)
	(1085,738)(1113,738)(1138,738)
	(1160,738)(1181,738)(1201,737)
	(1220,737)(1238,736)(1257,736)
	(1276,735)(1296,734)(1317,734)
	(1341,733)(1367,733)(1396,733)
	(1428,732)(1464,732)(1503,732)
	(1546,732)(1593,732)(1642,733)
	(1694,733)(1747,734)(1800,734)
	(1867,735)(1929,735)(1982,735)
	(2028,735)(2065,735)(2094,735)
	(2116,735)(2133,735)(2144,734)
	(2152,734)(2157,733)(2161,732)
	(2164,732)(2168,731)(2173,731)
	(2181,731)(2192,731)(2208,731)
	(2228,731)(2252,731)(2282,732)
	(2315,734)(2351,735)(2388,737)
	(2439,740)(2484,744)(2520,748)
	(2551,752)(2575,756)(2596,760)
	(2613,764)(2628,768)(2639,772)
	(2649,775)(2657,778)(2662,781)
	(2666,783)(2668,784)(2669,784)
\path(915,627)(917,627)(920,628)
	(926,629)(935,631)(947,633)
	(961,636)(978,638)(997,641)
	(1019,644)(1045,646)(1077,649)
	(1114,652)(1158,654)(1194,656)
	(1227,657)(1255,658)(1277,658)
	(1293,658)(1303,658)(1309,657)
	(1312,657)(1314,656)(1317,656)
	(1321,655)(1330,655)(1345,655)
	(1367,655)(1400,655)(1442,656)
	(1495,657)(1556,658)(1601,659)
	(1646,660)(1689,661)(1729,662)
	(1765,663)(1797,663)(1826,663)
	(1850,663)(1872,663)(1891,663)
	(1907,662)(1923,662)(1937,661)
	(1952,661)(1967,661)(1984,661)
	(2003,661)(2024,661)(2048,662)
	(2076,663)(2108,664)(2144,667)
	(2183,669)(2225,673)(2268,677)
	(2312,681)(2371,688)(2423,696)
	(2467,704)(2504,711)(2535,719)
	(2561,727)(2583,734)(2601,742)
	(2617,749)(2631,756)(2642,762)
	(2651,768)(2658,773)(2663,776)
	(2666,779)(2668,780)(2669,781)
\path(2628,619)(2672,801)
\put(395,428){\makebox(0,0)[lb]{\smash{{{\SetFigFont{10}{12.0}{\familydefault}{\mddefault}{\updefault}$\gamma_3$}}}}}
\put(887,421){\makebox(0,0)[lb]{\smash{{{\SetFigFont{10}{12.0}{\familydefault}{\mddefault}{\updefault}$\gamma_2$}}}}}
\put(833,898){\makebox(0,0)[lb]{\smash{{{\SetFigFont{10}{12.0}{\familydefault}{\mddefault}{\updefault}$\frac{4}{27}$}}}}}
\put(3400,544){\makebox(0,0)[lb]{\smash{{{\SetFigFont{10}{12.0}{\familydefault}{\mddefault}{\updefault}$Re$}}}}}
\put(273,455){\makebox(0,0)[lb]{\smash{{{\SetFigFont{10}{12.0}{\familydefault}{\mddefault}{\updefault}0}}}}}
\put(2335,1073){\makebox(0,0)[lb]{\smash{{{\SetFigFont{10}{12.0}{\familydefault}{\mddefault}{\updefault}$c$}}}}}
\put(2937,1200){\makebox(0,0)[lb]{\smash{{{\SetFigFont{10}{12.0}{\familydefault}{\mddefault}{\updefault}$b$}}}}}
\put(3248,1204){\makebox(0,0)[lb]{\smash{{{\SetFigFont{10}{12.0}{\familydefault}{\mddefault}{\updefault}$a$}}}}}
\put(406,2021){\makebox(0,0)[lb]{\smash{{{\SetFigFont{10}{12.0}{\familydefault}{\mddefault}{\updefault}$I_m$}}}}}
\put(2662,802){\makebox(0,0)[lb]{\smash{{{\SetFigFont{10}{12.0}{\familydefault}{\mddefault}{\updefault}$*_1$}}}}}
\put(2744,470){\makebox(0,0)[lb]{\smash{{{\SetFigFont{10}{12.0}{\familydefault}{\mddefault}{\updefault}$\gamma_1$}}}}}
\put(2540,409){\makebox(0,0)[lb]{\smash{{{\SetFigFont{10}{12.0}{\familydefault}{\mddefault}{\updefault}12}}}}}
\put(0,2262){\makebox(0,0)[lb]{\smash{{{\SetFigFont{12}{14.4}{\familydefault}{\mddefault}{\updefault}$\mathrm{H_{ii}}$.}}}}}

\put(2992,2002){\makebox(0,0)[lb]{\smash{{{\SetFigFont{10}{12.0}{\familydefault}{\mddefault}{\updefault}$*_{\mathrm{H_{ii}}}$}}}}}
\put(2762,2302){\makebox(0,0)[lb]{\smash{{{\SetFigFont{10}{12.0}{\familydefault}{\mddefault}{\updefault}$l_{*_1,\C}$}}}}}
\put(1262,1602){\makebox(0,0)[lb]{\smash{{{\SetFigFont{10}{12.0}{\familydefault}{\mddefault}{\updefault}$L_{\mathrm{H_{ii}},\C}$}}}}}
\end{picture}
}\qquad\qquad\\\\

\setlength{\unitlength}{0.00065333in}
\begingroup\makeatletter\ifx\SetFigFont\undefined%
\gdef\SetFigFont#1#2#3#4#5{%
  \reset@font\fontsize{#1}{#2pt}%
  \fontfamily{#3}\fontseries{#4}\fontshape{#5}%
  \selectfont}%
\fi\endgroup%
{\renewcommand{\dashlinestretch}{30}
\begin{picture}(3210,2624)(0,-10)
\put(2099,1378){\ellipse{16}{16}}
\put(2614,2363){\ellipse{62}{62}}
\put(3082,1539){\ellipse{16}{16}}
\put(2540,1404){\ellipse{16}{16}}
\path(2590,2337)(2089,1420)
\path(2593,2336)(2125,1418)
\path(2182,1606)(2176,1569)(2209,1587)
\path(2190,1557)(2216,1580)(2221,1545)
\path(2614,2326)(2568,1446)
\path(2635,2336)(3087,1580)
\path(2632,2336)(3069,1572)
\path(2594,2337)(2632,2386)
\path(2595,2385)(2632,2337)
\path(2074,1322)(2098,1334)(2079,1352)
\path(3002,1645)(3036,1620)(3040,1654)
\path(3054,1611)(3052,1644)(3085,1627)
\path(3072,1485)(3096,1497)(3077,1514)
\path(2521,1349)(2546,1361)(2526,1379)
\path(2608,2336)(2607,2335)(2608,2333)(2534,1450)
\path(2533,1583)(2542,1552)(2561,1578)
\path(2567,1541)(2576,1569)(2596,1546)
\path(2089,1423)(2088,1422)(2083,1419)
	(2075,1412)(2066,1404)(2059,1393)
	(2057,1380)(2058,1372)(2061,1364)
	(2063,1359)(2066,1355)(2068,1351)
	(2070,1348)(2074,1345)(2078,1341)
	(2085,1338)(2093,1336)(2102,1337)
	(2110,1339)(2116,1341)(2120,1344)
	(2124,1347)(2128,1349)(2132,1353)
	(2137,1358)(2142,1364)(2146,1372)
	(2148,1382)(2147,1391)(2143,1399)
	(2138,1406)(2133,1411)(2128,1416)(2126,1418)
\path(3067,1576)(3064,1575)(3059,1572)
	(3053,1569)(3046,1565)(3041,1559)
	(3037,1552)(3036,1543)(3038,1533)
	(3042,1525)(3045,1519)(3048,1514)
	(3051,1509)(3055,1504)(3063,1499)
	(3073,1496)(3082,1496)(3090,1498)
	(3096,1501)(3101,1503)(3106,1506)
	(3110,1509)(3114,1512)(3119,1517)
	(3123,1524)(3126,1532)(3124,1546)
	(3117,1559)(3107,1570)(3096,1578)
	(3090,1583)(3089,1584)
\path(2535,1449)(2534,1448)(2529,1445)
	(2520,1439)(2511,1431)(2504,1421)
	(2501,1408)(2503,1398)(2505,1389)
	(2508,1383)(2511,1379)(2514,1375)
	(2518,1370)(2525,1365)(2535,1362)
	(2546,1362)(2554,1364)(2561,1367)
	(2566,1370)(2570,1373)(2575,1378)
	(2582,1385)(2587,1394)(2589,1405)
	(2587,1414)(2584,1423)(2579,1431)
	(2574,1438)(2569,1443)(2567,1446)
\put(2014,1205){\makebox(0,0)[lb]{\smash{{{\SetFigFont{10}{12.0}{\familydefault}{\mddefault}{\updefault}$c$}}}}}
\put(2995,1357){\makebox(0,0)[lb]{\smash{{{\SetFigFont{10}{12.0}{\familydefault}{\mddefault}{\updefault}$a$}}}}}
\put(2433,1235){\makebox(0,0)[lb]{\smash{{{\SetFigFont{10}{12.0}{\familydefault}{\mddefault}{\updefault}$b$}}}}}
\put(969.027,1046.313){\arc{126.229}{0.5748}{1.7303}}
\put(1057.026,992.085){\arc{29.657}{1.0811}{3.8740}}
\put(542,934){\ellipse{92}{92}}
\put(1429,936){\ellipse{92}{92}}
\path(2707.480,904.620)(2809.000,930.000)(2707.480,955.380)
\path(2809,930)(246,930)
\path(1459.380,2195.480)(1434.000,2297.000)(1408.620,2195.480)
\path(1434,2297)(1434,12)
\path(2006,1768)(1845,1675)
\path(1876.259,1722.367)(1845.000,1675.000)(1901.649,1678.413)
\path(2006,2449)(3198,2449)(3198,1083)
	(2006,1083)(2006,2449)
\path(583,960)(963,985)
\path(1435,1628)(1019,1010)
\path(1379,948)(1065,977)
\path(1044,997)(1437,1628)
\put(1468,984){\makebox(0,0)[lb]{\smash{{{\SetFigFont{10}{12.0}{\familydefault}{\mddefault}{\updefault}0}}}}}
\put(504,760){\makebox(0,0)[lb]{\smash{{{\SetFigFont{10}{12.0}{\familydefault}{\mddefault}{\updefault}$\gamma_2$}}}}}
\put(1296,764){\makebox(0,0)[lb]{\smash{{{\SetFigFont{10}{12.0}{\familydefault}{\mddefault}{\updefault}$\gamma_1$}}}}}

\put(2706,2300){\makebox(0,0)[lb]{\smash{{{\SetFigFont{10}{12.0}{\familydefault}{\mddefault}{\updefault}$*_{\mathrm{H_{v}}}$}}}}}
\put(2606,2600){\makebox(0,0)[lb]{\smash{{{\SetFigFont{10}{12.0}{\familydefault}{\mddefault}{\updefault}$l_{*_1,\C}$}}}}}
\put(606,1790){\makebox(0,0)[lb]{\smash{{{\SetFigFont{10}{12.0}{\familydefault}{\mddefault}{\updefault}$L_{\mathrm{H_v},\C}$}}}}}
\put(1392,2318){\makebox(0,0)[lb]{\smash{{{\SetFigFont{10}{12.0}{\familydefault}{\mddefault}{\updefault}$I_m$}}}}}
\put(1407,1579){\makebox(0,0)[lb]{\smash{{{\SetFigFont{10}{12.0}{\familydefault}{\mddefault}{\updefault}$*_1$}}}}}
\put(2838,884){\makebox(0,0)[lb]{\smash{{{\SetFigFont{10}{12.0}{\familydefault}{\mddefault}{\updefault}$Re$}}}}}
\put(0,2453){\makebox(0,0)[lb]{\smash{{{\SetFigFont{12}{14.4}{\familydefault}{\mddefault}{\updefault}$\mathrm{H_v}$.}}}}}
\end{picture}
}
\setlength{\unitlength}{0.00065333in}
\begingroup\makeatletter\ifx\SetFigFont\undefined%
\gdef\SetFigFont#1#2#3#4#5{%
  \reset@font\fontsize{#1}{#2pt}%
  \fontfamily{#3}\fontseries{#4}\fontshape{#5}%
  \selectfont}%
\fi\endgroup%
{\renewcommand{\dashlinestretch}{30}
\begin{picture}(3380,2609)(0,-10)
\put(1704,988){\ellipse{90}{88}}
\put(1704,363){\ellipse{90}{88}}
\put(1959,708){\ellipse{90}{88}}
\put(1704,358){\ellipse{24}{24}}
\put(1956,708){\ellipse{24}{24}}
\put(1700,991){\ellipse{24}{24}}
\put(729,2052){\ellipse{54}{54}}
\put(474,1255){\ellipse{12}{12}}
\put(695,1156){\ellipse{12}{12}}
\put(988,1284){\ellipse{12}{12}}
\path(107,708)(3363,708)
\path(3264.220,683.305)(3363.000,708.000)(3264.220,732.695)
\path(1983.695,2135.220)(1959.000,2234.000)(1934.305,2135.220)
\path(1959,2234)(1959,12)
\path(646,735)(1671,967)
\path(1918,687)(641,698)
\path(639,691)(1663,363)
\path(601,1016)(613,839)
\path(585.021,886.606)(613.000,839.000)(634.298,889.947)
\path(709,2031)(465,1272)
\path(710,2029)(502,1280)
\path(729,2020)(726,1197)
\path(727,2025)(727,2024)(727,2024)(698,1197)
\path(749,2029)(996,1313)
\path(745,2029)(971,1296)
\path(711,2031)(745,2073)
\path(712,2072)(745,2031)
\path(691,1121)(713,1130)(694,1145)
\path(691,1366)(703,1347)(715,1362)
\path(716,1347)(724,1371)(740,1352)
\path(1192,2095)(1192,1029)(372,1029)
	(372,2095)(1192,2095)
\path(457,1203)(477,1214)(462,1229)
\path(508,1375)(530,1395)(535,1366)
\path(503,1415)(498,1383)(527,1398)
\path(901,1424)(927,1401)(933,1431)
\path(955,1403)(955,1428)(984,1411)
\path(975,1241)(996,1251)(979,1266)
\path(702,1197)(701,1196)(697,1194)
	(689,1188)(681,1181)(675,1173)
	(672,1161)(673,1152)(676,1145)
	(678,1140)(681,1136)(683,1133)
	(687,1129)(693,1125)(702,1122)
	(711,1122)(719,1124)(724,1127)
	(729,1129)(733,1132)(737,1136)
	(742,1142)(747,1150)(748,1162)
	(745,1172)(739,1182)(733,1189)
	(730,1193)(729,1194)
\path(462,1294)(459,1292)(455,1289)
	(449,1285)(444,1280)(439,1274)
	(435,1266)(434,1257)(435,1250)
	(437,1243)(440,1239)(442,1235)
	(443,1232)(446,1229)(448,1226)
	(452,1223)(458,1220)(465,1218)
	(474,1219)(482,1221)(487,1224)
	(492,1227)(495,1231)(500,1235)
	(505,1242)(510,1250)(512,1259)
	(510,1266)(507,1273)(503,1279)
	(499,1284)(495,1287)(493,1289)
\path(973,1309)(970,1308)(966,1306)
	(961,1303)(955,1299)(950,1295)
	(947,1289)(946,1281)(948,1273)
	(951,1266)(954,1260)(956,1256)
	(959,1252)(963,1248)(969,1244)
	(978,1241)(988,1242)(996,1244)
	(1002,1247)(1006,1250)(1011,1253)
	(1016,1258)(1020,1264)(1024,1273)
	(1023,1285)(1016,1296)(1007,1305)
	(998,1312)(993,1316)(992,1317)
\put(1865,2259){\makebox(0,0)[lb]{\smash{{{\SetFigFont{10}{12.0}{\familydefault}{\mddefault}{\updefault}$I_m$}}}}}
\put(1685,1056){\makebox(0,0)[lb]{\smash{{{\SetFigFont{10}{12.0}{\familydefault}{\mddefault}{\updefault}$\gamma_1$}}}}}
\put(1782,305){\makebox(0,0)[lb]{\smash{{{\SetFigFont{10}{12.0}{\familydefault}{\mddefault}{\updefault}$\gamma_3$}}}}}

\put(788,1922){\makebox(0,0)[lb]{\smash{{{\SetFigFont{10}{12.0}{\familydefault}{\mddefault}{\updefault}$*_{\mathrm{H_{vi}}}$}}}}}
\put(608,2222){\makebox(0,0)[lb]{\smash{{{\SetFigFont{10}{12.0}{\familydefault}{\mddefault}{\updefault}$l_{*_1,\C}$}}}}}
\put(3380,654){\makebox(0,0)[lb]{\smash{{{\SetFigFont{10}{12.0}{\familydefault}{\mddefault}{\updefault}$Re$}}}}}
\put(588,525){\makebox(0,0)[lb]{\smash{{{\SetFigFont{10}{12.0}{\familydefault}{\mddefault}{\updefault}$*_{1}$}}}}}

\put(2551,1701){\makebox(0,0)[lb]{\smash{{{\SetFigFont{10}{12.0}{\familydefault}{\mddefault}{\updefault}$L_{\mathrm{H_{vi}},\C}$}}}}}
\put(1774,766){\makebox(0,0)[lb]{\smash{{{\SetFigFont{10}{12.0}{\familydefault}{\mddefault}{\updefault}$\gamma_2$}}}}}
\put(2022,766){\makebox(0,0)[lb]{\smash{{{\SetFigFont{10}{12.0}{\familydefault}{\mddefault}{\updefault}0}}}}}
\put(760,1112){\makebox(0,0)[lb]{\smash{{{\SetFigFont{9}{10.8}{\familydefault}{\mddefault}{\updefault}$b$}}}}}
\put(990,1347){\makebox(0,0)[lb]{\smash{{{\SetFigFont{9}{10.8}{\familydefault}{\mddefault}{\updefault}$a$}}}}}
\put(400,1110){\makebox(0,0)[lb]{\smash{{{\SetFigFont{9}{10.8}{\familydefault}{\mddefault}{\updefault}$c$}}}}}
\put(0,2438){\makebox(0,0)[lb]{\smash{{{\SetFigFont{12}{14.4}{\familydefault}{\mddefault}{\updefault}$\mathrm{H_{vi}}$.}}}}}
\end{picture}
}\qquad\qquad \\\\

\centerline{
\setlength{\unitlength}{0.00068333in}
\begingroup\makeatletter\ifx\SetFigFont\undefined%
\gdef\SetFigFont#1#2#3#4#5{%
  \reset@font\fontsize{#1}{#2pt}%
  \fontfamily{#3}\fontseries{#4}\fontshape{#5}%
  \selectfont}%
\fi\endgroup%
{\renewcommand{\dashlinestretch}{30}
\begin{picture}(3215,2624)(0,-10)
\put(1210.500,1033.900){\arc{517.251}{0.9209}{1.1035}}
\put(1399.712,800.035){\arc{30.759}{1.1496}{3.6116}}
\put(305,736){\ellipse{90}{90}}
\put(954,741){\ellipse{90}{90}}
\put(1725,734){\ellipse{90}{90}}
\put(2118,1292){\ellipse{14}{14}}
\put(2821,1483){\ellipse{14}{14}}
\put(3051,1494){\ellipse{14}{14}}
\put(2628,2267){\ellipse{72}{72}}
\path(1746.125,2172.490)(1721.000,2273.000)(1695.875,2172.490)
\path(1721,2273)(1721,12)
\path(2800.462,698.986)(2901.000,724.000)(2800.518,749.236)
\path(2901,724)(175,727)
\path(334,778)(1347,857)(1652,977)
\path(998,763)(1320,803)
\path(1370,828)(1649,970)
\path(1698,755)(1388,783)
\path(1995,895)(1647,979)(1385,812)
\path(2307,1067)(2147,975)
\path(2178.038,1021.829)(2147.000,975.000)(2203.086,978.267)
\path(2605,2242)(2109,1335)
\path(2607,2240)(2145,1332)
\path(2201,1518)(2195,1482)(2228,1500)
\path(2209,1469)(2234,1493)(2240,1459)
\path(2650,2240)(3063,1538)
\path(2647,2240)(3046,1519)
\path(2608,2242)(2647,2291)
\path(2609,2290)(2647,2242)
\path(2094,1238)(2117,1249)(2099,1267)
\path(2623,2240)(2622,2239)(2623,2238)(2805,1517)
\path(2629,2239)(2837,1519)
\path(2022,2349)(3203,2349)(3203,1063)
	(2022,1063)(2022,2349)
\path(2773,1604)(2782,1573)(2800,1599)
\path(2811,1564)(2820,1591)(2839,1569)
\path(2988,1581)(3022,1555)(3026,1589)
\path(3036,1558)(3034,1591)(3067,1573)
\path(3044,1443)(3069,1455)(3050,1472)
\path(2816,1427)(2841,1437)(2821,1455)
\path(2109,1338)(2108,1337)(2103,1334)
	(2095,1328)(2086,1320)(2080,1309)
	(2077,1296)(2078,1288)(2081,1280)
	(2083,1275)(2085,1271)(2087,1267)
	(2090,1264)(2093,1260)(2097,1257)
	(2104,1253)(2112,1251)(2121,1251)
	(2128,1253)(2135,1256)(2139,1259)
	(2143,1262)(2147,1264)(2151,1268)
	(2156,1273)(2161,1279)(2165,1287)
	(2167,1297)(2166,1306)(2162,1313)
	(2157,1320)(2152,1326)(2147,1330)(2145,1332)
\path(2810,1525)(2809,1524)(2804,1521)
	(2795,1515)(2785,1507)(2778,1496)
	(2775,1483)(2777,1473)(2780,1465)
	(2783,1459)(2786,1454)(2789,1450)
	(2794,1446)(2801,1441)(2811,1438)
	(2820,1438)(2827,1439)(2833,1441)
	(2837,1443)(2841,1445)(2844,1448)
	(2848,1451)(2852,1455)(2857,1461)
	(2861,1469)(2863,1483)(2859,1495)
	(2853,1506)(2847,1515)(2843,1520)(2842,1521)
\path(3042,1530)(3039,1529)(3035,1526)
	(3029,1523)(3023,1518)(3017,1513)
	(3014,1506)(3013,1497)(3015,1487)
	(3019,1480)(3022,1474)(3024,1469)
	(3027,1464)(3032,1460)(3039,1455)
	(3049,1452)(3060,1453)(3069,1455)
	(3076,1458)(3081,1461)(3086,1465)
	(3091,1470)(3097,1477)(3101,1487)
	(3099,1501)(3092,1514)(3082,1524)
	(3072,1533)(3066,1538)(3065,1539)
\put(1595,2302){\makebox(0,0)[lb]{\smash{{{\SetFigFont{10}{12.0}{\familydefault}{\mddefault}{\updefault}$I_m$}}}}}
\put(184,535){\makebox(0,0)[lb]{\smash{{{\SetFigFont{10}{12.0}{\familydefault}{\mddefault}{\updefault}-3}}}}}
\put(397,544){\makebox(0,0)[lb]{\smash{{{\SetFigFont{10}{12.0}{\familydefault}{\mddefault}{\updefault}$\gamma_3$}}}}}
\put(1553,539){\makebox(0,0)[lb]{\smash{{{\SetFigFont{10}{12.0}{\familydefault}{\mddefault}{\updefault}$\gamma_1$}}}}}
\put(1978,814){\makebox(0,0)[lb]{\smash{{{\SetFigFont{10}{12.0}{\familydefault}{\mddefault}{\updefault}$*_1$}}}}}

\put(2700,2219){\makebox(0,0)[lb]{\smash{{{\SetFigFont{10}{12.0}{\familydefault}{\mddefault}{\updefault}$*_{\mathrm{H_{viii}}}$}}}}}
\put(3350,1719){\makebox(0,0)[lb]{\smash{{{\SetFigFont{10}{12.0}{\familydefault}{\mddefault}{\updefault}$l_{*_1,\C}$}}}}}
\put(750,1719){\makebox(0,0)[lb]{\smash{{{\SetFigFont{10}{12.0}{\familydefault}{\mddefault}{\updefault}$L_{\mathrm{H_{viii}},\C}$}}}}}

\put(2034,1121){\makebox(0,0)[lb]{\smash{{{\SetFigFont{10}{12.0}{\familydefault}{\mddefault}{\updefault}$c$}}}}}
\put(2914,673){\makebox(0,0)[lb]{\smash{{{\SetFigFont{10}{12.0}{\familydefault}{\mddefault}{\updefault}$Re$}}}}}
\put(1787,586){\makebox(0,0)[lb]{\smash{{{\SetFigFont{10}{12.0}{\familydefault}{\mddefault}{\updefault}0}}}}}
\put(2717,1314){\makebox(0,0)[lb]{\smash{{{\SetFigFont{10}{12.0}{\familydefault}{\mddefault}{\updefault}$b$}}}}}
\put(3000,1335){\makebox(0,0)[lb]{\smash{{{\SetFigFont{10}{12.0}{\familydefault}{\mddefault}{\updefault}$a$}}}}}
\put(1142,544){\makebox(0,0)[lb]{\smash{{{\SetFigFont{10}{12.0}{\familydefault}{\mddefault}{\updefault}$\gamma_2$}}}}}
\put(703,547){\makebox(0,0)[lb]{\smash{{{\SetFigFont{10}{12.0}{\familydefault}{\mddefault}{\updefault}$-\frac{32}{27}$}}}}}
\put(0,2453){\makebox(0,0)[lb]{\smash{{{\SetFigFont{12}{14.4}{\familydefault}{\mddefault}{\updefault}$\mathrm{H_{viii}}$.}}}}}
\end{picture}
} 
}
}
\begin{center}
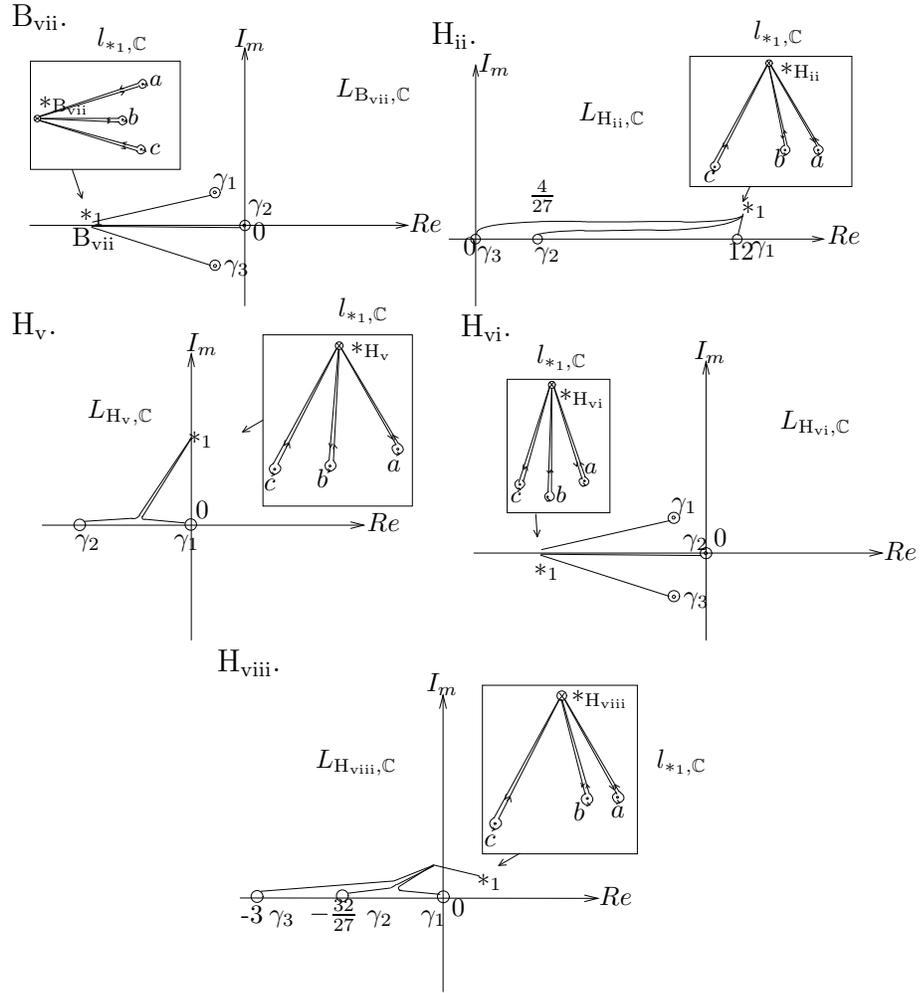
Figure 4. Complex line $B_{X,\C}$ and complex pencil $l_{*_1,\C}$.
\end{center}

For each type $X$  of 17 polynomials, applying  Zariski-van Kampen method to the generators
$a,\ b$ and $c$ explained above, we obtain the following presentations 
of the fundamental group 
$\pi_1(S_{X} \setminus D_{X}, *_{X})\cong 
\pi_1(H_{X} \setminus C_{X}, *_{X})$.

\newpage
\bf Table 1.
}
\vspace{-0.2cm}
\[
\begin{array}{lll}
\pi_1(S_{\mathrm{A_i}} \setminus D_{\mathrm{A_i}}, *_{\mathrm{A_i}}) 
&\!\! \cong 
\pi_1(H_{\mathrm{A_i}} \setminus C_{\mathrm{A_i}}, *_{\mathrm{A_i}})
\cong
\biggl{\langle}
a, b, c 
\biggl{|}
\begin{array}{cc}ab=ba, \\
 bcb=cbc, \\
  aca=cac
  \end{array}\
\biggl{\rangle}. \\

\pi_1(S_{\mathrm{A_{ii}}}\setminus D_{\mathrm{A_{ii}}},*_{\mathrm{A_{ii}}}) 
&\!\! \cong 
\pi_1(H_{\mathrm{A_{ii}}}\setminus C_{\mathrm{A_{ii}}},*_{\mathrm{A_{ii}}})
\cong
\Biggl{\langle}
a, b, c 
\Biggl{|}
\begin{array}{cc}ababab=bababa, \\
aba=bab, \\
 b=c
  \end{array}\
\Biggl{\rangle}.\\

\pi_1(S_{\mathrm{B_{i}}} \setminus D_{\mathrm{B_{i}}}, *_{\mathrm{B_i}}) 
&\!\! \cong 
\pi_1(H_{\mathrm{B_{i}}} \setminus C_{\mathrm{B_{i}}}, *_{\mathrm{B_i}}) 
\cong
\Biggl{\langle}
a, b, c 
\Biggl{|}
\begin{array}{cc}abab=baba, \\
 bc=cb, \\
 aca=cac, \\
 cbac=baca 
  \end{array}\
\Biggl{\rangle}.\\
 
\pi_1(S_{\mathrm{B_{ii}}}\setminus D_{\mathrm{B_{ii}}},*_{\mathrm{B_{ii}}}) 
&\!\! \cong 
\pi_1(H_{\mathrm{B_{ii}}}\setminus C_{\mathrm{B_{ii}}},*_{\mathrm{B_{ii}}})
\cong
\Biggl{\langle}
a, b, c 
\Biggl{|}
\begin{array}{cc}ababab=bababa,\\
 bc=ab, \\
 ac=ca
  \end{array}\
\Biggl{\rangle}. \\
 
\pi_1(S_{\mathrm{B_{iii}}}\!\setminus\! D_{\mathrm{B_{iii}}},*_{\mathrm{B_{iii}}})
&\!\!\! \cong 
\pi_1(H_{\mathrm{B_{iii}}}\!\setminus\! C_{\mathrm{B_{iii}}},*_{\mathrm{B_{iii}}})
\cong\!\!
\Biggl{\langle}\!\!
a, b, c 
\Biggl{|}\!\!
\begin{array}{cc}a=b, \\
 a=cbab^{-1}c^{-1}, \\
 b=cbacbc^{-1}a^{-1}b^{-1}c^{-1}, \\
 c\!=\!cbacbcb^{-\!1}\!c^{-\!1}\!a^{-\!1}\!b^{-\!1}\!c^{-\!1}
  \end{array}\!\!
\Biggl{\rangle}.\\

\pi_1(S_{\mathrm{B_{iv}}}\setminus D_{\mathrm{B_{iv}}},*_{\mathrm{B_{iv}}}) 
&\!\!\cong 
\pi_1(H_{\mathrm{B_{iv}}}\setminus C_{\mathrm{B_{iv}}},*_{\mathrm{B_{iv}}})
\cong
\Biggl{\langle}
a,b,c
\Biggl{|}
\begin{array}{cc}acb=cba,\\
 bcba=cbac,\\
 cbac=bacb,\\
 ab=ba
  \end{array}\
\Biggl{\rangle}.\\
 
\pi_1(S_{\mathrm{B_{v}}}\setminus D_{\mathrm{B_{v}}},*_{\mathrm{B_v}}) 
&\!\!\cong 
\pi_1(H_{\mathrm{B_{v}}}\setminus C_{\mathrm{B_{v}}},*_{\mathrm{B_v}})
\cong
\biggl{\langle}
a, b, c 
\biggl{|}
\begin{array}{cc}
 a=b=c
  \end{array}\
\biggl{\rangle}.\\
 
\pi_1(S_{\mathrm{B_{vi}}}\setminus D_{\mathrm{B_{vi}}},*_{\mathrm{B_{vi}}}) 
&\!\!\cong 
\pi_1(H_{\mathrm{B_{vi}}}\setminus C_{\mathrm{B_{vi}}},*_{\mathrm{B_{vi}}})
\cong
\Biggl{\langle}
a, b, c 
\Biggl{|}
\begin{array}{cc}aba=bab,\\
 aca=bac,\\
 acaca=cacac
  \end{array}\
\Biggl{\rangle}.
\end{array}
\vspace{-0.4cm}
\]
\leftline{
$\pi_1(S_{\mathrm{B_{vii}}}\setminus D_{\mathrm{B_{vii}}},*_{\mathrm{B_{vii}}}) 
\cong 
\pi_1(H_{\mathrm{B_{vii}}}\setminus C_{\mathrm{B_{vii}}},
*_{\mathrm{B_{vii}}})
$ 
}
\vspace{-0.2cm}
\[
\cong
\Biggl{\langle}
a, b, c 
\Biggl{|}
\begin{array}{cc}
 a=b^{-1}cbab^{-1}cbab^{-1}cbab^{-1}cba^{-1}b^{-1}c^{-1}ba^{-1}b^{-1}c^{-1}ba^{-1}b^{-1}c^{-1}b,\\
 c=bab^{-1}cbab^{-1}cbab^{-1}cbab^{-1}c^{-1}ba^{-1}b^{-1}c^{-1}ba^{-1}b^{-1}c^{-1}ba^{-1}b^{-1},\\
 a=ba^{-1}b^{-1}c^{-1}bab^{-1}cbab^{-1},\\
 cba=bab,
 cba=bcb,
 cba=bab^{-1}c^{-1}b^{-1}cbcb
  \end{array}\
\Biggl{\rangle}.
\vspace{-0.3cm}
\]
\[
\begin{array}{lll} 
\pi_1(S_{\mathrm{H_i}}\setminus D_{\mathrm{H_i}},*_{\mathrm{H_i}}) 
&\!\!\cong 
\pi_1(H_{\mathrm{H_i}}\setminus C_{\mathrm{H_i}},*_{\mathrm{H_i}})
\cong
\Biggl{\langle}
a,b,c
\Biggl{|}
\begin{array}{cc}ababa=babab,\\
 bc=cb,\\
 aca=cac
  \end{array}\
\Biggl{\rangle}.\\
 
\vspace{-0.2cm}
\pi_1(S_{\mathrm{H_{ii}}}\setminus D_{\mathrm{H_{ii}}},*_{\mathrm{H_{ii}}}) 
&\!\!\cong 
\pi_1(H_{\mathrm{H_{ii}}}\setminus C_{\mathrm{H_{ii}}},*_{\mathrm{H_{ii}}})
\cong
\Biggl{\langle}
a,b,c
\Biggl{|}
\begin{array}{cc}
 abab=baba,\\
 aca=bac,\\
 acaca=cacac
  \end{array}\
\Biggl{\rangle}.\\
 
\pi_1(S_{\mathrm{H_{iii}}}\setminus D_{\mathrm{H_{iii}}},*_{\mathrm{H_{iii}}}) 
&\!\!\cong 
\pi_1(H_{\mathrm{H_{iii}}}\setminus C_{\mathrm{H_{iii}}},*_{\mathrm{H_{iii}}})
\cong
\Biggl{\langle}
a,b,c
\Biggl{|}
\begin{array}{cc}aba=bab,\\
 bcba=cbac,\\
 cba=acb
  \end{array}\
\Biggl{\rangle}.\\
\vspace{-0.3cm}

\pi_1(S_{\mathrm{H_{iv}}}\setminus D_{\mathrm{H_{iv}}},*_{\mathrm{H_{iv}}}) 
&\!\!\cong 
\pi_1(H_{\mathrm{H_{iv}}}\setminus C_{\mathrm{H_{iv}}},*_{\mathrm{H_{iv}}})
\cong
\biggl{\langle}
a, b, c 
\biggl{|}
\begin{array}{cc}a=b=c
  \end{array}\
\biggl{\rangle}.\\
\end{array}
\vspace{-0.2cm}
\]
\[
\begin{array}{lll} 
\pi_1(S_{\mathrm{H_{v}}}\setminus D_{\mathrm{H_{v}}},*_{\mathrm{H_v}}) 
&\!\!\cong 
\pi_1(H_{\mathrm{H_{v}}}\setminus C_{\mathrm{H_{v}}},*_{\mathrm{H_v}})
\cong
\Biggl{\langle}
a,b,c
\Biggl{|}
\begin{array}{cc}acba=cbac,\\
bcbac=cbacb,\\
bacb=cbac,\\
bc=cb
  \end{array}\
\Biggl{\rangle}.\\
\end{array}
\]
\[ 
\begin{array}{lll} 
\pi_1(S_{\mathrm{H_{vi}}}\setminus D_{\mathrm{H_{vi}}},*_{\mathrm{H_{vi}}}) 
&\!\!\cong 
\pi_1(H_{\mathrm{H_{vi}}}\setminus C_{\mathrm{H_{vi}}},*_{\mathrm{H_{vi}}})
\cong\!
\Biggl{\langle}\!\!
a,b,c
\Biggl{|}\!\!
\begin{array}{cc}abababab=babababa,\\
 ba=cb,\\
 ac=ba
  \end{array}\!\!
\Biggl{\rangle}.\\
 
\pi_1(S_{\mathrm{H_{vii}}}\!\setminus\! D_{\mathrm{H_{vii}}},*_{\mathrm{H_{vii}}}) 
\!\!&\!\!\cong\! 
\pi_1(H_{\mathrm{H_{vii}}}\!\!\setminus\!\! C_{\mathrm{H_{vii}}},*_{\mathrm{H_{vii}}})
\!\!\cong\!
\Biggl{\langle}\!\!
a,b,c
\Biggl{|}\!\!\!
\begin{array}{cc}
 a=cbaca^{-1}b^{-1}c^{-1},\\
 b=cbacbc^{-1}a^{-1}b^{-1}c^{-1},\\
 c\!=\!cbacbab^{-\!1}\!c^{-\!1}\!a^{-\!1}\!b^{-\!1}\!c^{-\!1}\!\!,\\
 b=c
  \end{array}\!\!\!
\Biggl{\rangle}.\\
 
\pi_1(S_{\mathrm{H_{viii}}}\!\setminus\! D_{\mathrm{H_{viii}}},*_{\mathrm{H_{viii}}}) 
\!&\!\!\!\cong 
\!\pi_1(H_{\mathrm{H_{viii}}}\!\setminus\! C_{\mathrm{H_{viii}}},*_{\mathrm{H_{viii}}})
\cong
\Biggl{\langle}
\!a,b,c
\Biggl{|}\!
\begin{array}{cc}
 abababa=bababab,\\
 ab=bc,\\
 ac=ca
  \end{array}\!\!\!\!
\Biggl{\rangle}.
\end{array}
\]
\vspace{-0.3cm}

\section{Positive Homogeneous Presentation}

In the present section, we rewrite the 
presentations of the fundamental groups in section 3 to a positive homogeneous form.
We, first,  prepare some terminology.

\begin{definition}
1. Let 
$G = \langle L \mid R\rangle$
be a presentation of a group $G$, where $L$ is the set of generators 
(called alphabets) and 
$R$ is the set of relations. We call that the presentation is {\it positive 
homogeneous}, if  $R$ consists of relations of the form $R_i\!=\!S_i$ 
where 
$R_i$ and $S_i$ are positive words in the letters $L$ 
(i.e.~words consisting of only 
non-negative powers of the letters in $L$) of the same length.

2. If a positive homogeneous presentation $\langle L \mid R\rangle$ 
of a group $G$ is given,
then we associate a monoid $M$ defined as the quotient 
of free monoid $L^*$ generated by $L$ by the equivalence relation $\simeq$
defined as follows: 

1) two words $U$ and $V$ in $L^*$ are called {\it elementarily 
equivalent} if either $U=V$ or $V$ is obtained from $U$ by substituting 
a substring $R_i$ of $U$ by $S_i$ where $R_i\!=\!S_i$ is a relation of $R$ 
($S_i=R_i$ is also a relation if $R_i=S_i$ is a relation), 

2) two words $U$ and $V$ in $L^*$ are called {\it equivalent}, denoted by $U\simeq V$, if there exists 
a sequence $U\!=\!W_0, W_1,\cdots, W_n\!=\!V $ of words in $L^*$ for $n\!\in\!\Z_{\ge0}$
such that $W_i$ is elementarily equivalent to $W_{i-1}$ for $i=1,\cdots,n$.

3. The natural homomorphism $M\to G$ will be called the {\it localization 
morphism}. The image of the localization homomorphism is denoted by $G^+$.

\end{definition}

\noindent
{\it Note.} 1. The monoid $G^+$ depends on the choice of the generators for the group $G$. Even if we choose the same generators for the same group $G$, the monoid $M$ depends on the choice of the relations $R$. 

2.  Due to the homogeneity of the relations, one defines a homomorphism: 

\centerline{
$l\ :\ G\ \longrightarrow\ \Z$
}

\noindent
by associating 1 to each letter in $L$. The restriction of the homomorphism on $G^+$ and its pull-back to $M$ by the localization homomorphism are called {\it length functions}. Length functions have the additivity: $l(UV)=l(U)+l(V)$ and the conicity: $l(U)=1$ implies $U=1$. 
The existence of such length functions implies that the monoids $M$ and $G^+$ are {\it atomic} (\cite[\S2]{[D-P]}).

\bigskip
\noindent
{\bf Theorem 1.} {\it The fundamental group in {\bf Table 1.}\ of type $X$ is naturally isomorphic to  the following positive homogeneously presented group $G_X$ by identifying the generators $\{a,b,c\}$ in the both groups. }

\[\begin{array}{rlll}
\mathrm{A_{i}}&:  G_{\mathrm{A_{i}}}\ :=\
\biggl{\langle}
a,b,c
\biggl{|}
\begin{array}{cc}ab=ba,\\
 bcb=cbc,\\
  aca=cac
  \end{array}\
\biggl{\rangle}\ . \\
\mathrm{A_{ii}}& :  G_{\mathrm{A_{ii}}}\ :=\
\biggl{\langle}
a,b,c
\biggl{|}
\begin{array}{cc}aba=bab,\\
 b=c
  \end{array}\
\biggl{\rangle} \ .\\

\mathrm{B_{i}} &:  G_{\mathrm{B_{i}}}\ :=\
\biggl{\langle}
a,b,c
\biggl{|}
\begin{array}{cc}abab=baba,\\
 bc=cb,\\
 aca=cac
  \end{array}\
\biggl{\rangle}\ .
 \vspace{0.0cm}\\

\mathrm{B_{ii}} &:   G_{\mathrm{B_{ii}}}\ :=\
\biggl{\langle}
a,b,c
\biggl{|}
\begin{array}{cc}cbb=bba,\\
 bc=ab,\\
 ac=ca
  \end{array}\
\biggl{\rangle}\ .\\

\mathrm{B_{iii}} &:  G_{\mathrm{B_{iii}}}\ :=\
\biggl{\langle}
a,b,c
\biggl{|}
\begin{array}{cc}a=b,\\
 ac=ca
  \end{array}\
\biggl{\rangle}\ .\\

\mathrm{B_{iv}} &:  G_{\mathrm{B_{iv}}}\ :=\
\biggl{\langle}
a,b,c
\biggl{|}
\begin{array}{cc}ab=ba,\\
 bcb=cbc,\\
 ac=ca
  \end{array}\
\biggl{\rangle}\ .\\

\mathrm{B_{v}} &:  G_{\mathrm{B_{v}}}\ :=\
\biggl{\langle}
a,b,c
\biggl{|}
\begin{array}{cc}
 a=b=c
  \end{array}\
\biggl{\rangle}\ .\\

\mathrm{B_{vi}} &:  G_{\mathrm{B_{vi}}}\!\! :=\!
\Biggl{\langle}\!\!
a,b,c
\biggl{|}
\begin{array}{lll}
 aba\!=\!bab,\ bcb\!=\!cbc,\ aca\!=\!bac,\ cab\!=\!bca,\ acb\!=\!cac,\\
abb=bbc,\ bcca=ccac,\ bbac=caab,\ cbbb=bbba,\\
\!\!acbcb\!=\!bccca, accbb\!=\!bccba, accaa\!=\!ccaac, caacc\!=\!aacca,\!\!\!\!\!\!\\
acccc=bcccb,\ bbaac=cbaab,\ caaab=abaac,\\
 a^5=b^5=c^5,\ ccbaac=accbaa
  \end{array}\
\Biggl{\rangle}\ .\\

\mathrm{B_{vii}} &:  G_{\mathrm{B_{vii}}}\ :=\
\biggl{\langle}
a,b,c
\biggl{|}
\begin{array}{cc}
 a=b=c
  \end{array}\
\biggl{\rangle}\ .\\

\mathrm{H_i} &:  G_{\mathrm{H_{i}}}\ :=\
\biggl{\langle}
a,b,c
\biggl{|}
\begin{array}{cc}ababa=babab,\\
 bc=cb,\\
 aca=cac
  \end{array}\
\biggl{\rangle}\ .
 \vspace{0.0cm}\\

\mathrm{H_{ii}} &:  G_{\mathrm{H_{ii}}}\ :=\
\biggl{\langle}
a,b,c
\biggl{|}\quad 
R_{\mathrm{H_{ii}}}\quad
\biggl{\rangle}\ \ \text{ ($R_{\mathrm{H_{ii}}}$ is given at the end of present Table)}.
 \vspace {0.0cm}\\

\mathrm{H_{iii}} &:  G_{\mathrm{H_{iii}}}\!\! :=\!
\Biggl{\langle}\!\!
a,b,c
\biggl{|}
\begin{array}{lll}
 aba\!=\!bab,\ aca\!=\!cac,\ bcb\!=\!abc,\ cba\!=\!acb,\ bca\!=\!cbc,\\
baa=aac,\ accb=ccbc,\ aabc=cbba,\ caaa=aaab,\\
\!\!bcaca\!=\!acccb, bccaa\!=\!accab,bccbb\!=\!ccbbc, cbbcc\!=\!bbccb,\!\!\!\!\!\!\!\!\!\!\\
bcccc=accca,\ aabbc=cabba,\ cbbba=babbc,\\
 a^5=b^5=c^5,\ ccabbc=bccabb
  \end{array}\
\Biggl{\rangle}\ .\\

\mathrm{H_{iv}} &:  G_{\mathrm{H_{iv}}}\ :=\
\biggl{\langle}
a,b,c
\biggl{|}
\begin{array}{cc}a=b=c
  \end{array}\
\biggl{\rangle}\ .\\

\mathrm{H_{v}}& :  G_{\mathrm{H_{v}}}\ :=\
\biggl{\langle}
a,b,c
\biggl{|}
\begin{array}{cc}a=b=c
  \end{array}\
\biggl{\rangle}\ .\\

\mathrm{H_{vi}} &:  G_{\mathrm{H_{vi}}}\ :=\
\biggl{\langle}
a,b,c
\biggl{|}
\begin{array}{cc}
 a=b=c
  \end{array}\
\biggl{\rangle}\ .\\

\mathrm{H_{vii}} &:   G_{\mathrm{H_{vii}}}\ :=\
\biggl{\langle}
a,b,c
\biggl{|}
\begin{array}{cc}
 a=b=c
  \end{array}\
\biggl{\rangle}\ .\\

\mathrm{H_{viii}} &:  G_{\mathrm{H_{viii}}}\ :=\
\biggl{\langle}
a,b,c
\biggl{|}
\begin{array}{cc}
 a=b=c
  \end{array}\
\biggl{\rangle}\ .
\end{array}
\]
 \[
R_{\mathrm{H_{ii}}}:=
\Biggl\{
\begin{array}{l}
a b a b = b a b a,
a c a = b a c,
b c b c = c b c b,
a c b = c a c,
b b c a b a = a b c c a c,\\
a b b b c a = b a a a a c,
b b b a b b = a b b a a a,
b a a a a b a = a b b b b a b,\\
b a a b b b = a a a b a a,
a b c c c = c c c a b,
b b c b a b = c c c a a c,\\
c c c b c a a = b b b c c a b,
b c c b b b = c c c b c c,
b b c c a b = c a a c c c,\\
c c a a c = b c c a a,
c c a a b = a c c a a,
c c a b a a c = a c c b c a a,\\
c a a c c a b = b c a a c c a,
a a b a a a = b b b a a b,
b b b a a a = a a a b b b,\\
a b a a a a b = b a b b b b a,
a a a b b a = b b a b b b,
b a a b b a a = a a b b a a b,\\
b a a b a a b a a = a b b a b b a b b,
a a b b a a c = b a b b c a a,
a a a b c = b c a a a,\\
a b b a a b a a c = b a b b a b b c a,
c c c a a a = a a a c c c,
c c c b b b = b b b c c c,\\
c a a c a a c = a a b c c b a,
b b b c b b = c b b c c c,
a b a c b c = c b c a b a,\\
c b b b b c b = b c c c c b c,
c a b b b c = a c c c c b,
b c c c c c a a = c b b b c a a c,\\
c c b c c c = b b b c c b,
c b c a a a b = b c c c a b a,
c a a b c b = b a c c c a,\\
b c b a a b = a a c c b a,
b a a c c b b c = c a c c a b c b,
b c c a b b = a c c a a a,\\
b a b c b a b = c a b c a c a,
c a a b b b b c b = b a c c c c c c a,
c b a a c c = b c c b a b,\\
a b c b a a = c c b a b b,
b c b b a a = c c b b a b,
c a a c a c = b a b c c a,\\
c b b a a a a c c = a c a c b b c b a,
c a a a a c c = a a b c c c a,
b c a b b c c = a a b b c b b,\\
b b c a a b c = c c c a a b b,
c b b c a a b = b c c a b b a,
b b a a b b a = a b b a a b b,\\
a b a a b c c = b b a b b c b,
b a c b c a b = c a b c a b a,
c b c a b c a = b c a a c a b,\\
c a a c c b b a = b c a b c a a b,
b a b b c b b = a a b c c b c,
b b c b b b = c c c b b c,\\
b c b b b b c = c b c c c c b,
b c c b b a b b c = a b c a b c c b a,
b b a b c b b a b = c b b a b b c c c,\\
c a b a a c c c = a b b c b b a b,
b a c a b c = c b c a b b,
b c a b a a b = a b c a b a a,\\
a a c c b c a b = b c a b a a c c,
c b a a b c c = b a c c b c a,
c c c b a a b c = b a a c c a b a,\\
b c c b a a b c = c a b a c b c a,
a b a a b c a b a = b b a a b c a b b,
c c b b a a a = a a a c c b b,\\
c c b b a a b c a = a b c a b b a a c,
b a a b c a b b a = a b a c a b a a b,
b c a a b b = a a a c c a,\\
a c c b b c c = c c a b b c b,
b b c a b b c c c = a b b a b b c c b,
b c a a c c b c = a b c a b c c a,\\
c a b a a b c c = b a b c c b c a,
b a b c c b a = c b b a b c b
\end{array}
\Biggr\}
\]
\begin{proof} Except for the types 
$\mathrm{B_{ii}},\mathrm{B_{vi}},\mathrm{H_{ii}},\mathrm{H_{iii}},\mathrm{H_{viii}}$,
the relations are obtained by elementary reductions of the Zariski-van Kampen relations, and we omit details.

Some new relations for the cases of types 
$\mathrm{B_{ii}},\mathrm{B_{vi}},\mathrm{H_{ii}},\mathrm{H_{iii}}$
are obtained by  cancelling common factors from the left or from the right of equivalent expressions of the same fundamental elements (introduced in \S6 6.1. See \S7 Definition 7.1), where these equivalent expressions of a fundamental element are obtained by the help of Hayashi's computer program (see http://www.kurims.kyoto-u.ac.jp/~saito/SI/). 
In the following, we sketch how some of them are obtained by hand calculations. In the proof,  ``the first relation, the second relation, \ldots'', 
mean ``the relation which is at the first place, the second place, 
\ldots  in Table 1. of Zariski-van Kampen relations in \S3''. 

The case for the type $\mathrm{H_{viii}}$ needs to be treated separately because its calculations are non-trivial.
Detailed verifications  are left to the reader.

\medskip
$\mathrm{B_{ii}}$: Using $ab=bc$, rewrite the LHS $ababab$ (resp.\ RHS $bababa$) 
of the first relation to $bcabbc$ (resp.\ $babbca$). Then, using 
the commutativity of $a$ and $c$, we cancel $ba$ from left and $c$ from right 
so that we obtain a new relation $cbb=bba$.

\medskip
$\mathrm{B_{vi}}$:  Using $aca=bac$, rewrite the LHS $acaca$  
of the third relation to $acbac$  so that the relation  turns to 
$acbac=cacac$. We cancel $ac$ from right and obtain a new relation $acb=cac$.
Using this, one has $bcbac=bcaca=bacba=acaba=cacab=cbacb=cbcac$.  
We cancel $ac$ from right and obtain $bcb=cbc$.
Using this, one has $acabc=bacbc=babcb=abacb=abcac$.
Cancelling $a$ and $c$ for left and right, we obtain a new relation $cab=bca$.
Using this, one has $cabba=bcaba=bcbab=cbcab=cbbca$.
Cancelling $c$ and $a$ for left and right, we obtain a new relation $abb=bbc$.
The last relation of length 4 is obtained by cancelling $a$ from left of the equality: $abbac=abbac=bbcac=bbacb=bacab=acaab$.


\medskip
$\mathrm{H_{ii}}$:  Using $aca=bac$, rewrite the LHS $ acaca$  
of the third relation to $acbac$  so that the relation  turns to 
$acbac=cacac$. We cancel $ac$ from right so that we obtain  a new relation $acb=cac$.

\medskip
$\mathrm{H_{iii}}$: Multiply $b$ to the second relation from the right, 
and rewire the LHS to $bcaba$ (by a use of $bab=aba$ and rewrite the RHS 
to $cbcba$ (by a use of $acb=cba$). Cancelling by $ba$ from  right,  
we obtain a new relation $bca=cbc$. 

Using the length 3 relations, on has 
$acabc = acbcb = cbacb = cbcba = cabca = cacbc$.
Cancelling by $bc$ from right, we obtain a new relation $aca=cac$.

Using the length 3 relations, on has 
$bcaac= cbcac = cbaca = acbca = abcaa = bcbaa$.
Cancelling by $bc$ from left, we obtain a new relation $aac=baa$.

In the above sequence, the middle term $acbca$ is also equivalent to $accbc$.
Thus, cancelling $c$ from right, we obtain a new relation $accb\!=\!cbca (\!=\!bcaa)$.

$\mathrm{H_{viii}}$: From the defining relations, we have 
$abababa=bcbcbca$, $bababab=bbcbcbc$, and, hence, $bcbcbca=bbcbcbc$. 
Dividing by $b$ from the left, we get $cbcbca=bcbcbc$. The left hand side of this equality is equivalent to $cabbca=acbbca$, and the right hand side of the equality is equivalent to $abbcbc$ so that $acbbca\simeq abbcbc$. dividing by $a$ from the left, we get $bbcbc\simeq cbbca \simeq cbbac$. Dividing by $c$ from the light, we get $cbba\simeq bbcb (\simeq babb)$ (1). Multiplying $cbcb$ from the right, we get 
 $cbbacbcb\simeq bbcbcbcb$. The right hand side is equivalent to $bbcbcbcb\simeq bcbcbcbc\simeq cbcbcbcc \simeq cbcbcabc\simeq cbcbacbc$. The left hand side is equivalent to $ cbbacbcb\simeq cbbcabcb\simeq cbababcb$, and hence $cbababcb\simeq cbcbacbc.$ dividing by $cb$ from the left, we get $cbacbc\simeq ababcb$. The left hand side is equivalent to $cbacab\simeq cbaacb$. Dividing by $cb$ from the right, we get $abab\simeq cbaa$ (2). Mutiplying $b$ from the right, the left hand side  is equivalent to $acbba\simeq cabba\simeq cbcba$ so that $cbcba\simeq cbaab$. Dividing by $cb$ from the left, we get $cba\simeq aab$ (3). 
Applying (3) to the equality (2), we get $abab\simeq cbaa\simeq aaba$. Dividing by $a$ from the left, we get $bab\simeq aba\simeq bca$. Dividing by $b$ from the left, we get $ab\simeq  ca=ac$, and hence $b=c$. 
\end{proof}

\noindent
{\bf Notation.} {\it For each type $X\in\{\mathrm{A_i,A_{ii},B_i,B_{ii},B_{iii},B_{iv},B_v,B_{vi},}$ $\mathrm{B_{vii},H_i,H_{ii},H_{iii},}$ $\mathrm{H_{iv},H_v},$ $\mathrm{H_{vi},H_{vii},H_{viii}}\}$,  we denote by $G_X$, $M_X$ and $G_X^+$ the {\bf group},
the {\bf monoid} and the {\bf image of localization}: $M_X\!\to\! G_X$, respectively, associated with the 
positive homogeneous relations of type $X$ given in {\bf Theorem.\! 1}}.

\smallskip
From the presentations, 
we immediately observe the followings.

\medskip
\noindent
{\bf Corollary.} i) {\it For the type $X\in \{{\mathrm{A_i}},{\mathrm{A_{ii}}}, {\mathrm{B_i}},{\mathrm{B_{iv}}}, {\mathrm{H_i}}\}$, the monoid $M_X$ and the group $G_X$ is an  Artin monoid and an Artin group of type $\mathrm{A_3},\ \mathrm{A_2},\ \mathrm{B_3},\ \mathrm{A_3},\ \mathrm{A_1}\times\mathrm{A_2}$ and $\mathrm{H_3}$, respectively. We have 
the natural isomorphisms: $M_X\simeq G_X^+$.}

\medskip
ii) {\it For the type $X\in \{{\mathrm{B_v}}, {\mathrm{B_{vii}}},{\mathrm{H_{iv}}}, {\mathrm{H_{v}}}, {\mathrm{H_{vi}}}, {\mathrm{H_{vii}}}, {\mathrm{H_{viii}}}\}$, the monoid $M_X$ and the group $G_X$ is the infinite cyclic monoid $\Z_{\ge0}$ and group $\Z$, respectively. The monoid $M_{\mathrm{B_{iii}}}$ and the group $G_{\mathrm{B_{iii}}}$ is a free abelian monoid $(\Z_{\ge0})^2$ and group $\Z^2$ of rank 2. We have the natural isomorphisms: $M_X\simeq G_X^+$.
 }

\medskip
iii) {\it The correspondence: $\{ a\mapsto b, b\mapsto a, c\mapsto c\}$ induces an isomorphism:
\[ \vspace{-0.1cm}
M_{\mathrm{B_{vi}}} \quad \simeq \quad M_{\mathrm{H_{iii}}}
\]
and, hence, also the isomorphisms:
$G_{\mathrm{B_{vi}}}  \simeq  G_{\mathrm{H_{iii}}}$
and 
$G^+_{\mathrm{B_{vi}}} \simeq G^+_{\mathrm{H_{iii}}}$.
} 

\noindent
 ({\it Proof.} We can show that the Zariski-van Kampen relations of one of the two types can be deduced, up to the transposition of $a$ and $b$, from that of the other type.\ $\Box$)
 \ Note that the isomophism does not identify the Coxeter elements.

\bigskip
As the consequence of {\bf Corollary}, in the rest of the present paper, we shall focus our attention to the remaining 4 types $\mathrm{B_{ii}, B_{vi}, H_{ii}}$ and $\mathrm{H_{iii}}$ together with the ``constraint $\mathrm{B_{vi} \simeq H_{iii}}$''.


\begin{rem} The group $G_X$ is naturally isomorphic to the fundamental group, which does not depend on the choice of Zariski-van Kampen generators $\{a,b,c\}$, but the monoid $G_X^+$ depends on that choice (see next Remark 4.2).  

Further more, the monoid $M_X$, a priori, depends on the choice of relations in Theorem 1. The isomorphism $M_X\simeq G_X^+$ in the above corollary follows from cancellation conditions on $M_X$ (see \cite{[B-S]}). We shall show that, also for $M_{\mathrm{B_{ii}}}$ in \S7, the cancellation condition holds, implying $M_{\mathrm{B_{ii}}}\!\!\simeq\!\! G_{\mathrm{B_{ii}}}^+$. Thus, for these cases as a consequence of the cancellation condition, $M_X$ does not depend on the choice of relations in Theorem 1. However, for the remaining types  ${\mathrm{B_{vi}}}$,\! ${\mathrm{H_{ii}}}$ and ${\mathrm{H_{iii}}}$, it may be still possible that we need more relations in order to obtain the isomorphism $M_X\!\!\simeq\!\! G_X^+$.\!
\end{rem}
\begin{rem} Recall that, in the present paper, the generators $a,b,c$ are presented by the paths, which start from the base point $*_X$ and move along the intervals connecting $*_X$ and the three points $D_X\cap l_{*1}$ in the pencil $l_{*1,\C}$ and turn once counterclockwise the points $D_X\cap l_{*1}$ and then return to $*_X$ along the interval (see Fig. 2). Then, the set of the tuples generator system $a,b,c$ explained in \S3.3 admits the action of the braid group $B(3)$ of three strings, which changes associated relations. Here is a remarkable observation.

\bigskip
\noindent
{\bf Assertion.} 
{\it  Recall the projection $\pi:S_{\mathrm{B_{ii}}}\simeq \C^3\to T_{\mathrm{B_{ii}}}\simeq \C^2$. Then,  for any choice of   Zariski-van Kampen generator system $\{a,b,c\}$ (up to a permutation) in a pencil with respect to $\pi$ (i.e.\ a fiber of $\pi$) admit only one of the following two presentations I. and II.
}

\[
\begin{array}{lll}
\mathrm{I:}&\qquad
 \biggl{\langle}
a,b,c
\biggl{|}
\begin{array}{cc}cbb=bba,\\
 bc=ab,\\
 ac=ca
  \end{array}\
\biggl{\rangle}\ .\\


\mathrm{II:}&\qquad
 \biggl{\langle}
a,b,c
\biggl{|}
\begin{array}{cc}ababab=bababa,\\
 b=c,\\
 aabab=baaba
  \end{array}\
\biggl{\rangle}\ .
\vspace{0.5cm}\\
\end{array}
\]
\end{rem}

\medskip
\noindent
{\bf Corollary.} 
{\it The groups $G_{\mathrm{B_{vi}}}$ and $G_{\mathrm{H_{iii}}}$ do not admit Artin group presentation with respect to any Zariski-van Kampen type generator system.}
\begin{proof} Due to Theorem 1., both groups  have the relations: $a^5=b^5=c^5$, which are invariant by the change of generator system by the braid group $B(3)$.
\end{proof}

\section{Non-division property of the monoid $G^+_X$}

In the present section, we show that none of the monoids $G_X^+$ of the four types 
 $\mathrm{B_{ii}}, \mathrm{B_{vi}},\mathrm{H_{ii}}$ 
and $\mathrm{H_{iii}}$ does admit the 
divisibility theory (\cite[\S4]{[B-S]}), and therefore the monoid is neither Gaussian, Garside nor Artin.

\medskip
We first recall some terminology and concepts on the monoid $G^+$.

\noindent
An element $U\!\in\! G^+$ is said 
to {\it divide} $V\!\in\! G^+$ from the left (resp.~right), 
denoted by $U|_lV$ (resp.~$U|_rV$), if there exists $W\!\in\! G^+$ 
such that $V\!=\! UW$ (resp.~$V\!=\! WU$). 
We also say $V$ is {\it left-divisible} by $U$, or $V$ is 
a {\it left-multiple} of $U$. 

We say that $G^+$ {\it admits the left {\rm (resp.} right{\rm )} 
divisibility theory}, if 
for any two elements $U,V$ of $G_X^+$, there always exists their 
left (resp.~right) least common multiple, i.e.~a left (resp.~right) common multiple which divides any other left (resp.~right) common multiple,
denoted by $\mathrm{lcm_l}(U,V)$ (resp.~$\mathrm{lcm_r}(U,V)$). 

\medskip
\noindent
{\bf Theorem 2.}
{\it The monoids $G_{\mathrm{B_{ii}}}^+,G_{\mathrm{B_{vi}}}^+, G_{\mathrm{H_{ii}}}^+, G_{\mathrm{H_{iii}}}^+$  admits neither the left-divisibility theory nor the right divisibility theory.}
\begin{proof} 
We claim a fact, which shall be proven in \S8 Theorem 5 ii) independent of the results of \S5, 6 and 7.

\medskip
\noindent
{\bf Fact.}  None of the groups $G_{\mathrm{B_{ii}}}, G_{\mathrm{B_{vi}}}, G_{\mathrm{H_{ii}}}$ and $G_{\mathrm{H_{iii}}}$ is abelian.

\medskip
Assuming that the monoid $G^+_X$ admits the left division theory, we show that $G_X$ becomes an abelian group: a  contradiction! to {\bf Fact}.  The case for the right-division theory can be shown similarly.

\medskip

1) $G_{\mathrm{B_{ii}}}^+$: 
%
It is immediate to see $l(\mathrm{lcm}_l(b,c))>2$ from the defining relations in Theorem 1. Then, $bba=cbb$ is a common multiple of $b$ and $c$ of the shortest length 3, and, hence, should be equal to $\mathrm{lcm}_l(b,c)$. On the other hand, we have the following sequence of elementary equivalent words: $bcba,\ bbba,\ acbb,\ cabb$. That is, $bcba=cabb$ in $G_{\mathrm{B_{ii}}}^+$ is another common left-multiple of $b$ and $c$. 
If  $bba=cbb$ divides $bcba=cabb$ from the left, there exists $d\in\{a,b,c\}$ such that $bcba=bbad$. So, in $G_{\mathrm{B_{ii}}}^+$, we have $cba=bad$ which is again a common left-multiple of $b$ and $c$. Thus, we have the equality: $cba=cbb$ in $G_{\mathrm{B_{ii}}}^+$. That is, $a=b$ in $G_{\mathrm{B_{ii}}}^+$. 
By adding this relation  $a=b$ to the set of the defining relations of the group $G_{\mathrm{B_{ii}}}$, 
we get $G_{\mathrm{B_{ii}}}\!\simeq\! \Z$.
A contradiction!

\medskip
2) $G_{\mathrm{B_{vi}}}^+$:  
Due to the first defining relation in Theorem 1., we have $l(\mathrm{lcm}_l(a,b))$ 
\noindent
$\le 3$.
Let us consider 3 cases:

i) $l(\mathrm{lcm}_l(a,b))=1$. This means $l(\mathrm{lcm}_l(a,b))=a=b$.
By adding this relation to the defining relation of the group 
$G_{\mathrm{B_{vi}}}$, we get $G_{\mathrm{B_{vi}}}\simeq\Z$.
A contradiction!

ii) $l(\mathrm{lcm}_l(a,b))=2$. This means that there exists $u,v\in\{a,b,c\}$ 
such that $l(\mathrm{lcm}_l(a,b))=au=bv$. 
Depending on each choice of $u$ and $v$, one can show that this assumption 
leads to a contradictory conclusion $G_{\mathrm{B_{vi}}}\simeq\Z$.
Details are left to the reader.

iii) $l(\mathrm{lcm}_l(a,b))=3$. In view of the first two defining relations 
in Theorem 1., one has $aba\!=\!bab\!=\!aca\!=\!bac$. By adding this relation to the set of the defining relations of the group $G_{\mathrm{B_{vi}}}$, 
we get 
$G_{\mathrm{B_{vi}}}\!\simeq\! \Z$. A contradiction!.

\medskip
3) $G_{\mathrm{H_{ii}}}^+$: 
Due to the second defining relation in Theorem 1., we have $l(\mathrm{lcm}_l(a,b))$

\noindent
$\le 3$.
Let us consider 3 cases:

i) $l(\mathrm{lcm}_l(a,b))=1$. This means $l(\mathrm{lcm}_l(a,b))=a=b$.
By adding this relation to the defining relation of the group 
$G_{\mathrm{H_{ii}}}$, we get a contradiction $G_{\mathrm{H_{ii}}}\simeq\Z$.

ii) $l(\mathrm{lcm}_l(a,b))=2$. This means that there exists $u,v\in\{a,b,c\}$ 
such that $l(\mathrm{lcm}_l(a,b))=au=bv$.
Depending on each choice of $u$ and $v$, one can show that this assumption 
leads to a contradictory conclusion $G_{\mathrm{H_{ii}}}\simeq\Z$.
Details are left to the reader.

iii) $l(\mathrm{lcm}_l(a,b))=3$.
In view of the first two defining relations, 
one has $\mathrm{lcm}_l(a,b)\!=\!aca\!=\!bac$, and it divides 
$abab\!=\!baba$ (from left). 
This means that there exist $d\!\in\!\{a,b,c\}$ such that 
$cd\!=\!ba$ in $G_{\mathrm{H_{ii}}}$. 
For each case $d\!=\!a, b$ or $c$ separately, one can show that
 $G_{\mathrm{H_{ii}}}\!\simeq\! \Z$. A contradiction!.

\medskip
4) $G_{\mathrm{H_{iii}}}^+$: 
Due to the first defining relation in Theorem 1., we have $l(\mathrm{lcm}_r(a,b))\!\le\! 3$.
Let us consider 3 cases:

i) $l(\mathrm{lcm}_r(a,b))=1$. This means $l(\mathrm{lcm}_r(a,b))=a=b$.
By adding this relation to the defining relation of the group 
$G_{\mathrm{H_{iii}}}$, we get a contradiction $G_{\mathrm{H_{iii}}}\simeq\Z$.

ii )$l(\mathrm{lcm}_r(a,b))=2$. This means that there exists $u,v\in\{a,b,c\}$ 
such that $l(\mathrm{lcm}_r(a,b))=ua=vb$.
Depending on each choice of $u$ and $v$, one can show that this assumption 
leads to a contradictory conclusion $G_{\mathrm{H_{iii}}}\simeq\Z$. 
Details are left to the reader.

iii) $l(\mathrm{lcm}_r(a,b))=3$. 
In view of the first two defining relations, 
one has $\mathrm{lcm}_r(a,b)\!=\!aba\!=\!bab\!=\!cba\!=\!acb$. 
This leads to a conclusion $G_{\mathrm{H_{iii}}}\!\simeq\! \Z$, which is a contradiction!.
%
\quad  These complete the proof of Theorem 2.
\end{proof}
\noindent
{\bf Corollary 5.1.} {\it The monoids $G_{\mathrm{B_{ii}}}^+,G_{\mathrm{B_{vi}}}^+, G_{\mathrm{H_{ii}}}^+, G_{\mathrm{H_{iii}}}^+$ 
are not Gaussian, where a monoid is Gaussian if it is atomic, cancellative and admits divisibility theory (\cite[\S2]{[D-P]}). Hence, they are neither Artin groups nor Garside groups.}

\section{Fundamental elements of the monoid $M_X$}

Artin monoid of finite type has a particular element, denoted by $\Delta$ and called the {\it fundamental element} ([B-S] \S6). We want to generalize the concept for our new setting. However, in view of Theorem 2, we cannot employ the original definition: the left and right least common multiple of the generators. Analyzing equivalent defining properties of the fundamental element for Artin monoid case, we consider two classes of elements in the monoid $M$:
quasi-central elements and fundamental elements, forming submonoids 
$\mathcal{QZ}(M)$ and $\mathcal{F}(M)$ in $M$, respectively, with 
$\mathcal{F}(M)\subset \mathcal{QZ}(M)$.
The goal of the present section is to show $\mathcal{F}(M_X)\not=\emptyset$ for all types $X$, implying also $\mathcal{F}(G^+_X)\not=\emptyset$ for all types $X$.

\medskip
Let $M$ be a monoid given in \S4, i.e.\ defined by a positive homogeneous relations on a generator set $L$. Let us denote by $L/\!\sim$ the quotient set of $L$ divided by the equivalence relation generated by the equalities between two alphabets (in the relation set $R$).
An element $\Delta\in M$ is called {\it quasi-central} ([B-S] 7.1), if there exists a permutation $\sigma_\Delta$ of $L/\!\!\sim$ such that 
\[
a\cdot\Delta= \Delta\cdot \sigma_\Delta(a)
\vspace{-0.2cm}
\]
holds for all generators $a\in L/\!\sim$. 
 The set of all quasi-central elements is denoted by $\mathcal{QZ}(M)$.
The following is an immediate consequence of the definition.

\smallskip
\noindent
{\bf Fact 2.} {\it The $\mathcal{QZ}(M)$ is closed under the product.
For two elements $\Delta_1,\Delta_2\in \mathcal{QZ}(M)$,
we have $\sigma_{\Delta_1\cdot\Delta_2}=\sigma_{\Delta_2}\cdot\sigma_{\Delta_1}$.
}

\smallskip
According to {\bf Fact 2.}, we introduce an anti-homomorphism:
\[
\begin{array}{lll}
\sigma\ :& \mathcal{QZ}(M) \ \longrightarrow \ \mathfrak{S}(L/\!\sim), 
\quad \Delta \mapsto \sigma_\Delta.
\end{array}
\]
The kernel of $\sigma$ is the center $\mathcal{Z}(M)$ of the monoid $M$. 

\medskip
Next, we introduce the concept of a fundamental element.

\begin{defn}
An element $\Delta\in M$ is called {\it fundamental}
if there exists a permutation $\sigma_\Delta$ of $L/\!\sim$ such that, for any $a\in L/\!\sim$, there exists 
$\Delta_a\in G_X^+$ satisfying the following relation:
\[
\Delta\ =\ a\cdot\Delta_a\ =\ \Delta_a \cdot \sigma_\Delta(a).
\] 
We denote by $\mathcal{F}(M)$ the set of all fundamental elements of $M$. 
Note that $1\in \mathcal{QZ}(M)$ but $1\not\in \mathcal{F}(M)$
\end{defn}

\medskip
\noindent
{\bf Fact 3.} {\it The $\mathcal{F}(M)$ is an {\it idealistic 
submonoid} of $\mathcal{QZ}(M)$. That is, the following two properties hold.

{\rm i)} A fundamental element is a quasi-central element:  $\mathcal{F}(M)\subset \mathcal{QZ}(M)$. The associated permutation of $L/\!\sim$ as a fundamental element coincides with that as a quasi-central element.


{\rm ii} Products $\Delta\cdot \Delta'$ and  $\Delta'\cdot \Delta$ 
of a fundamental element $\Delta$  and a quasi-central element 
$\Delta'$ are again fundamental elements whose permutation of $L/\!\sim$ is given in {\bf Fact 2}.
We have  $(\Delta\Delta')_a=\Delta_a\Delta'$,
and  $(\Delta'\Delta)_a=\Delta'\Delta_{\sigma_{\Delta'}(a)}$.
}
\[
\mathcal{F}(M)\mathcal{QZ}(M)=\mathcal{QZ}(M) \mathcal{F}(M)= \mathcal{F}(M).
\]
\begin{proof}
i) We have 
$a\cdot \Delta= a\!\cdot\! \Delta_a\!\cdot\! \sigma_\Delta(a)\!=\! \Delta\!\cdot\! \sigma_\Delta(a)$ for all $a\!\in\! L/\!\!\sim$.

ii) We prove only the case $\Delta\cdot \Delta'$.

On one side, one has:
 
\centerline{
$\Delta\cdot \Delta'\simeq (a\cdot \Delta_a)\cdot\Delta'
\simeq a\cdot (\Delta_a\cdot\Delta')$.
}

On the other side, one has:

\quad$\Delta\cdot \Delta'
\simeq (\Delta_a\cdot \sigma_\Delta(a))\cdot \Delta'
\simeq \Delta_a \cdot (\sigma_\Delta(a)\cdot \Delta')
\simeq \Delta_a\cdot ( \Delta' \cdot \sigma_{\Delta'}(\sigma_\Delta(a)))$

\qquad \qquad $\simeq (\Delta_a\cdot  \Delta') \cdot \sigma_{\Delta'}(\sigma_\Delta(a))\simeq (\Delta_a\cdot  \Delta') \cdot \sigma_{\Delta\Delta'}(a))$.
\end{proof}

One basic property of a fundamental element is that it can be a universal denominator for the localization morphism (c.f.\ \S7 Lemma7.2\ 2.). 

\medskip
\noindent
{\bf Fact 4.}
{\it Let $\Delta$ be a fundamental element of $M$. 
Then, for any $U\in M$, $U$ divides $\Delta^{l(U)}$ from the left and from the right.}
\begin{proof} We prove only for the left division. Right division 
can be shown similarly.
We show the statement by induction on $l(U)$, where the case $l(U)=1$ 
follows from the definition of a fundamental element.
Let $l(U)>1$ and $U\simeq U'\cdot a$. By induction hypothesis, we 
have
$\Delta^{l(U)-1}\simeq U'\cdot V$ for some $V$. Then, multiplying 
$\Delta$ from the right, we have 
$\Delta^{l(U)}
\simeq U'\cdot V\cdot \Delta 
\simeq U'\cdot \Delta \cdot \sigma(V)
\simeq U'\cdot a\cdot  \Delta_a \cdot \sigma(V)$. 
Here, if $V$ is a word $v_1\cdots v_n$ then $\sigma(V)$ is a word  $\sigma(v_1)\cdots \sigma(v_n)$
\end{proof}

\begin{rem}
If $M$ is an indecomposable Artin monoid (of finite type), then any non-trivial quasi-central element is fundamental (\cite{[B-S]} 5.2 and 7.1). That is, one has the ``opposite'' inclusion: $(\mathcal{QZ}(M)\!\setminus\{1\}) \subset\mathcal{F}(M)$.
\vspace{-0.4cm}
\end{rem}¡¡
\begin{rem}
By the definition, any fundamental element is divisible from both left and right by all generators in $L$. However, (non-trivial) quasi-central element in general may not have this property.

(i) $b^3 \in\mathcal{QZ}(M_{\mathrm{B_{ii}}})$ is central. However, it is not divisible by $a$ and $c$ from the left and right.


(ii) $ababa\in M_{\mathrm{B_{ii}}}$ is divisible by all generators from both sides, but it does not belong to $\mathcal{QZ}(M_{\mathrm{B_{ii}}})$.
\end{rem}

We state the second main result of the present paper.

\medskip
\noindent
{\bf Theorem 3.} {\it 
The following elements $\Delta_X$ belong to $\mathcal{F}(M_X)$ for any type $X$. 
\[
\begin{array}{ll ll}
\mathrm{A_i}: \qquad & \Delta_{\mathrm{A_i}}\!&\!:=\ (cba)^2  \qquad & \sigma: \big(^{a,\ b,\ c}_{c,\ b,\ a}\big) \\

\mathrm{A_{ii}}:\qquad & \Delta_{\mathrm{A_{ii}}}\!&\!:=\ aba \qquad & \sigma: \big(^{a,\ b=c}_{b=c,\ a}\big)\\

\mathrm{B_i}: \qquad & \Delta_{\mathrm{B_i}}\!&\!:=\ (cba)^3 \qquad & \sigma: \big(^{a,\ b,\ c}_{a,\ b,\ c}\big) \\

\mathrm{B_{ii}}: \qquad & \Delta_{\mathrm{B_{ii1}}}\!&\!:=\ (ab)^3 \qquad & \sigma: \big(^{a,\ b,\ c}_{a,\ b,\ c}\big)\\
& \Delta_{\mathrm{B_{ii2}}}\!&\!:=\ (bcc)^3 \simeq (cba)^3 & \sigma: \big(^{a,\ b,\ c}_{a,\ b,\ c}\big)\\

\mathrm{B_{iii}}: \qquad &  \Delta_{\mathrm{B_{iii}}}\!&\!:=\ ac \qquad & \sigma: \big(^{a= b,\ c}_{a= b,\ c}\big)\\

\mathrm{B_{iv}}:  \qquad & \Delta_{\mathrm{B_{iv}}}\!&\!:=\ abcb \qquad & \sigma: \big(^{a,\ b,\ c}_{a,\ c,\ b}\big)\\

\mathrm{B_{v}}:  \qquad & \Delta_{\mathrm{B_{v}}}\!&\!:=\ a \qquad & \sigma: \big(^{a= b= c}_{a= b= c}\big)\\

\mathrm{B_{vi}}:  \qquad &\Delta_{\mathrm{B_{vi}1}}\!&\!:=\ a^5\simeq b^5\simeq c^5&\sigma: \big(^{a,\ b,\ c}_{a,\ b,\ c}\big)\\
&\Delta_{\mathrm{B_{vi}2}}\!&\!:=\ (aba)^2 \qquad & \sigma: \big(^{a,\ b,\ c}_{a,\ b,\ c}\big) \\
&\Delta_{\mathrm{B_{vi}3}}\!&\!:=\ bccabcb \qquad & \sigma: \big(^{a,\ b,\ c}_{a,\ b,\ c}\big) \\
&\Delta_{\mathrm{B_{vi}4}}\!&\!:=\ (bbac)^2&\sigma: \big(^{a,\ b,\ c}_{a,\ b,\ c}\big)\\
& \Delta_{\mathrm{B_{vi}5}}\!&\!:=\ (acaca)^2 \qquad & \sigma: \big(^{a,\ b,\ c}_{a,\ b,\ c}\big) \\
& \Delta_{\mathrm{B_{vi}6}}\!&\!:=\ (cba)^3 \qquad & \sigma: \big(^{a,\ b,\ c}_{a,\ b,\ c}\big) \\
& \Delta_{\mathrm{B_{vi}7}}\!&\!:=\ (cab)^5 \qquad & \sigma: \big(^{a,\ b,\ c}_{a,\ b,\ c}\big) \\
\end{array}
\]
\[
\begin{array}{ll ll}
\mathrm{B_{vii}}:  \qquad & \Delta_{\mathrm{B_{vii}}}\!\!&\!:=\ a \qquad & \sigma: \big(^{a= b= c}_{a= b= c}\big)\\

\mathrm{H_i}:  \qquad & \Delta_{\mathrm{H_i}}\!&\!:=\ (cba)^5 \qquad & \sigma: \big(^{a,\ b,\ c}_{a,\ b,\ c}\big)\\

\mathrm{H_{ii}}:  \qquad & \Delta_{\mathrm{H_{ii1}}}\!&\!:=\ (acaca)^2 \simeq (ac)^5  & \sigma: \big(^{a,\ b,\ c}_{a,\ b,\ c}\big)\\
& \Delta_{\mathrm{H_{ii2}}}\!&\!:=\ (babac)^3 \simeq (cba)^5& \sigma: \big(^{a,\ b,\ c}_{a,\ b,\ c}\big)\\
\mathrm{H_{iii}}: \qquad &\Delta_{\mathrm{H_{iii}1}}\!&\!:=a^5\simeq b^5\simeq c^5 \qquad & \sigma: \big(^{a,\ b,\ c}_{a,\ b,\ c}\big) \\
&\Delta_{\mathrm{H_{iii}2}}\!&\!:=(aba)^2 \qquad & \sigma: \big(^{a,\ b,\ c}_{a,\ b,\ c}\big) \\
&\Delta_{\mathrm{H_{iii}3}}\!&\!:=accbaca \qquad & \sigma: \big(^{a,\ b,\ c}_{a,\ b,\ c}\big) \\
&  \Delta_{\mathrm{H_{iii}}4}\!&:=\ (bcba)^2 \qquad & \sigma: \big(^{a,\ b,\ c}_{a,\ b,\ c}\big)\\
&\Delta_{\mathrm{H_{iii}5}}\!&\!:=\ (bcbcb)^2\simeq (bc)^5 \quad & \sigma: \big(^{a,\ b,\ c}_{a,\ b,\ c}\big) \\
&\Delta_{\mathrm{H_{iii}6}}\!&\!:=(abc)^3 \qquad & \sigma: \big(^{a,\ b,\ c}_{a,\ b,\ c}\big) \\
&\Delta_{\mathrm{H_{iii}7}}\!&\!:=(cba)^5 \qquad & \sigma: \big(^{a,\ b,\ c}_{a,\ b,\ c}\big) \\

\mathrm{H_{iv}}:  \qquad & \Delta_{\mathrm{H_{iv}}}\!&:=\ a \qquad & \sigma: \big(^{a= b= c}_{a= b= c}\big)\\

\mathrm{H_{v}}:  \qquad & \Delta_{\mathrm{H_{v}}}\!&\!:=\ a \qquad & \sigma: \big(^{a= b= c}_{a= b= c}\big)\\

\mathrm{H_{vi}}:  \qquad & \Delta_{\mathrm{H_{vi}}}\!&\!:=\ a \qquad & \sigma: \big(^{a= b= c}_{a= b= c}\big)\\

\mathrm{H_{vii}}:  \qquad & \Delta_{\mathrm{H_{vii}}}\!\!&\!:=\ a \qquad & \sigma: \big(^{a= b= c}_{a= b= c}\big)\\

\mathrm{H_{viii}}: \qquad &  \Delta_{\mathrm{H_{viii}}}\!\!\!&\!:=\ a \qquad & \sigma: \big(^{a= b= c}_{a= b= c}\big)\\
\end{array}
\]
}
\begin{proof} 
Since the cases for an Artin monoid or a free abelian monoid 
are classical, we show only the 4 exceptional cases.

\medskip
\noindent
$\mathrm{B_{ii}}$ :

$\Delta_{\mathrm{B_{ii}1}}:=ababab.$

$\Delta_{\mathrm{B_{ii}1}}=a(babab) $,

$\Delta_{\mathrm{B_{ii}1}}\simeq bcabab\simeq bacbab\simeq bacbbc \simeq babbac 
\simeq babbca\simeq (babab)a$.

\smallskip
$\Delta_{\mathrm{B_{ii}1}}\simeq b(ababa)$,

$\Delta_{\mathrm{B_{ii}1}} = (ababa)b $.

\smallskip
$\Delta_{\mathrm{B_{ii}1}}\simeq bababa\simeq bbcaba\simeq bbacba\simeq c(bbcba)$,

$\Delta_{\mathrm{B_{ii}1}}\simeq bcbcbc \simeq bcabbc \simeq bacbbc \simeq babbac \simeq (bbcba)c$.

\smallskip
$\Delta_{\mathrm{B_{ii}2}}:=(bcc)^3.$

\smallskip
$\Delta_{\mathrm{B_{ii}2}}=b(ccbccbcc) \simeq abcbccbcc \simeq aabbccbcc \simeq aabbccabc $

\qquad 
$\simeq aabbaccbc \simeq aacbbccbc\simeq caabbccbc\simeq caabbccab $

\qquad
$\simeq caabbaccb \simeq caacbbccb \simeq ccaabbccb \simeq ccabcbccb$

\qquad 
$ \simeq ccbccbccb=( ccbccbcc)b.$

\smallskip
$\Delta_{\mathrm{B_{ii}2}} \simeq a(bcbccbcc) \simeq bccbccbcc \simeq bccabcbcc \simeq bccaabbcc $

\qquad
$\simeq bcaacbbcc \simeq bcaabbacc \simeq bcabcbacc \simeq bcbccbacc $

\qquad 
$\simeq bcbccbcca =(bcbccbcc)a.$

\smallskip
$\Delta_{\mathrm{B_{ii}2}} \simeq c(aacbbccb) \simeq aaccbbccb \simeq aacbbaccb \simeq aacbbcacb$

\qquad
 $\simeq aacbbccab \simeq aacbbccbc=(aacbbccb)c.$

\medskip
\noindent
$\mathrm{B_{vi}}$ : Due to {\bf Remark.} after {\bf Theorem 1.} in \S4,  we may reduce the proof to the case $\mathrm{H_{iii}}$.

\medskip
\noindent
$\mathrm{H_{ii}}$ :  First, let us show a relation: $acaca=cacac$ ($acaca\simeq acbac\simeq cacac$), which shall be used in the sequel. 

\smallskip
$\Delta_{\mathrm{H_{ii}1}}:=acacaacaca.$

$\Delta_{\mathrm{H_{ii}1}}=a(cacaacaca)\simeq cacacacaca \simeq (cacaacaca)a$.

$\Delta_{\mathrm{H_{ii}1}}\simeq c(acacacaca)\simeq acacaacaca\simeq (acacacaca)c $.

\smallskip
$\Delta_{\mathrm{H_{ii}1}}\simeq acacaacaca\simeq b(accaacaca)\simeq acacaacaca\simeq acacacacac$

\qquad 
$\simeq accacaccac\simeq accaacbcac\simeq accaacbacb\simeq (accaacaca)b$,

$\Delta_{\mathrm{H_{ii}2}} := babacbabacbabac \simeq ababcbabacbabac$.

\smallskip
$\Delta_{\mathrm{H_{ii}2}}=a(babcbabacbabac)\simeq bababcbabacbabac\simeq babcacabacbabac$

\qquad 
$\simeq babcbacbacbabac\simeq babcbacacababac \simeq babcbaacbababac$

\qquad 
$\simeq babcbaacababbac\simeq babcbabacbabbac\simeq (babcbabacbabac)a$.

$\Delta_{\mathrm{H_{ii}2}}=b(abacbabacbabac)\simeq ababcbabacbabac\simeq ababcababcbabac$ 

\qquad 
$\simeq ababcababcbaaca\simeq ababcbabacbaaca\simeq ababacbaacabaaca$ 

\qquad 
$\simeq ababcbaacababac\simeq ababcbaacbabaac\simeq ababcbacacabaac$

\qquad 
$\simeq ababcacaacabaac\simeq ababacbaacabaac\simeq abaacabaacabaac$

\qquad 
$\simeq abaacababacbaac\simeq abaacbabaacbaac\simeq abacacabaacbaac$

\qquad 
$\simeq abacbacbaacbaac\simeq abacbacbacacaac\simeq abacbacacaacaac$

\qquad 
$\simeq abacbaacbaacac\simeq abacbaacbabacac\simeq abacbaacababcac$

\qquad 
$\simeq abacbabacbabcac\simeq( abacbabacbabac)b  $.

$\Delta_{\mathrm{H_{ii}2}}=abacbabacbabacb \simeq aacababacbabacb\simeq  aacbabaacbabacb$ 

\qquad 
$\simeq acacabaacbabacb\simeq acbacbaacbabacb\simeq c(acacbaacbabacb)$

\qquad
$\simeq acbacbaacbabacb\simeq acacabaacababcb\simeq acacababacbabcb$

\qquad 
$\simeq acacbabaacbabcb\simeq acacbabacacabcb \simeq acacbaacaacabcb$

\qquad 
$\simeq acacb aacabacbcb\simeq acacbaacababcbc\simeq (acacbaacbabacb)c$.

\medskip
\noindent
$\mathrm{H_{iii}}$ :

$\Delta_{\mathrm{H_{iii}}2}:= (aba)^2.$

$\Delta_{\mathrm{H_{iii}}2}= a(baaba)\simeq bababa
\simeq (baaba)a.$

$\Delta_{\mathrm{H_{iii}}2}=b(ababa)\simeq abaaba\simeq (ababa)b.$

$\Delta_{\mathrm{H_{iii}}2}=abaaba\simeq aaacba\simeq aacbaa\simeq aacaac$

\smallskip
$\Delta_{\mathrm{H_{iii}}3}:=accbaca.$

$\Delta_{\mathrm{H_{iii}}3}=a(ccbaca)\simeq
cbcaaca\simeq
ccbcaca\simeq
(ccbaca)a.$
     
 $\Delta_{\mathrm{H_{iii}}3}= accbaca\simeq
   cbcaaca\simeq
    b(caaaca).$

$\Delta_{\mathrm{H_{iii}}3}= accbaca\simeq
    cbcaaca\simeq
    ccbcaca\simeq
   caccbca\simeq
   cacbcaa$

\qquad 
$\simeq
    caaccba\simeq
   caacacb\simeq
  ( caaaca)b.$

 $\Delta_{\mathrm{H_{iii}}3}=   accbaca\simeq
     c(bcaaca).$
    
 $\Delta_{\mathrm{H_{iii}}3}=      accbaca\simeq
     cbcaaca\simeq
     bcaaaca\simeq
    =(bcaaca)c.$

\smallskip
$\Delta_{\mathrm{H_{iii}}4}:=bcbabcba.$

$\Delta_{\mathrm{H_{iii}}4}\simeq a(bcabcba)\simeq bcabacba \simeq (bcabcba)a$.

\smallskip
$\Delta_{\mathrm{H_{iii}}4}=b(cbabcba)\simeq bcabacba\simeq cbcbacba\simeq cbacbcba$

\qquad $\simeq cbabcaba \simeq (cbabcba)b $.

\smallskip
$\Delta_{\mathrm{H_{iii}}4}\simeq bcabacba\simeq c(bcbacba)\simeq bcbaabca \simeq bcbaacbc \simeq bcbacbac$.

\medskip
$\Delta_{\mathrm{H_{iii}}5}:= bcbcbbcbcb$.

$\Delta_{\mathrm{H_{iii}}5}\simeq abccbbcbcb\simeq abccbcbcbc\simeq a(bcbcbcbbc)$,

$\Delta_{\mathrm{H_{iii}}5}\simeq bcbcbcbcbc\simeq (bcbcbcbbc)a.$

\smallskip
$\Delta_{\mathrm{H_{iii}}5} = b(cbcbbcbcb)\simeq cbcbcbcbcb\simeq (cbcbbcbcb)b$.

\smallskip
$\Delta_{\mathrm{H_{iii}}5} \simeq c(bcbcbcbcb) \simeq (bcbcbcbcb)c.$

\smallskip

$\Delta_{\mathrm{H_{iii}6}}:=(abc)^3$

$\Delta_{\mathrm{H_{iii}6}}= a(bcabcabc)\simeq bcbabcab\simeq bcabacabc\simeq bcabcacbc\simeq (bcabcabc)a$

$\Delta_{\mathrm{H_{iii}6}}\simeq (abcabcab)c\simeq abcabcbcb\simeq abcabbcab \simeq acbcbbcab$

\quad \qquad $\simeq acabcbcab\simeq cacbcbcab \simeq c(abcabcab)$.

\smallskip
$\Delta_{\mathrm{H_{iii}7}}:=(cba)^5$

$\Delta_{\mathrm{H_{iii}7}}=(cba)^5 \simeq (acb)^5 \simeq (bac)^5$. 
\end{proof}



As a consequence of Theorem 3, we have the folloing fact.

\smallskip
\noindent
{\bf Fact 5.} 
{\it There exists a positive integer $k\in\Z_{>0}$ such that the $k$-th power of the Coxeter element $C:=cba$ (= a homotopy class which turns once around all the three points $C_X\cap l_{*_1,\C}$ counterclockwise) is a fundamental element.}

\medskip
Finally, we ask a few questions related to the fundamental elements.

Let $M$ be a monoid defined by positive homogeneous relations. Recall (\S4 Definition) that $G^+$ is the image of $M$ in the group $G$ by the localization homomorphism.  We define quasi-central elements and fundamental elements of $G^+$ exactly by the same defining relations for $M^+$. Let us denote by $\mathcal{QZ}(G^+)$ and $\mathcal{F}(G^+)$ the set of quasi-central elements and fundamental elements in $G^+$, respectively. Then, the localization morphism induces homomorphisms: $\mathcal{QZ}(M)\to \mathcal{QZ}(G^+)$ and $\mathcal{F}(M)\to \mathcal{F}(G^+)$, which may be neither injective nor surjective. However, Theorem 3 implies the following fact.

\medskip
\noindent
{\bf Fact 6.} {\it For any type $X$, the set of fundamental elements $\mathcal{F}(G_X^+)$ is non-empty.}

\smallskip
We observe that  $\mathcal{F}(G_X^+)$ may not be singly generated. On the other hand, the list in Theorem 3 may not be sufficient to generate whole  $\mathcal{F}(M_X)$ or $\mathcal{F}(G_X^+)$.

\medskip
\noindent
{\bf Question 1.} Is $\mathcal{F}(M)$ (resp.\ $\mathcal{F}(G^+)$) finitely generated over $\mathcal{QZ}(M)$ (resp.\ $\mathcal{QZ}(G^+)$)? 
That is, are there finitely many elements $\Delta_1,\cdots, \Delta_k\in \mathcal{F}(M)$ (resp. $\mathcal{F}(G^+)$) such that following holds?
\[\begin{array}{ll}
\mathcal{F}(M) \  = \ \mathcal{QZ}(M)\Delta_1 \ \cup\ \cdots\ \cup\  \mathcal{QZ}(M)\Delta_k.\\
\mathcal{F}(G^+) = \mathcal{QZ}(G^+)\Delta_1 \ \cup\ \cdots\ \cup\  \mathcal{QZ}(G^+)\Delta_k.
\end{array}
\]

\noindent
{\bf Question 2.}  The following five cases 1, 2, 3, 4, and 5. give or may give example of {\it an indecomposable logarithmic free divisor such that the local fundamental group of its compliment admits positive homogeneous presentation by a suitable choice of Zariski-van Kampen generators and a power of the Coxeter element gives a fundamental element of the monoid generated by them}. 
We ask whether this property holds  for any indecomposable logarithmic free divisor or not (for a more precise formulation of the question, see \cite{[S-I2]}).

1. The discriminant of a finite irreducible reflection group (\cite{[B-S],[S2],[S3]}).

2. The discriminant of a finite irreducible complex refrection group (except for type $G_{31}$) (\cite{[B-M-R],[Be]}).

3. The Sekiguchi polynomials (Theorems 1. and 3. of the present paper).

4.  A plane curve is locally logarithmic free (see \cite{[S1]}). The local fundamental group of the complement of a plane curve seems  to be presented by positive homogeneous relations in \cite{[K]}. It seems likely that a power of the Coxeter element is a fundamental element of the associated monoid (to be confirmed yet).

5. The discriminant of elliptic Weyl group is a free divisor (\cite{[S4]}II). A Zariski-van Kampen presentation of the fundamental group of the complement of the divisor is not yet given. However, the hyperbolic Coxeter element in the elliptic Weyl group (\cite{[S4]}I,III) may (conjecturally) be lifted to the fundamental group, whose power of order $m_\Gamma$, gives a fundamental element.

\section{Cancellation conditions on $M_X$}

In the present section, we study the {\it cancellation condition} on a monoid $M$.
In the first half, we show some general consequences on the monoid $M$ under the cancellation condition, or under its weaker version: a {\it weak cancellation condition}. In the latter half, we prove that the monoid $M_{\mathrm{B_{ii}}}$ satisfies the cancellation condition, however, we do not know whether the monoids 
 $M_{\mathrm{B_{vi}}}$, $M_{\mathrm{H_{ii}}}$ and  $M_{\mathrm{H_{iii}}}$ satisfy it or not.

\smallskip

\begin{defn} A monoid $M$ is said to {\it satisfy the cancellation condition}, if 
an equality $AXB\!=\!AYB$ for $A,B,X,Y\!\in\!M$ implies $X\!=\!Y$.
\end{defn}

It is well-known that an Artin monoid satisfies the cancellation condition \cite[Prop.2.3]{[B-S]}.
Let us state some  important consequences of the cancellation condition on a monoid defined by positive homogeneous relations.

\medskip
\begin{lem}
{\it Let $M$ be a monoid defined by positive homogeneous relations. Suppose it satisfies the cancellation condition. Then, we have the followings.

\smallskip
\noindent
{\bf 1. } For any $\Delta\in  \mathcal{QZ}(M)$, the associated permutation $\sigma_\Delta$ of $L/\!\sim$ extends to an isomorphism, denoted by the same $\sigma_\Delta$, of $M$. The correspondence: $\Delta\mapsto \sigma_\Delta$ induces an anti-homomorphism:
\[
\mathcal{QZ}(M) \longrightarrow Aut(M).
\]
{\bf 2.}  If $\mathcal{F}(M)\not=\emptyset$, then the localization homomorphism is injective and, hence, one has an isomorphism:
\[
M \ \simeq \ G^+.
\]
{\bf 3.} For any element $A\!\in\! G$ and any $\Delta\!\in\!\mathcal{F}(M)$, there exists $B\!\in\! G^+$ and $n\!\in\!\Z_{\ge0}$ such that, in $G$, one has equalities:
}
\smallskip

\centerline{
$A=B\cdot (\Delta)^{-n}=(\Delta^{-n})\cdot\sigma_{\Delta}^n(B)$.
}
\end{lem}
\begin{proof}
{1.} First, we note that the permutation $\sigma_\Delta$ induces an isomorphism of the free monoid $(L/\!\sim)^*$, denote by the same $\sigma_\Delta$. Let $U$ and $V$ be words in $(L/\!\sim)^*$ which are equivalent by the relations $R$ (i.e.\ 
give the same element in $M$).
Then, by definition, $U\Delta \simeq\Delta \sigma_\Delta(U)$ and $V\Delta\simeq\Delta\sigma_\Delta(V)$ are equivalent. That is, $\Delta \sigma_\Delta(U)$  and $\Delta \sigma_\Delta(V)$ give the same element in $M$.  Then, cancelling $\Delta$ from the left, we see that $\sigma_\Delta(U)$ and $\sigma_\Delta(V)$ give the same element in $M$. Thus $\sigma_\Delta$ induces a homomorphism from $M$ to $M$. The homomorphism is 
invertible, since a finite power of it is an identity. By the definition, 
for any $U\in M$ and $\Delta_1,\Delta_2\in\mathcal{QZ}(M)$, one has:
\[
U \cdot \Delta_1\Delta_2 \ \simeq \ \Delta_1 \cdot \sigma_{\Delta_1}(U)\cdot\Delta_2\simeq \Delta_1\Delta_2 \cdot \sigma_{\Delta_2} (\sigma_{\Delta_1}(U)).
\]

{\bf 2.} For a localization morphism to be injective, it is sufficient to show 
that the monoid satisfies the cancellation condition and that any two elements of the monoid have (at least) one 
(left and right) common multiple (\"Ore's condition, see [C-P]). In view of Fact 4. in \S6, for any two elements $U,V\in M$ and $\Delta\in\mathcal{F}(M)$, $\Delta^{\max\{l(U),l(V)\}}$ is a common multiple of $U$ and $V$ from both sides. 

{\bf 3.} Owing to the previous {\bf 2.}, it is sufficient to show that, for any element $A\in G$ and  any $\Delta\in  \!\mathcal{F}(M)$, there exists $k\in \Z_{\ge0}$ such that $\Delta^k\cdot A\in G^+$. This can be easily shown by an induction on $k(A)\in \Z_{\ge0}$ where $k(A)$ is the (minimal) number of letters of negative power in a word expression of $A$ in $(L\cup L^{-1})^*$. Details are left to the reader.
\end{proof}

Next, we formulate a weak cancellation condition and its consequences.

\begin{defn}
An element $\Delta\in M$ is called 
left (resp.\ right) {\it weakly cancellative}, 
if an equality $\Delta=U\cdot V=U\cdot W$ (resp.\ $\Delta=V\cdot U=W\cdot U$) holds in $M$ for some $U,V,W\in M$, then $V=W$ holds in $M$.
%
\end{defn} 

\noindent
{\it Notation.} For an element $\Delta\in M$, we put

 $Div_l(\Delta):=\{U\in M :\ U\mid_l \Delta\}$ \ and \ $Div_r(\Delta):=\{U\in M :\ U\mid_r \Delta\}$.

\medskip
\noindent
{\bf Fact 7.}  {\it Let  a fundamental element $\Delta\in\mathcal{F}(M)$ be left weakly cancellative. 

\noindent
Then the following {\rm i), ii), iii)} and {\rm iv)} hold.

{\rm i)} For any element $U\in Div_l(\Delta)$, let $\tilde U\in (L/\!\sim)^*$ be a lifting to a word. 

\quad Then, the class of $\sigma_\Delta(\tilde U)$ in $M$ depends only on the class $U$ but not 

\quad on the lifting $\tilde U$. Let us denote the class in $M$ by $\sigma_\Delta(U)$. 

{\rm ii)} The divisor set $Div_l(\Delta)$ is invariant under the action of $\sigma_\Delta$.  In particular, 

\quad the unique longest element $\Delta$ is fixed by $\sigma_\Delta$.

{\rm iii)} The fundamental element $\Delta$ is right weakly cancellative.

{\rm iv)}  We have the equality:
\quad 
$
Div_l(\Delta) \ =\ Div_r(\Delta).
$
}
\begin{proof} i) 
Suppose one has a decomposition $\Delta\simeq U\cdot V$ for $U,V\in M$, and 
let $\tilde U$ be a lifting of $U$ into a word in $L/\!\sim$. Then, 
$\sigma_\Delta(\tilde U)$ is well-defined as a word and hence induce an element in $M$, which we denote by the same $\sigma_\Delta(\tilde U)$.
We claim that  $\Delta$ is equivalent to $V\cdot \sigma_\Delta(\tilde U)$. This is shown by induction on $l(U)$. If $l(U)=1$, this is the definition of 
quasi-centrality. Let $l(U)>1$, $\tilde U= \tilde U'\cdot a$ and 
$\Delta\simeq \tilde U'\cdot a \cdot V$. By induction hypothesis, we 
have $\Delta\simeq a\cdot V\cdot \sigma_\Delta(\tilde U')$.  
Due to the weak cancellativity, $V\cdot \sigma_\Delta(\tilde U')$ is equivalent 
to $\Delta_a$. Then, by definition of quasi-centrality, $\Delta$ 
is equivalent to $V\cdot \sigma_\Delta(\tilde U')\cdot \sigma_\Delta(a)
\simeq $ $V\cdot \sigma_\Delta(\tilde U)$. 

Let $\tilde U_1$ and $\tilde U_2$ be liftings of $U$.  Then, applying the above result, we see that $\Delta$ is  equal to $V\cdot \sigma_\Delta(\tilde U_1)$
and $V\cdot \sigma_\Delta(\tilde U_2)$. Then, applying the weak cancellativity of $\Delta$, we see that $\sigma_\Delta(\tilde U_1)$ and $\sigma_\Delta(\tilde U_2)$ define the same  element in $M$, which we shall denote by $\sigma_\Delta(U)$.

ii) In the proof of i), taking $U=\Delta$ and $V=1$, we obtain $\Delta=\sigma_\Delta(\Delta)$. Then, $\sigma_\Delta(Div_l(\Delta))=Div_l(\sigma_\Delta(\Delta))=Div_l(\Delta)$.

iii) Suppose $\Delta=V\cdot U=W\cdot U$. then according to i), we have $\Delta=U\cdot \sigma_\Delta(V)=U\cdot \sigma_\Delta(W)$. Then the left cancellation condition implies $\sigma_\Delta(V)=\sigma_\Delta(W)$. On the other hand, according to ii), $\sigma_\Delta(V)=\sigma_\Delta(W)$ are again elements of $Div_l(\Delta)$ so that we can apply $\sigma_\Delta$ to the equality. Since $\sigma_\Delta$ is of finite order, after repeating this several times, we obtain the equality $V=W$.

iv) If $\Delta$ is left divisible by $U$, 
$\Delta$ is right divisible by $\sigma_\Delta(U)$. 
That is, the set $Div_r(\Delta)$ of the right divisors of $\Delta$ is 
equal to $\sigma_\Delta(Div_l(\Delta))=Div_l(\Delta)$.
\end{proof}


\medskip
\noindent
{\bf Conjecture.}  Let $C^k$ of the element in \S6 {\bf Fact 5.} If $C^{k\cdot \mathrm{ord}(\sigma_{C^k})}$  is weakly cancellative, then $M$ satisfies the cancellation condition.

\medskip
{\small
The following Theorem shows that we have already enough relations for type $\mathrm{{B_{ii}}}$. 
}

\medskip
\noindent
{\bf Theorem 4.} {\it The monoid  $M_{\mathrm{B_{ii}}}$ satisfies the cancellation condition.}
\begin{proof}
We, first, remark  the following.

\medskip
\noindent
 {\bf Fact 8.}\  {\it The left cancellation condition on $\mathrm{M_{B_{ii}}}$ implies the right cancell.\ condition.}
\begin{proof} 
  Consider a map $\varphi:\mathrm{M_{B_{ii}}}\rightarrow\mathrm{M_{B_{ii}}}$,  $W\mapsto \varphi(W):=\sigma$$(rev(W))$, where $\sigma$ is a permutation $\big(^{a\,\, b\,\, c}_{c\,\, b\,\, a}\big)$ and $rev(W)$ is the reverse of the word $W=x_1 x_2 \cdots x_t$ ($x_i$ is a letter or an inverse of a letter) given by the word  $x_t  \cdots x_2 x_1$. In view of the defining relation of $\mathrm{M_{B_{ii}}}$ in Theorem 1., $\varphi$ is well defined and is an anti-isomorphism. If $\beta \alpha \simeq \gamma \alpha$, then $\varphi(\beta \alpha) \simeq \varphi(\gamma \alpha)$, i.e.,  $\varphi(\alpha) \varphi(\beta) \simeq \varphi(\alpha)\varphi(\gamma)$. Using left cancellation condition, we obtain $\varphi(\beta)=\varphi(\gamma)$ and, hence, $\beta \simeq \gamma $. 
\end{proof}

The following is sufficient to  show the left cancellation condition on $\mathrm{M_{B_{ii}}}$.

\medskip
\noindent
{\bf Proposition.}  {\it Let $X$ and $Y$ be positive words in ${\mathrm{M_{B_{ii}}}}$ of word-length $r\in \Z_{\ge0}$.
\smallskip
\\
{\rm (i)}\ If $uX \simeq uY$ for some $u \in \{a,b,c\}$, then $X \simeq Y$.\\
{\rm  (ii)}\ If $aX \simeq bY$, then $X \simeq bZ$, $Y \simeq cZ$ for some positive word $Z$.\\
{\rm (iii)}\ If $aX \simeq cY$, then $X \simeq cZ$, $Y \simeq aZ$ for some positive word $Z$.\\
{\rm (iv)}\ If $bX \simeq cY$, then there exists an integer $k$ ($0\!\le\! k\!<\!r\!-\!1$) and a word $Z$ 

\quad such that $X\simeq c^kbaZ$ and $Y\simeq a^kbbZ$.

\medskip




}
\vspace{-0.3cm}
\begin{proof} Let us denote by $H(r,t)$ the statement in Theorem for all pairs of words $X$ and $Y$ such that their word-lengths are $r$ and for all $u,v\in\{a,b,c\}$ such that $uX\simeq vY$ and the number of elementary transformasions to bring $uX$ to $vY$ is less or equal than $t$. It is easy to see that $H(r,t)$ is true if $r\le1$  or $t\le1$. 

For $r,t\in \Z_{>1}$, we prove $H(r,t)$ under the induction hypotheis that $H(s,u)$ holds for $(s,u)$ such that either $s< r$ and arbitrary $u$ or $s=r$ and $u< t$.

Let $X, Y$ be of word-length $r$, and let $u_1 X\simeq u_2W_2 \simeq\cdots \simeq u_{t}W_{t} \simeq u_{t+1}Y$ be a sequence of elementary transformations of $t$ steps, where $u_1,\cdots,u_{t+1}\in\{a,b,c\}$ and $W_2,\cdots,W_{t}$ are positive words of length $r$. 
By assumption $t>1$, there exists an index $i\! \in\! \{2,...,t\}$ so that we decompose the sequence into two steps $u_1 X \simeq u_i W_i \simeq u_{t+1}Y$, where each step satsifies the induction hypothesis. 

If there exists $i$ such that $u_i$ is equal either to $u_1$ or $u_{t+1}$, then by induction hypothesis, $W_i$ is equivalent either to $X$ or to $Y$. Then, again, applying the induction hypothesis to the remaining step,
we obtain the statement for the $u_1X\simeq u_{t+1}Y$. Thus, we assume from now on $u_i\not= u_1,u_{t+1}$ for $1<i\le t$.

Suppose $u_1\!=\!u_{t+1}$. If there exists $i$ such that $\{u_1\!=\!u_{t+1},u_i\}\!\not=\!\{b,c\}$, then each of the equvalence says the existence of $\alpha,\beta\!\in\!\{a,b,c\}$ and words $Z_1,Z_2$ such that $X\simeq \!\alpha\! Z_1$, $W_i\!\simeq\! \beta Z_1\!\simeq\! \beta Z_2$ and $Y\!\simeq\! \alpha Z_2$. Applying the induction hypothesis for $r$ to $\beta Z_1 \!\simeq\! \beta Z_2$, we get $Z_1 \!\simeq\! Z_2$ and, hence, we obtained the statement $X\!\simeq\! \alpha Z_1\!\simeq\! \alpha Z_2\!\simeq\! Y$.
Thus, we exclude these cases from our considerations. Next, we consider the case $\{u_1\!=\!u_{t+1},u_i\}\!=\!\{b,c\}$. However, due to the above consideration, we have only the case $u_2\!=\!u_3\!=\!\cdots\!=\!u_t$.  Then, by induction hypothesis, we have $W_2\!\simeq\!\cdots\!\simeq\!W_t$. On the other hand, since the equivalences $u_1X\!\simeq\!u_2W_2$ and $u_{t+1}Y\!\simeq\!u_{t}W_t$ are the elementary trasformations at the tops of the words, there exists again  $\alpha,\beta\!\in\!\{a,b,c\}$ and words $Z_1,Z_2$ with the similar descriptions as above hold, implying again $X\!\simeq\! Y$.

To complete the proof, we have to examine three more cases $(u_1, u_2, u_{3}) = (a, b, c), (a, c, b)$ and $(b,a,c)$ for $t=2$, where we shall put $W:=W_2$.

\smallskip
 $(\mathrm{I})$ \quad Case $(a, b, c)$.  We have $aX\simeq bW\simeq cY$. \\
Since the equivalences are single elementary transformations, there exists words $Z_1$ and $Z_2$ such that 
$X\simeq bZ_1$, $W\simeq cZ_1\simeq baZ_2$ and $Y\simeq bbZ_2$.
Applying the induction hypothesis for $r$ to the two equivalent expressions of $W$, we see that there exists $k$ and a word $Z_3$ such that $0\le k<r-2$ such that $Z_1\simeq a^kbbZ_3$ and $aZ_2\simeq c^kbaZ_3$. We can apply $k$-times the induction hypothesis to the last two equivalent expressions and we see that there exists a word $Z_4$ such that $Z_2\simeq c^kZ_4$ and $baZ_3\simeq aZ_4$. 
Applying again the induction hypothesis to the last equivalence relation, there exists a word $Z_5$ such that $Z_4\simeq bZ_5$ and $aZ_3\simeq cZ_5$. Once again applying the induction hypothesis to the last equivalence relation, we finally obtain $Z_3\simeq cZ_6$ and $Z_5\simeq aZ_6$ for a word $Z_6$. Reversing the procedure, obtain the descriptions:
\[\begin{array} {lll}
X& \simeq & bZ_1\simeq ba^kbbZ_3\simeq ba^kbbcZ_6 , \\
Y& \simeq & bbZ_2\simeq bb c^kZ_4\simeq bbc^k bZ_5\simeq bbc^kbaZ_6.
\end{array}
\]

By using the relations of ${\mathrm{M_{B_{ii}}}}$, we can  show $ba^kbbc \simeq cbbc^kb$ and  
$ bbc^kba \simeq abbc^kb$. So, we conclude that $X \simeq cZ, Y \simeq aZ$ for $Z\simeq bbc^kbZ_6$.

\smallskip
$(\mathrm{II})$ \quad Case $(a, c,b)$. We have $aX\simeq cW\simeq bY$.\\
Since the equivalences are single elementary transformations, there exists words $Z_1$ and $Z_2$ such that 
$X\simeq cZ_1$, $W\simeq aZ_1\simeq bbZ_2$ and $Y\simeq baZ_2$.
Applying the induction hypothesis for $r$ to the two equivalent expressions of $W$, we see that there exists a word $Z_3$ such that $Z_1\simeq bZ_3$ and $bZ_2\simeq cZ_3$. Again applying the induction hypothesis to the last two equivalent expressions, we see that there exists an integer $k$ with $0\le k<r-3$ and a word $Z_4$ such that $Z_2\simeq c^kbaZ_4$ and $Z_3\simeq a^kbbZ_4$.  
Reversing the procedure, obtain the descriptions:
\[
X\ \simeq \ cZ_1\ \simeq\ cbZ_3\ \simeq \ cba^kbbZ_4 \quad \text{ and }\quad 
Y\ \simeq \ baZ_2\ \simeq\ ba c^kbaZ_4.
\]
It is not hard to show the equivalences $cba^kbb\simeq bbac^kb$ and 
$bac^kba\simeq cbac^kb$. Thus, we obtain $X \simeq bZ, Y \simeq cZ$ for  $Z:= bac^kbZ_4$.

\smallskip
$(\mathrm{III})$ \quad Case $(b,a, c)$. We have $bX\simeq aW\simeq cY$.\\
By induction hypothesis, there exists words $Z_1$ and $Z_2$ such that 
$X\simeq cZ_1$, $W\simeq bZ_1\simeq cZ_2$ and $Y\simeq aZ_2$.
Applying the induction hypothesis for $r$ to the two equivalent expressions of $W$, we see that here exists $k$ and a word $Z_3$ such that $0\le k<r-2$ such that $Z_1\simeq c^kbaZ_3$ and $Z_2\simeq a^kbbZ_3$. Thus, we obtain the descriptions:
\[
X\ \simeq \ cZ_1\ \simeq\ cc^kbaZ_3  \quad \text{ and }\quad
Y\ \simeq \  aZ_2\ \simeq\  a a^kbbZ_3.
\]
This\! is\! the\! conclusion\! in\! Proposition\! (iv)\! with\! $0\!\le \! k\!\!+\!\!1\!\!<\!\!r\!-\!1$,\! which\! we\! looked\! for.\

This completes the proof of Proposition.
\end{proof}
\vspace{-0.3cm}
This completes the proof of Theorem 4.
\end{proof}

\section{$2\times 2$-matrix representation of the group $G_X$}

 We construct non-abelian representations of the groups $G_\mathrm{B_{ii}},G_\mathrm{B_{vi}}, G_\mathrm{H_{ii}}, 
G_\mathrm{H_{iii}}$.

\bigskip
\noindent
{\bf Theorem 5.} {\it 
For each type $X\in\{\mathrm{B_{ii}}, \mathrm{B_{vi}}, \mathrm{H_{ii}}, 
\mathrm{H_{iii}}\}$, consider matrices $A,B,C$ $\in \mathrm{GL}(2,\C)$
listed below.  Then we have the following i) and ii).

i) The correspondence $a\mapsto A, b\mapsto B, c\mapsto C$ induces representations 

\centerline
{$\rho\ : \ G_X\longrightarrow \mathrm{GL}(2,\C)$.}


ii) The image $\rho(G_X)$ is not an abelian group if $l^2\not=1$.

\medskip
\noindent
1. Type $\mathrm{B_{ii}}$:\quad
{\small     
$
 A=u\left( \begin{array}{cc}
1&l^2\\0&1
\end{array}
\right)
, \quad 
B=v\left( \begin{array}{cc}
l&0\\0&l^{-1}
\end{array}
\right), \quad 
C=u\left( \begin{array}{cc}
1&1\\0&1
\end{array}
\right),
$
}

\noindent
where $ l^6=1$ and  $u,v \in \C^\times$.

\medskip
\noindent
2. Type $\mathrm{B_{vi}}$: \quad 
{\small
$
A=u\left( \begin{array}{cc}
l&0\\0&l^{-1}
\end{array}
\right), 
\quad 
B=u\left( \begin{array}{cc}
a&b\\c&d
\end{array}
\right)
, 
\quad 
C=u\left( \begin{array}{cc}
p&q\\r&s
\end{array}
\right),
$
}

\noindent
where $l^{10}=1$ $(l^2\not=1)$ and $u \in \C^\times$

{\small
\[
a=-\frac{1}{l(l^2-1)}, \quad b c=\frac{-l^4+l^2-1}{(1-l^2)^2}, 
\quad d=\frac{l^3}{l^2-1}
\]
\[
p=-l^4a, \quad q=-\frac{b}{l^4}, 
\quad r=-l^4c, \quad s=-\frac{d}{l^4}
\]
}

\medskip
\noindent
3. Type $\mathrm{H_{ii}}$: \quad
{\small
$
A=
u\left( \begin{array}{cc}
l&0\\0&l^{-1}
\end{array}
\right)
, 
\quad
B=
u\left( \begin{array}{cc}
a&b\\c&d
\end{array}
\right)
, 
\quad
C=
u\left( \begin{array}{cc}
p&q\\r&s
\end{array}
\right)
$
}

\noindent
where $u \in \C^\times$ and one of the following two cases holds.

i) $l^2+l+1=0$ and $3p^2+3p+2=0$
{\small
\[
a=\frac{l-1}{3},\   d=\frac{-l-2}{3},\  bc=-\frac{2}{3},\ 
q=\frac{-b(l+2)}{3p},\ r=\frac{p(1-l)}{3b} ,\ s=\frac{2}{3p} .
\]
}

ii) $l^2-l+1=0$ and $3p^2-3p+2=0$.
{\small
\[
a= \frac{l+1}{3},\   d=\frac{-l+2}{3} ,  bc=-\frac{2}{3} ,\
q= \frac{b(-l+2)}{3p},\ r=\frac{-p(l+1)}{3b} ,\ s=\frac{2}{3p}.
\]
}

\medskip
\noindent
4. Type $\mathrm{H_{iii}}$:\quad 
{\small
$
A=
u\left( \begin{array}{cc} 
l&0\\0&l^{-1} 
\end{array} 
\right) 
, 
\quad
B=
u\left( \begin{array}{cc} 
a&b\\c&d 
\end{array} 
\right) 
, 
\quad
C=
u\left( \begin{array}{cc} 
p&q\\r&s 
\end{array} 
\right) ,
$
}

\noindent
where $l^{10}=1$ $(l^2\not=1)$ and $u \in \C^\times$

{\small
\[
a=-\frac{1}{l(l^2-1)}, \quad bc=\frac{-l^4+l^2-1}{(l^2-1)^2}, 
\quad d=\frac{l^3}{l^2-1} 
\]
\[
p=a, \quad q=\frac{b}{l^4}, 
\quad r=l^4c, \quad s=d
\vspace{-0.3cm}\]
}
}
\begin{proof}
It is sufficient to prove only for the case $u=v=1$.

We present the matrices $A,B$ and $C$ by the indeterminates $a,b,c,d,p,q,r,s$ and $l$ as in Theorem, and then solve the polynomial equation on them defined by the relations listed in Theorem 1. It is unnecessary to check all relations, since some relations are included in the list because of the cancellation condition, whereas $\mathrm{GL}(2,\C)$ is a group where the cancellation condition is automatically satisfied. However, as we shall see, it is sometimes convenient to take these ``superfluous'' relations in account. Detailed calculations are left to the reader.

\medskip
1.Type $\mathrm{B_{ii}}$: 
We need to show $CBB=BBA, BC=AB, AC=CA$, whose verifications are left to the reader. We have $\det(A)\!=\! \det(C)\!=\!u^2\!\not=\!0, \det(B)\!=\!v^2\!\not=\!0$.
Since 
{\small 
$ABA^{-1}B^{-1}=\left( \begin{array}{cc} 
1&l^2(1-l^2)\\0&1 
\end{array} 
\right)$
}
and
{\small 
$BCB^{-1}C^{-1}=\left( \begin{array}{cc} 
1&l^2-1\\0&1 
\end{array} 
\right)$
},
$\rho(G_{\mathrm{B_{ii}}})$ is abelian if and only if $l^2=1$.

2. Type $\mathrm{B_{vi}}$: 
We need to show $ABA=BAB, ACA=BAC, ACB=CAC$. Actually, solving the (1,1) entry of the equation 
$ABA=BAB$, $\mathrm{tr}( A)=\mathrm{tr}( B)$ and $\det(B)=1$, w obtain the expressions for $a,b,c,d$. 
Then, using the relation $C=ABA^{-2}B$, we obtain the expressions for $p,q,r,s$. 
Further more, comparing the $(1,1)$-entry of $A^5=B^5$, we get $l^8+l^6+l^4+l^2+1=0$.

3. Type $\mathrm{H_{ii}}$:
We need to show $ABAB\!=\!BABA, ACA\!=\!BAC, ACB\!=\!CAC$.

{\small
\[
ABAB=
\left( \begin{array}{cc} 
bc+a^2l^2&bd+abl^2\\ac+cd/l^2&bc+d^2/l^2 
\end{array} 
\right) 
, BABA= 
\left( \begin{array}{cc} 
bc+a^2l^2&ab+bd/l^2\\cd+acl^2&bc+d^2/l^2 
\end{array}
\right)
\]
}
By these calculations, we have $d+al^2=0$. 
By $TrA=TrB=TrC$ and $detA=detB=detC$, 
we have $a+d=l+l^{-1}=p+s$, $ad-bc=ps-qr=1$.
{\footnotesize
\[
a=\frac{l^2+1}{l(1-l^2)}, 
d=\frac{l(l^2+1)}{l^2-1}, 
bc=\frac{-2(l^4+1)}{(l^2-1)^2} 
\]
}
{\footnotesize
\[
ACA= 
\left( \begin{array}{cc}
l^2p&q\\r&s/l^2
\end{array}
\right) 
, BAC= 
\left( \begin{array}{cc} 
alp+br/l&alq+bs/l\\clp+dr/l&clq+ds/l 
\end{array} 
\right) 
\]
}
{\footnotesize
\[
q=\frac{b}{lp(1-l^2)}, 
r=\frac{l(l^4+1)p}{b(l^2-1)}, 
s=\frac{-2l^2}{(1-l^2)^2p} 
\]
}
{\footnotesize
\[
ACB=CAC 
\Leftrightarrow 
(1-l+l^2)=0 {\rm \ and \ } 3p^2-3p+2=0,
{\rm \ or \ }, (1+l+l^2)=0 {\rm \ and \ } 3p^2+3p+2=0
\]
}
We calculate each cases separately and obtain the result.

4. Type $\mathrm{H_{iii}}$:
We need to show $ABA\!=\!BAB, CBA\!=\!ACB, BCA\!=\!CBC$. 
As in case of $\mathrm{B_{vi}}$, already the first relation 
$ABA=BAB$ (in particular $\mathrm{tr}( A)=\mathrm{tr}(B)$ and $\det(B)=1$) implies the expressions for $a,b,c,d$. 
Using further the relation $ACA=CAC$, we obtain $a=p$, $d=s$ and $bc=qr$.  Then applying the relation $A^2C=BA^2$, we get $q=l^4b$ and $r=l^{-4}c$. Further, using the relation $CA^3=A^3B$, we obtain $l^{10}=1$.
\end{proof}

\noindent
{\it Corollary.  For $X\!\in\!\{\mathrm{B_{ii}}$, $\mathrm{B_{vi}}$, $\mathrm{H_{ii}}, \mathrm{H_{iii}}\}$, $\sigma(\mathcal{QZ}(G_X^+))$ consists only of the identity.}

\medskip
\noindent
{\it Sketch of Proof.} For $\sigma\in \mathfrak{S}(L)$, we consider a matrix $X\in \mathrm{Mat}(2,\C)$ satisfying the equations: $AX=X\sigma(A),\ BX=X\sigma(B),\ CX=X\sigma(C)$. If  $\sigma=1$, then the solutions are $constant\times \mathrm{id}$. If $\sigma\not=1$, then $X=0$.

\medskip
\noindent
{\it Acknowledgement.}\! We thank\! to David Bessis for pointing out that 
the group\!  $G_{\mathrm{H_{vi}}}$\!\! is an infinite cyclic group, and to 
Masaki Kashiwara for interesting discussions.

\end{document}